\UseRawInputEncoding
\documentclass[a4paper]{article}

\usepackage[latin1]{inputenc}    
\usepackage[T1]{fontenc}       
\usepackage{lmodern}           

\usepackage{amsmath, amssymb, amsfonts}

\usepackage{graphicx}          
\usepackage{subcaption}        

\usepackage{paralist}          
\usepackage{geometry}          
\usepackage{setspace}          
\usepackage{authblk}           

\usepackage[
    colorlinks=true,
    linkcolor=blue,
    citecolor=red,
    urlcolor=cyan,
    unicode=true
]{hyperref}

\newtheorem{example}{Example}[section]
\newtheorem{remark}{Remark}[section]
\numberwithin{equation}{section}    
\numberwithin{figure}{section}     
\numberwithin{table}{section}      

\title{A simple, high-order and compact WENO limiter based on control volume for spectral volume method}

\author[1,2]{Na Liu}
\author[3]{Jianxian Qiu\thanks{Corresponding author, Email: jxqiu@xmu.edu.cn}}
\affil[1]{CAEP Software Center for High Performance Numerical Simulation, Beijing 100088, PR China}
\affil[2]{Institute of Applied Physics and Computational Mathematics, Beijing 100088, PR China}
\affil[3]{School of Mathematical Sciences and Fujian Provincial Key Laboratory of Mathematical Modeling and\\
          High-Performance Scientific Computing, Xiamen University, Xiamen, Fujian 361005, PR China}

\date{}

\begin{document}
\maketitle






\begin{abstract}
The spectral volume(SV) method constructs a high-order polynomial for SV based on the average value of control volume(CV),
but for discontinuous problems, a limiter is required to mitigate oscillations. This paper presents a
novel CV-based high-resolution limiter to effectively suppress oscillations and maintain CV
 resolution. Drawing inspiration from the SWENO method \cite{zhu2019swenoRKDG}, we utilize a nonlinear weighting approach to reconstruct a
    novel high-order polynomial for the target control volume by combining the high-order reconstructed
    polynomial and linear polynomials which are reconstructed by the cell average of the target CV and its neighboring CVs. The new high-order polynomial breaks the continuity in the SV, thus the utilization of numerical flux at the boundaries of troubled CVs and the SV boundaries.
     However, at other boundaries of CVs where physical quantities remain continuous, direct calculation of flux
      based on the values of physical variables is feasible. The limiter is simple, as it only requires
      several linear polynomials in the limitation process. Moreover, it still maintains the compactness
       of the SV method and preserves the resolution of CV. Numerical results for one- and two-dimensional
        scalar and system of conservation laws verified that the high-order property of the CV-SWENO limiter
         is effective in solving both smooth and strongly discontinuity problems.
\end{abstract}

%
%
%

\section{Introduction}
In recent decades, high-order methods and their applications in refined simulation have become research hotspots. Notable achievements have been realized
alongside the progressive enhancement of computational floating-point capabilities.
High-order methods possess significant advantages over low-order methods in controlling numerical errors when applied to problems such as turbulent flow that require fine simulation \cite{wang2013High}. For discontinuous flows, such as strong shock waves and contact discontinuities, high-order methods have
demonstrated superior performance \cite{shu1998enoweno,cockburn1998dg,zhu2019swenoRKDG,wang2002SV1,liu2014hgks,liu2017svNOK}. Over the past thirty years, various high-order methods have been developed, including essentially non-oscillatory (ENO) method \cite{Harten1987eno}, weighted essentially non-oscillatory (WENO) finite volume method \cite{Osher1994weno,Ren2012kexactWENO,shu1998enoweno}, discontinuous Galerkin (DG) method \cite{cockburn2000dg,cockburn1989tvb2,cockburn1998dg,zhu2019swenoRKDG,zhong2013swenoDG,Balsara2007hwenodg,zhu2013dg,qiu2005dg,zhu2008dg}, spectral and spectral-element method \cite{1999Stable},
spectral volume (SV) method \cite{wang2002SV1,wang2006SV5,wang2006SV6,liu2017svNOK,liu2017svLF,Yousefi2019CSV,Yousefi2022CSV}, spectral difference (SD) method \cite{wang2006SD,wang2007SD}, and correction procedure via reconstruction (CPR) method \cite{Huynh2007CPR,wang2009CPR,Du2015swenoCPR}. For the current research progress in high-order methods in CFD, readers can refer to \cite{2005High,wang2007High,wang2013High}.

The SV method was first proposed by Wang in 2002 \cite{wang2002SV1} for solving conservation laws on unstructured grids. By introducing  the concept of "spectral volume", the computational domain is divided into multiple cells referred to as SVs, with each SV further partitioned into subcells known as CVs.
The method employs the averages of CVs within each SV to reconstruct an approximate solution of the SV, with each reconstruction being independent of the others. Thus, the SV method is a compact high-order method. When the solution is "smooth" in SV schemes, numerical fluxes (exact or approximate Riemann solutions) are only used at the SV interface, whereas at the CV interface, continuous fluxes are used directly. Similar to the DG method, the SV method can flexibly handle complex geometries with h-p adaptivity \cite{wang2009SV8} and achieve high parallel implementation efficiency. More importantly, it provides higher efficiency in CPU time, particularly for higher-order accuracy schemes \cite{Abeele2007SV2}, and improved resolution, especially in cases of discontinuous solutions \cite{wang2006SV5}. Based on the analysis of the dissipation and dispersion of the SV method, Van den Abeele et al. \cite{Abeele2007SV1,Abeele2008SV3} discovered that the stability of this method is contingent upon cell division, with certain cell divisions potentially leading to marginal instability. To address this issue, Harris et al. \cite{Harris2009SV} studied optimized cell division in two-dimensional grids. Currently, the SV method has been successfully applied to solve various conservation laws \cite{wang2002SV1,wang2006SV5,Yousefi2019CSV}, shallow-water wave equations \cite{Yousefi2022CSV}, multi-medium problems \cite{liu2017svNOK,liu2017svLF}, and viscous fluid dynamics \cite{wang2006SV6}.

Solutions of nonlinear fluid equations often contain discontinuities. To avoid numerical oscillations that arise from high-order reconstruction, nonlinear limiters are commonly used to detect and mitigate spurious oscillations in the vicinity of the discontinuities. There are various types of limiters, including total variation bounded (TVB) limiters \cite{cockburn1989tvb2}, moment limiters \cite{Biswas1994,Burbeau2001}, multi-dimensional optimal order detection (MOOD) limiters \cite{Clain2011mood,Clain2012mood}, ENO/WENO limiters \cite{Harten1987eno,Osher1994weno}, Hermite-WENO limiters \cite{qiu2003hwenodg,qiu2009hwenodg3,Balsara2007hwenodg}. In \cite{wang2002SV1,wang2006SV5}, the SV method utilized TVB and total variation diminishing (TVD) limiters.
Although these limiters effectively control spurious numerical oscillations near discontinuities, they can often compromise accuracy in smooth regions of the solution. An ideal limiter should maintain high-order accuracy in smooth regions, control spurious numerical oscillations near discontinuities, and preserve high-resolution properties. High-resolution limiters have been developed for DG methods, such as WENO-based limiters designed by Qiu et al. \cite{qiu2005dg} and Zhu et al. \cite{zhu2008dg}. However, such limiters require reconstruction stencils match the width of the WENO method. This requirement involves not only neighboring cells but also the neighbors of those neighbors, thereby reducing the compactness of DG methods. To enhance the compactness of the stencils, Qiu et al. \cite{qiu2003hwenodg,qiu2009hwenodg3} developed a Hermite-WENO limiter. However, high-order methods utilizing this limiter still require information from second-neighbor cells. Later, Zhong and Shu \cite{zhong2013swenoDG} proposed a simplified WENO limiter, which constructs a nonlinear combination of the target cell and its neighboring cells based on primitive reconstruction. The stencil is compact, only including the target cell and its immediate neighboring cells. Recently, Zhu et al. \cite{zhu2016sweno,zhu2017sweno2,zhu2017sweno3,zhu2018sweno4} proposed a new simplified WENO method that uses only a high-order polynomial and a few linear polynomials for reconstruction, thereby enhancing both simplicity and efficiency compared to the classical WENO method. Zhu et al. \cite{zhu2019swenoRKDG} further applied this limiter to DG methods, achieving promising numerical results.

Both the SV and DG methods are compact schemes that rely on the idea of constraining the high-order polynomial reconstructions within the primitive cell. All limiters utilized in the DG method can be extended to the SV method, as demonstrated by Yousefi et al. \cite{Yousefi2019CSV}, who successfully adapted Zhong and Shu's simple WENO limiter \cite{zhong2013swenoDG} for use in the SV method. Similar to the DG method, the limiter operates on the entire SV and ultimatedly generates a constrained reconstructed polynomial for the SV. This was the first time that a high-resolution limiter was incorporated into the SV method for solving conservation laws and shallow water equations \cite{Yousefi2022CSV}. Du et al. \cite{Du2015swenoCPR} extended the limiter to the CPR method. Unlike the DG method, the SV method stores CV averages and achieves the resolution of CVs. According to this property, limiters can be applied to the SV elements or the CVs, referred to as SV-wise limiters and CV-wise limiters, respectively. Reference \cite{wang2002SV1} compared the resolutions of SV-wise and CV-wise TVB and TVD limiters, indicating that CV-wise limiters exhibit higher resolution than those SV-wise limiters.

In this paper, we develop a CV-wise high-resolution limiter using the WENO methodology for the SV method. This approach is characterized by its simplicity, efficiency, and compactness. The limiter utilizes the simple WENO limiter framework developed by Zhu, Qiu et al., achieving limitation through the combination of high-order polynomial and several linear polynomials. This approach ensures that the limitation is applied to CVs, thereby preserving the resolution of the method. This paper is organized as follows. Section 2  briefly reviews the SV method, while Section 3 presents the newly developed WENO limiter. Section 4 provides extensive one- and two- dimensional numerical results to validate the accuracy and efficiency of the proposed limiter. Conclusions are drawn in Section 5, and the last section is the appendix.

\section{Review of the spectral volume method}
Considering hyperbolic conservation laws,
\begin{align}\label{EQ:sec2hyperbolic}
u_t+\nabla \cdot f(u)=0, \\
u(\bold x,0)=u_0(\bold x),
\end{align}
where $u$ is either a scalar or a vector, $f(u)$ is the corresponding flux, $\bold x$ represents the spatial coordinate in region $\Omega$, and $t$ denotes time with the interval $[0,T]$. In this section, we briefly review the aforementioned SV method for solving hyperbolic conservation laws. For detailed implementation of the scheme, please refer to \cite{wang2002SV1,liu2017svNOK,Yousefi2019CSV}.
\subsection{Spacial discretization}
Without loss of generality, this section focuses on the reconstruction of SV on two-dimensional rectangular grids. The computational domain $\Omega$ is partitioned into $N=N_x \times N_y$ non-overlapping sub-intervals denoted as $\{S_{i,j}\}_{i=1\ldots N_x, j=1\ldots N_y}$, where each $S_{i,j}=[x_{i-\frac12},x_{i+\frac12}]\times [y_{j-\frac12},y_{j+\frac12}]$ represents a spectral volume. The length of $S_{i,j}$ in the $x$ and $y$ directions
are given by $h^x_{i}=x_{i+\frac12}-x_{i-\frac12}$  and $h^y_{j}=y_{j+\frac12}-y_{j-\frac12}$, respectively. Assuming a $\bold{k}$-th order accurate scheme is employed, each $S_{i,j}$ is further subdivided into $\bold{k}^2$ control volumes, i.e.,
\begin{equation*}
C_{i,m,j,n}=(x_{i,m-\frac12},x_{i,m+\frac12})\times (y_{j,n-\frac12},y_{j,n+\frac12}),m,n=1,\ldots,\bold{k},
\end{equation*}
where $x_{i,\frac12}=x_{i-\frac12}$, $x_{i,\bold{k}+\frac12}=x_{i+\frac12}$, $y_{j,\frac12}=y_{j-\frac12}$, $y_{j,\bold{k}+\frac12}=y_{j+\frac12}$. Let $h^x_{i,m}=x_{i,m+\frac12}-x_{i,m-\frac12}$ and $h^y_{j,n}=y_{j,n+\frac12}-y_{j,n-\frac12}$ denote the lengths of $C_{i,m,j,n}$ in the $x$ and $y$ directions, respectively.
The time domain $[0, T]$ is partitioned as follows: $\{t_{s+1} = t_s + \Delta t_s,t_0 = 0,  s \geq 0\}$ ,
where the time step $\Delta t_s$ is determined by the CFL condition.
The cell average of $C_{i,m,j,n}$ is defined as
\begin{equation}\label{EQ:sec2cellaverage2d}
\overline{u}_{i,m,j,n}^s=\frac{1}{V}\int_{C_{i,m,j,n}} u(x,y,t_s) dV, m,n=1,\ldots,\bold{k},
\end{equation}
where $V$ denotes the volume of $C_{i,m,j,n}$. Integrating the governing
 Eq. \eqref{EQ:sec2hyperbolic} over $C_{i,m,j,n}$ yields a semi-discrete scheme:
\begin{equation}\label{EQ:sec2semischeme2d}
\frac{d\overline{u}_{i,m,j,n}^s}{dt}=-\frac{1}{V}(\int_{y_{j,n-\frac12}}^{y_{j,n+\frac12}}(\widehat{F}_{i,m+\frac{1}{2}}(y)-\widehat{F}_{i,m-\frac{1}{2}}(y))dy+
\int_{x_{i,m-\frac12}}^{x_{i,m+\frac12}}(\widehat{G}_{j,n+\frac{1}{2}}(x)-\widehat{G}_{j,n-\frac{1}{2}}(x))dx),
\end{equation}
where $\widehat{F}$ and $\widehat{G}$ represent the numerical fluxes of the two-dimensional flux $(f,g)$ in the x and y directions, respectively. The integral can be approximated using Gaussian quadrature. The subsequent task is to achieve the two-dimensional SV reconstruction of $u(x,y)$ and to compute the numerical flux at the Gaussian quadrature points based on this reconstruction.

\textbf{Spectral volume reconstruction}
According to the average of $\bold{k}^2$ control volumes, a two-dimensional bilinear polynomial, of degree at most $(\bold{k}-1)$, can be reconstructed in each $S_{i,j}$,
\begin{equation}\label{EQ:sec2recon2d}
u_h^s(x,y):=p_{i,j}^s(x,y)
=u_{i,j}(x,y,t_s)+O(h^\bold{k}), (x,y)\in S_{i,j}.
\end{equation}
For simplicity, the superscript $s$ is omitted in the following.
$p_{i,j}(x,y)$, a $\bold{k}$-th order accurate approximation to $u(x,y)$ in $S_{i,j}$,
can be expressed as
\begin{equation}\label{EQ:sec2recon2d1}
p_{i,j}(x,y)=\sum_{l=0}^{\bold{k}-1}\sum_{r=0}^{\bold{k}-1}\overline{W}_{i,j}^{(l,r)}\phi^{(l)}(x)\phi^{(r)}(y),
\end{equation}
where $\phi^{(l)}(x)$, $l=0,1,\ldots,\bold{k}-1$, denotes a set of basis functions
that belong to the space of polynomials of degree at most $\bold{k}-1$.
The cell-wise average derivatives $\overline{W}_{i,j}^{(l,r)}$, $l,r=0,1,\ldots,\bold{k}-1$ can be obtained by solving the following linear algebraic equation system established using the CV averages in $S_{i,j}$,
\begin{align}\label{EQ:sec2Wcoeff2}
\sum_{l=0}^{\bold{k}-1}\sum_{r=0}^{\bold{k}-1} \overline{W}_{i,j}^{(l,r)}\int_{CV}\phi^{(l)}(x)\phi^{(r)}(y)dxdy=
\int_{CV} p_{i,j}(x,y)dxdy=V \overline{u}_{i,m,j,n}, m,n=1,\ldots \bold{k}.
\end{align}

\begin{remark}\label{remark2}
The reconstruction function $u_h(x,y)$ for the entire computational domain has the following properties:\\
1. The reconstruction polynomial function is continuous within each SV.\\
2. The reconstruction polynomial function exhibits discontinuity at the boundaries of the SVs.
\end{remark}

\textbf{Numerical flux}
The numerical flux at the SV boundary can be computed using the Riemann solver,
\begin{align}\label{EQ:sec2fluxd}
\widehat{F}_{i,\frac{1}{2}}(u(x_{i,\frac{1}{2}},y_{j,\beta}))=f_{Riemann}(u(x_{i-1,\bold{k}+\frac{1}{2}},y_{j,\beta}),u(x_{i,\frac{1}{2}},y_{j,\beta}),\\
\widehat{F}_{i,\bold{k}+\frac{1}{2}}(u(x_{i,\bold{k}+\frac{1}{2}},y_{j,\beta}))=f_{Riemann}(u(x_{i,\bold{k}+\frac{1}{2}},y_{j,\beta}),u(x_{i+1,\frac{1}{2}},y_{j,\beta})),
\end{align}
where $x_{i,\alpha}$ and $y_{j,\beta}$ are the Gaussian integration points on $(x_{i,m-\frac12},x_{i,m+\frac12})$ and $(y_{j,n-\frac12},y_{j,n+\frac12})$, respectively.
Commonly-used numerical flux calculation methods include Lax-Friedrichs, HLLC, Roe, NOK \cite{liu2017svNOK}, BGK \cite{liu2014hgks}, etc. For simplicity, this paper adopts the Lax-Friedrichs numerical flux.

Due to Remark \ref{remark2}, $u_h(x,y)$ has a continuous internal state at the CV boundary inside $S_{i,j}$.
The numerical flux satisfies locally Lipschitz condition and is consistent with the flux function $f(u)$,
\begin{align}
f_{Riemann}(u,u)=f(u).
\end{align}
Therefore, the numerical flux at the internal CV boundaries in SV degenerates to the analytical flux,
\begin{align}\label{EQ:sec2fluxc2d}
\widehat{F}_{i,m+\frac{1}{2}}(u(x_{i,m+\frac{1}{2}},y_{j,\beta}))=f(u(x_{i,m+\frac{1}{2}},y_{j,\beta})),m=1,...,\bold{k}-1.
\end{align}
The numerical flux in the y direction follows a similar formulation as described above.

\subsection{Time discretization}
The semi-discrete schemes Eq. \eqref{EQ:sec2semischeme2d} can be written into the following ordinary differential equation:
\begin{align}\label{EQ:sec2time}
\frac{du}{dt}=L(u),
\end{align}
 where $L(u)$ is the spatial discrete operator. The $\bold{k}$-th order SV method is discretized in the time direction by employing the $\bold{k}$-th order Runge-Kutta method \cite{shuosher1988,shuosher1989,cockburn1989tvb2,zhu2019swenoRKDG}:
\begin{align}\label{EQ:sec2RK}
u^{(i)}=\sum_{l=0}^{i-1}(\alpha_{i,l}u^{(l)}+\beta_{i,l}\Delta tL(u^{(l)})),i=1,...,\bold{k},\\
u^{(0)}=u^s,u^{(s+1)}=u^{(\bold{k}+1)}.
\end{align}
The values of coefficients $\alpha_{i,l}, \beta_{i,l}$ are listed in Table \ref{table:sec2RK1}.
\begin{table}[htbp]
 \begin{center}   \caption{Parameters of Runge-Kutta scheme.}
 \label{table:sec2RK1}
   \begin{tabular}{c|ccccc|ccccc}
     \hline
order &         &$\alpha_{i,l}$ &       &        &           &    &$\beta_{i,l}$ &   &     &
     \\
       \hline
       & 1          &           &               &           & 1         &    &      &   \\
     2 & $\frac12$  & $\frac12$ &               &           & 0         & $\frac12$  &      &   \\
     \hline
       & 1          &           &               &           & 1         &    &      &   \\
     3 & $\frac34$  & $\frac14$ &               &           & 0         &  $\frac14$  &      &   \\
       & $\frac13$  & 0         & $\frac23$     &           & 0         &  0  &$\frac23$ &   \\
     \hline
       & 1          &           &               &           &$\frac12$  &     &   &   \\
     4 & 1          & 0         &               &           & 0         &  $\frac12$   &   &   \\
       & 1          & 0         & 0             &           & 0         &  0           & 1 &   \\
       & $-\frac13$ & $\frac13$ & $\frac23$     & $\frac13$ & 0         &  0           & 0 &  $\frac16$ \\
      \hline
       & 1          &           &       &        &           &$\frac15$  &     &   &   &\\
       & 1          & 0         &       &        &           & 0         &  $\frac14$   &   &  & \\
     5 & 1          & 0         & 0     &        &           & 0         &  0           & $\frac13$ &  & \\
       & 1          & 0         & 0     &        & 0         & 0         &  0           & 0 &  $\frac12$ &\\
       & 1          & 0         & 0     &        & 0         & 0         &  0           & 0 &  0 &1\\
      \hline
   \end{tabular}
 \end{center}
 \end{table}

\section{A new CV-wise WENO limiter for spectral volume method}
The polynomials reconstructed using Eq. \eqref{EQ:sec2recon2d1} may introduce numerical oscillations.
To mitigate these oscillations, various limiters have been employed in the reconstruction of SV schemes \cite{wang2002SV1,wang2006SV5}, including control volume total variation diminishing (CVTVD), control volume total variation bounded (CVTVB), spectral volume total variation diminishing (SVTVD), and spectral volume total variation bounded (SVTVB) limiters. The CV-wise limiters are found to outperform those SV-wise counterparts in suppressing numerical oscillations while maintaining high resolution. In numerical simulations, our objective is to achieve reliable and stable solutions with minimal loss of resolution, avoiding overly restrictive limiting. Hence, the idea of "tagging" and "limiting" \cite{qiu2005dg,zhu2019swenoRKDG} is employed. Specifically, elements that may cause oscillations are first identified, and then a high-order precision limiter, specifically the SWENO limiter, is applied to these identified elements. This approach relies only on the SV element they belong to and the neighboring SV, thereby preserving the compactness inherent in SV method.

\subsection{One-dimensional case}
First of all, the modified TVB minmod function \cite{cockburn1989tvb2} is used to identify "troubled cells". Given that the cell average
and function values at the CV boundaries are known for each  $C_{i,m}$, we denote:
\begin{eqnarray}
\delta_+u_{i,m}=u_{i,m+\frac12}-\overline{u}_{i,m},\delta_-u_{i,m}=\overline{u}_{i,m}-u_{i,m-\frac12},\\
\Delta_+\overline{u}_{i,m}=\overline{u}_{i,m+1}-\overline{u}_{i,m},
\Delta_-\overline{u}_{i,m}=\overline{u}_{i,m}-\overline{u}_{i,m-1}.
\end{eqnarray}
We further denote
\begin{eqnarray}
\delta_+u_{i,m}^{(mod)}=\widetilde{m}(\delta_+u_{i,m},\Delta_+\overline{u}_{i,m},\Delta_-\overline{u}_{i,m}),\\
\delta_-u_{i,m}^{(mod)}=\widetilde{m}(\delta_-u_{i,m},\Delta_+\overline{u}_{i,m},\Delta_-\overline{u}_{i,m}),
\end{eqnarray}
with the modified TVB minmod function defined as
\begin{eqnarray}\label{EQ:chap03TVBminmod}
\widetilde{m}(a_1,a_2,a_3)=\left\{ \begin{array}{lc}
a_1, & |a_1|\leq M h_{i,m}^2,\\
m(a_1,a_2,a_3), & otherwise,
\end{array} \right.
\end{eqnarray}
where the constant $M$ is chosen depending on the upper bound of the absolute value of the second derivative of the solution at the local extrema \cite{cockburn1989tvb2},
and $m(a,b,c)$ is the minmod function, defined as
\begin{eqnarray}
m(a,b,c)=\begin{cases}
s\cdot \min\{|a|,|b|,|c|\}, & \mbox{if } \mbox{sign}(a)= \mbox{sign}(a)=\mbox{sign}(c)=:s,
\\
0, &\mbox{otherwise}.
\end{cases}
\end{eqnarray}
 If any minmod function in $\delta_+u_{i,m}^{(mod)}, \delta_-u_{i,m}^{(mod)}$ is triggered, that is, if either $\delta_+u_{i,m}^{(mod)}\neq\delta_+u_{i,m}$  or  $\delta_-u_{i,m}^{(mod)}\neq\delta_-u_{i,m}$, the CV will be marked as a troubled cell. Consequently, its reconstructed polynomial function $u_h^n(x), x\in C_{i,m}$ requires further limitation. It should be noted that the TVB detector is one of the three detectors suggested for troubled cell identification in \cite{qiu2005indicators}, and it may not necessarily be the optimal choice. For high-order SV schemes, both the KXRCF detector and the Harten detector fail to detect CVs with internal discontinuities since the initial reconstruction in the SV is continuous. Therefore, we will not elaborate on these detectors and instead opt for the TVB detector due to its simplicity and practicality.

Assuming $C_{i,m}$ is a troubled cell, we now introduce the SWENO limiting procedure outlined in \cite{zhu2019swenoRKDG}
to reconstruct a new solution polynomial. First, select $C_{i,m}$ and $\lfloor \bold{k}/2 \rfloor$ control volumes on either side of $C_{i,m}$ to form
a large stencil $T_0$ centered on $C_{i,m}$. Subsequently, based on the average value of each control volume in $T_0$, we can
reconstruct a $(\bold{k}-1)$th order polynomial $p^0_{i,m}(x)$ satisfying the least square rules,
\begin{equation}\label{EQ:sec3p01d}
\frac{1}{V}\int_{C_{i,m}} p^0_{i,m}(x)dx=\overline{u}_{i,m}, min \sum_{C_{l,r}\in \widetilde{T}_0}(\frac{1}{V}\int_{C_{l,r}} p^0_{i,m}(x)dx-\overline{u}_{l,r})^2,
\end{equation}
where $\widetilde{T}_0=T_0-C_{i,m}$.
For $\bold{k}=2,3$, the stencil $T_0=\{C_{i,m-1},C_{i,m},C_{i,m+1}\}$. For $\bold{k}=4,5$, the stencil $T_0=\{C_{i,m-2},C_{i,m-1},C_{i,m},C_{i,m+1},C_{i,m+2}\}$.
Choose two small stencils, namely  $T_1=\{C_{i,m-1},C_{i,m}\}$ and $T_2=\{C_{i,m},C_{i,m+1}\}$.
We can easily get two linear polynomials $p^1_{i,m}(x)$ and $p^2_{i,m}(x)$ from stencil $T_1$ and $T_2$, which satisfy
\begin{eqnarray}\label{EQ:sec3p11d}
\int_{C_{i,r}} p^1_{i,m}(x) dx=\overline{u}_{i,r}^n,r=m-1,m,\\
\int_{C_{i,r}} p^2_{i,m}(x) dx=\overline{u}_{i,r}^n,r=m,m+1.
\end{eqnarray}

Similar to the ideas in \cite{zhu2016sweno,zhu2017sweno2,zhu2017sweno3,zhu2018sweno4,zhu2019swenoRKDG}, for given linear weights $\gamma_0, \gamma_1$ and $\gamma_2$, we can rewrite $p^0_{i,m}(x)$  into the following form, where subscripts $i, m$ are omitted for simplicity.
\begin{eqnarray}\label{EQ:sec3p01d2}
p^0(x)=\gamma_0(\frac{1}{\gamma_0}p^0(x)-\frac{\gamma_1}{\gamma_0}p^1(x)-\frac{\gamma_2}{\gamma_0}p^2(x))+\gamma_1 p^1(x)+\gamma_2 p^2(x).
\end{eqnarray}
where $\gamma_0, \gamma_1$ and $\gamma_2$ are linear weights, which can be any set of positive numbers that add up to 1.
Reference \cite{zhu2019swenoRKDG} analyzed the numerical performance of different linear weights and concluded that those weights do not significantly impact  the accuracy of the new scheme. Based on the practical experience from \cite{zhu2016sweno,zhu2017sweno2,zhu2017sweno3,zhu2018sweno4,zhu2019swenoRKDG}, we select the linear weights as $\gamma_0=0.8, \gamma_1=\gamma_2=0.1$. Subsequently, we compute the smoothness indicators,
\begin{eqnarray}\label{EQ:sec3beta}
\beta_l=\sum_{q=1}^{\bold{k}}\int_{C_{i,m}}(h_{i,m})^{2q-1}(\frac{\partial^q \widetilde{p}^l(x)}{\partial x^q})^2dx,l=0,1,2.
\end{eqnarray}
where  $\widetilde{p}^0(x)=\frac{1}{\gamma_0}p^0(x)-\frac{\gamma_1}{\gamma_0}p^1(x)-\frac{\gamma_2}{\gamma_0}p^2(x)$,
  $\widetilde{p}^1(x)=p^1(x)$ and $\widetilde{p}^2(x)=p^2(x)$ are the three linearly weighted polynomials of $p^0_{i,j}(x)$ in Eq. \eqref{EQ:sec3p01d2}. The parameters $\beta_l, l=0,1,2,$  quantify the smoothness of $\widetilde{p}^l(x)$  within the target cell. Let
\begin{eqnarray}\label{EQ:sec3tau}
\tau=(\frac{\mid\beta_0-\beta_1\mid+\mid\beta_0-\beta_2\mid}{2})^2,
\end{eqnarray}
so the nonlinear weight can be defined as follows,
\begin{eqnarray}\label{EQ:sec3weight}
\omega_l=\frac{\widetilde{\omega}_l}{\sum_{r=0}^2\widetilde{\omega}_l},
\widetilde{\omega}_l=\gamma_l(1+\frac{\tau}{\beta_l+\epsilon}),l=0,1,2.
\end{eqnarray}
Here, $\epsilon$ is a small positive constant introduced to prevent the denominator from becoming zero. Consequently, the reconstructed polynomial for the target cell can be expressed as
\begin{eqnarray}\label{EQ:sec3p0new1d}
p^{0}_{new}(x)=\omega_0\widetilde{p}^0(x)+\omega_1 p^1(x)+\omega_2p^2(x).
\end{eqnarray}


\begin{remark}\label{remark3}
After applying the limiter, the reconstruction function across the entire computational domain exhibits the following properties:\\
1. The reconstructed polynomial remains continuous at the CV boundaries that are not around the troubled CVs.\\
2. The reconstructed polynomial becomes discontinuous at the SV boundaries and at the CV boundaries where the limiter is active.
\end{remark}

Therefore, the state value at the CV boundary upon the limitation interface can be calculated using Eq. \eqref{EQ:sec3p0new1d}.
According to property \ref{remark3}, an analytical flux is employed at the CV boundary where the internal state remain continuous.
However, at the SV boundary or the boundaries of troubled CVs, the numerical flux is calculated using the Riemann solver.


\subsection{Two-dimensional case}
The two-dimensional limiting procedure on the rectangular meshes is as simple as the one-dimensional case, wherein we employ the
modified TVB minmod function to identify troubled cells. In fact, we can fix variables in one direction, and deal with variables in another direction in the
same way as in the one-dimensional case. The details are shown below.
For two-dimensional control volume, let
\begin{eqnarray}
\delta_+^xu_{i,m,j,n}=u_{i,m+\frac12,j,n}-\overline{u}_{i,m,j,n},\quad\quad
\delta_-^xu_{i,m,j,n}=\overline{u}_{i,m,j,n}-u_{i,m-\frac12,j,n},\\
\Delta_+^x\overline{u}_{i,m,j,n}=\overline{u}_{i,m+1,j,n}-\overline{u}_{i,m,j,n},\quad\quad
\Delta_-^x\overline{u}_{i,m,j,n}=\overline{u}_{i,m,j,n}-\overline{u}_{i,m-1,j,n},
\end{eqnarray}
and
\begin{eqnarray}
\delta_+^yu_{i,m,j,n}=u_{i,m,j,n+\frac12}-\overline{u}_{i,m,j,n},\quad\quad
\delta_-^yu_{i,m,j,n}=\overline{u}_{i,m,j,n}-u_{i,m,j,n-\frac12},\\
\Delta_+^y\overline{u}_{i,m,j,n}=\overline{u}_{i,m,j,n+1}-\overline{u}_{i,m,j,n},\quad\quad
\Delta_-^y\overline{u}_{i,m,j,n}=\overline{u}_{i,m,j,n}-\overline{u}_{i,m,j,n-1}.
\end{eqnarray}
where $u_{i,m+\frac12,j,n}$ and $u_{i,m-\frac12,j,n}$ denote the interface values of the control volume
reconstructed from the average value of the control volumes in the x-direction, with the y-direction (j,n) held constant,
and  $u_{i,m,j,n+\frac12}$ and $u_{i,m,j,n-\frac12}$ represent the interface values in the y-direction.
For simplicity, we drop subscripts $i,m,j,n$ and define
\begin{eqnarray}
\delta_+^xu^{(mod)}=\widetilde{m}(\delta_+^xu,\Delta_+^x\overline{u},\Delta_-^x\overline{u}),\quad\quad
\delta_-^xu^{(mod)}=\widetilde{m}(\delta_-^xu,\Delta_+^x\overline{u},\Delta_-^x\overline{u}),
\end{eqnarray}
and
\begin{eqnarray}
\delta_+^yu^{(mod)}=\widetilde{m}(\delta_+^yu,\Delta_+^y\overline{u},\Delta_-^y\overline{u}),\quad\quad
\delta_-^yu^{(mod)}=\widetilde{m}(\delta_-^yu,\Delta_+^y\overline{u},\Delta_-^y\overline{u}).
\end{eqnarray}
If any minmod function in  $\delta_+^xu^{(mod)}$, $\delta_-^xu^{(mod)}$, $\delta_+^yu^{(mod)}$, and $\delta_-^yu^{(mod)}$ is triggered, the control volume will be flagged as a troubled cell and its reconstructed polynomial function $u_h^n(x,y), (x,y)\in C_{i,m,j,n}$ will undergo further limitation.

Assuming $C_{i,m,j,n}$ is a troubled cell, similar to the one-dimensional case,
we first select $C_{i,m,j,n}$ and $(\lfloor \bold{k}/2 \rfloor)^2$ control volumes on either side of $C_{i,m,j,n}$
in x-direction and y-direction to form a large stencil $T_0$ centered on $C_{i,m,j,n}$.
We can reconstruct a $(\bold{k}-1)$th order polynomial $p^0_{i,m,j,n}(x,y)$ satisfying the least square rules,
here we drop subscripts $i,m,j,n$,
\begin{equation}\label{EQ:sec3p02d}
\frac{1}{V}\int_{C_{i,m,j,n}} p^0(x,y)dxdy=\overline{u}_{i,m,j,n}, min \sum_{C_{l,r,s,t}\in \widetilde{T}_0}(\frac{1}{V}\int_{C_{l,r,s,t}} p^0(x,y)dxdy-\overline{u}_{l,r,s,t})^2,
\end{equation}
where $\widetilde{T}_0=T_0-C_{i,m,j,n}$.
The neighbors of $C_{i,m,j,n}$ are selected to constitute four stencils, namely
$T_1=\{C_{i,m-1,j,n},$ $C_{i,m,j,n},C_{i,m,j,n-1}\}$, $T_2=\{C_{i,m,j,n},C_{i,m+1,j,n},C_{i,m,j,n-1}\}$,
 $T_3=\{C_{i,m-1,j,n},C_{i,m,j,n},$  $C_{i,m,j,n+1}\}$, and $T_4=\{C_{i,m,j,n},C_{i,m+1,j,n},C_{i,m,j,n+1}\}$.
By utilizing CV averages, the linear polynomials  $p^1(x,y)$, $p^2(x,y)$, $p^3(x,y)$ and $p^4(x,y)$ corresponding to these four stencils can be readily derived.
Subsequently, we express $p^0(x,y)$  as
\begin{eqnarray}\label{EQ:sec3p0new2d}
p^0(x,y)=\gamma_0\widetilde{p}^0(x,y)+\gamma_1 p^1(x,y)+\gamma_2p^2(x,y)+\gamma_3p^3(x,y)+\gamma_4p^4(x,y).
\end{eqnarray}
where $\gamma_0, \gamma_1, \gamma_2, \gamma_3,$ and $\gamma_4$ are linear weights.
We define $\widetilde{p}^0(x,y)=\frac{1}{\gamma_0}p^0(x,y)-\frac{\gamma_1}{\gamma_0}p^1(x,y)-\frac{\gamma_2}{\gamma_0}p^2(x,y)
-\frac{\gamma_3}{\gamma_0}p^3(x,y)-\frac{\gamma_4}{\gamma_0}p^4(x,y)$.
 Here we set $\gamma_0=0.8, \gamma_1=\gamma_2=\gamma_3=\gamma_4=0.05$. Subsequently, we calculate the smoothness indicators of the five linearly combined polynomials on $C_{i,m,j,n}$,
\begin{eqnarray}\label{EQ:sec3beta2d}
\beta_l=\sum_{q=1}^{\bold{k}}\int_{C_{i,m,j,n}}(h_{i,m}^x)^{2q_1-1}(h_{j,n}^y)^{2q_2-1}(\frac{\partial^q \widetilde{p}^l(x,y)}{\partial x^{q1}\partial y^{q2}})^2dxdy, q=q1+q2,l=0,...,4.
\end{eqnarray}
where $\widetilde{p}^l(x,y)=p^l(x,y)$ for $l=1,2,3,4$.
Let
\begin{eqnarray}\label{EQ:sec3tau2d}
\tau=(\frac{\mid\beta_0-\beta_1\mid+\mid\beta_0-\beta_2\mid+\mid\beta_0-\beta_3\mid+\mid\beta_0-\beta_4\mid}{4})^2,
\end{eqnarray}
non-linear weights can be defined as follows:
\begin{eqnarray}\label{EQ:sec3weight2d}
\omega_l=\frac{\widetilde{\omega}_l}{\sum_{r=0}^4\widetilde{\omega}_l},\quad
\widetilde{\omega}_l=\gamma_l(1+\frac{\tau}{\beta_l+\epsilon}),\quad l=0,1,2,3,4.
\end{eqnarray}
Here, $\epsilon$ is a small positive constant to prevent the denominator from becoming zero. The reconstructed polynomial for the target cell can be expressed as
\begin{eqnarray}\label{EQ:sec3p0new2d1}
p^{0,new}(x,y)=\omega_0\widetilde{p}^0(x,y)+\omega_1 p^1(x,y)+\omega_2p^2(x,y)+\omega_3p^3(x,y)+\omega_4p^4(x,y).
\end{eqnarray}

For both one- and two-dimensional system cases, local characteristic field decomposition is employed to improve the stability
of the limiter.
\section{Numerical test}
In this section, we will employ the proposed higher order spectral volume scheme to solve several fluid dynamics problems in order to evaluate the accuracy and effectiveness of the methods. The numerical performance of the second- through fifth- order spectral volume schemes is assessed. For simplicity, {\tt SV-cvMSWNO2}-\texttt{SV-cvMSWENO5} are used to denote the second- through fifth- order spectral volume schemes, respectively. In a $\bold{k}$th-order scheme, the computational domain $\Omega$ is uniformly partitioned into $N$ spectral volume elements, and within each element, $\bold{k}+1$ Gauss-Lobatto points are selected to form $\bold{k}$ control volumes, i.e
\begin{equation}
x_{i,m+\frac12}=x_{i-\frac12}+\frac{h_i}{2}(1-\cos(\frac{m\pi}{\bold{k}})),\quad m=0,\cdots,\bold{k}.
\end{equation}
For two-dimensional case, the same kind of division is carried out for the y direction.
In this chapter, unless otherwise stated, the specific heat ratio $\gamma$ in Euler equation is set to 1.4.
The $CFL$ number is $0.5$ for all simulation. Additionally, the constant $\epsilon$ in the nonlinear power calculation formula  Eq. \eqref{EQ:sec3weight} and  Eq. \eqref{EQ:sec3weight2d} is set to $10^{-6}$.

\subsection{Accuracy and convergence results for smooth solutions}
In this section, we examine the error and convergence order for one- and two-dimensional linear convection equations as well as smooth solutions of Euler equations. We define $l^1,l^2$ and $l^\infty$ errors as $\delta^p_N$ and the corresponding convergence order as $R^p_N$,
\begin{align*}
\delta^p_N&=\sum_{i=1}^{N}\sum_{m=1}^{\bold{k}}\|\rho_{i,j}-\rho_{i,j}^{exact}\|_p, \quad p=1,2\\
\delta_N^\infty&=\max_{i=1...N}\max_{m=1...\bold{k}}|\rho_{i,j}-\rho_{i,j}^{exact}|
\end{align*}
\begin{equation*}
R^p_N=\log_2(\frac{\delta^p_{N/2}}{\delta^p_N}),\quad p=1,2,\infty.
\end{equation*}

\begin{example}[Accuracy test of the one-dimensional linear convection equation]\label{chap03:sin1d}\rm
Consider the linear convection equation,
\begin{align}
 u_t+u_x=0,\quad\quad -1\le x\le 1,\\
 u(x,0)=\sin(\pi x).
\end{align}
This is a commonly employed example for assessing accuracy, which illustrates a one-dimensional fluid sine wave propagating at a constant velocity of 1 within the domain $\Omega=[-1,1]$. The problem possesses an exact solution given by the following form,
\begin{equation*}
u(x,t)=\sin\big(\pi(x-t)\big).
\end{equation*}

In the numerical calculation, periodic boundary conditions are applied at the boundaries, and the computation terminations at $t=1$. In order to evaluate the accuracy of the limiter, Table \ref{table:chap03sin1d1} shows the $l^1,l^2$ and $l^\infty$ errors along with the corresponding convergence rates for each order scheme when fully restricted by the cvMSWEMO limiter, meaning all elements are subjected to the limiter. The results indicate that schemes of 3 to 5 maintain an ideal convergence rate. For the second order scheme, the effect of the limiter resembles that of a TVD limiter, leading to a slight decrease in the convergence rate at the extreme points. Overall, the cvMSWEMO limiter preserves high order accuracy effectively. Table \ref{table:chap03sin1d2} presents the errors and corresponding convergence rates for each order scheme when the parameter $M$ in TVB minmod function Eq. \eqref{EQ:chap03TVBminmod} is set to 2. At this setting, no cell is identified as troubled cells. The results indicate that all order schemes achieve their theoretical convergence rates. Table \ref{table:chap03sin1d3} displays the errors and corresponding convergence rates when the parameter $M$ is set to be 0.01, along with the proportion of troubled cells. From these results, it can be observed that the 4th-order scheme exhibits a slight reduction in convergence rate, the $l^\infty$ error of the 2nd-order scheme shows a minor decrease, while the 3rd-order and 5th-order schemes maintain their theoretical convergence rates.

\begin{table}[htbp]
 \begin{center}   \caption{Example \ref{chap03:sin1d}: Numerical error and convergence rate of {\tt SV-cvMSWNO2}-\texttt{SV-cvMSWENO5} schemes
 completely restricted by cvMSWEMO limiter for at $t=1$.}
 \label{table:chap03sin1d1}
   \begin{tabular}{cccccccc}
     \hline
  & $N$   & $\delta_N^1$ & $R_N^1$& $\delta_N^2$ & $R_N^{2}$& $\delta_N^\infty$ & $R_N^{\infty}$
     \\
     \hline

{\tt SV-cvMSWNO2}   &5&    2.89E-01&   --- &   3.16E-01&   --- &            3.65E-01&   --- \\
  &10&            1.11E-01&   1.38&            1.18E-01&   1.42&            1.52E-01&   1.26\\
  &20&            4.03E-02&   1.46&            4.43E-02&   1.41&            6.56E-02&   1.22\\
  &40&            1.14E-02&   1.82&            1.43E-02&   1.63&            2.74E-02&   1.26\\
  &60&            5.45E-03&   1.83&            7.31E-03&   1.65&            1.63E-02&   1.28\\
  &80&            3.18E-03&   1.87&            4.54E-03&   1.66&            1.12E-02&   1.29\\
 &100&            2.10E-03&   1.86&            3.13E-03&   1.66&            8.40E-03&   1.30\\

{\tt SV-cvMSWNO3}&   5&    2.16E-02&   --- &   2.37E-02&   --- &            3.28E-02&   --- \\
 & 10&            2.72E-03&   2.99&            2.77E-03&   3.10&            3.27E-03&   3.33\\
 & 20&            3.46E-04&   2.98&            3.47E-04&   3.00&            3.51E-04&   3.22\\
 & 40&            4.35E-05&   2.99&            4.35E-05&   3.00&            4.36E-05&   3.01\\
 & 60&            1.29E-05&   3.00&            1.29E-05&   3.00&            1.29E-05&   3.00\\
 & 80&            5.44E-06&   3.00&            5.44E-06&   3.00&            5.45E-06&   3.00\\
 &100&            2.79E-06&   3.00&            2.79E-06&   3.00&            2.79E-06&   3.00\\

{\tt SV-cvMSWNO4} &  5&   5.38E-03&   --- &    6.04E-03&   --- &            1.04E-02&   --- \\
  &10&            6.01E-05&   6.49&            8.62E-05&   6.13&            1.76E-04&   5.89\\
  &20&            2.41E-06&   4.64&            2.52E-06&   5.10&            3.65E-06&   5.59\\
  &40&            1.46E-07&   4.04&            1.55E-07&   4.03&            1.85E-07&   4.30\\
  &60&            2.86E-08&   4.02&            3.05E-08&   4.00&            3.75E-08&   3.94\\
  &80&            9.02E-09&   4.02&            9.63E-09&   4.00&            1.20E-08&   3.96\\
 &100&            3.68E-09&   4.01&            3.94E-09&   4.00&            4.95E-09&   3.97\\

{\tt SV-cvMSWNO5}&   5&    6.00E-03&   --- &   7.72E-03&   --- &            1.58E-02&   --- \\
  &10&            1.67E-05&   8.49&            2.57E-05&   8.23&            5.50E-05&   8.17\\
  &20&            1.01E-07&   7.37&            1.45E-07&   7.47&            4.80E-07&   6.84\\
  &40&            2.41E-09&   5.38&            2.44E-09&   5.89&            3.94E-09&   6.93\\
  &60&            3.17E-10&   5.00&            3.18E-10&   5.02&            3.58E-10&   5.92\\
  &80&            7.53E-11&   5.00&            7.53E-11&   5.00&            7.72E-11&   5.33\\
 &100&            2.47E-11&   5.00&            2.47E-11&   5.00&            2.51E-11&   5.03\\

      \hline
   \end{tabular}
 \end{center}
 \end{table}

 \begin{table}[htbp]
 \begin{center}   \caption{Example \ref{chap03:sin1d}: Numerical error and convergence rate of {\tt SV-cvMSWNO2}-\texttt{SV-cvMSWENO5} schemes at $t=1$ with $M=2$.}
 \label{table:chap03sin1d2}
   \begin{tabular}{cccccccc}
     \hline
    & $N$ & $\delta_N^1$ & $R_N^1$& $\delta_N^2$ & $R_N^{2}$& $\delta_N^\infty$ & $R_N^{\infty}$
     \\
     \hline
{\tt SV-cvMSWNO2} &   5&    1.42E-01&   --- &    1.46E-01&   --- &            1.66E-01&   --- \\
&  10&            4.08E-02&   1.80&            4.15E-02&   1.82&            4.46E-02&   1.90\\
&  20&            1.10E-02&   1.89&            1.10E-02&   1.92&            1.16E-02&   1.95\\
&  40&            2.81E-03&   1.96&            2.81E-03&   1.97&            2.90E-03&   2.00\\
&  60&            1.26E-03&   1.98&            1.26E-03&   1.98&            1.29E-03&   2.01\\
&  80&            7.11E-04&   1.99&            7.11E-04&   1.99&            7.22E-04&   2.01\\
& 100&            4.56E-04&   1.99&            4.56E-04&   1.99&            4.62E-04&   2.01\\
{\tt SV-cvMSWNO3}  &  5&    1.57E-02&   --- &  1.60E-02&   --- &            1.89E-02&   --- \\
 & 10&            1.95E-03&   3.01&            2.06E-03&   2.96&            2.65E-03&   2.84\\
 & 20&            2.47E-04&   2.98&            2.59E-04&   2.99&            3.62E-04&   2.87\\
 & 40&            3.10E-05&   2.99&            3.26E-05&   2.99&            4.62E-05&   2.97\\
 & 60&            9.21E-06&   3.00&            9.66E-06&   3.00&            1.38E-05&   2.97\\
 & 80&            3.89E-06&   3.00&            4.08E-06&   3.00&            5.86E-06&   2.99\\
 &100&            1.99E-06&   3.00&            2.09E-06&   3.00&            3.01E-06&   2.99\\
{\tt SV-cvMSWNO4}&   5&    1.06E-03&   --- &   1.28E-03&   --- &            2.28E-03&   --- \\
&  10&            6.85E-05&   3.96&            8.22E-05&   3.97&            1.61E-04&   3.83\\
&  20&            4.42E-06&   3.95&            5.08E-06&   4.01&            9.71E-06&   4.05\\
&  40&            2.79E-07&   3.99&            3.18E-07&   4.00&            6.14E-07&   3.98\\
&  60&            5.53E-08&   3.99&            6.29E-08&   4.00&            1.22E-07&   3.99\\
&  80&            1.75E-08&   3.99&            1.99E-08&   4.00&            3.85E-08&   4.00\\
& 100&            7.19E-09&   4.00&            8.16E-09&   4.00&            1.58E-08&   4.00\\
{\tt SV-cvMSWNO5} &    5&  6.29E-05&   --- &   8.12E-05&   --- &            1.82E-04&   --- \\
&  10&            1.90E-06&   5.05&            2.52E-06&   5.01&            6.58E-06&   4.79\\
&  20&            5.98E-08&   4.99&            8.09E-08&   4.96&            2.17E-07&   4.92\\
&  40&            1.98E-09&   4.92&            2.65E-09&   4.94&            7.00E-09&   4.95\\
&  60&            2.77E-10&   4.85&            3.52E-10&   4.98&            8.65E-10&   5.16\\
&  80&            6.69E-11&   4.93&            8.08E-11&   5.11&            1.81E-10&   5.44\\
& 100&            2.16E-11&   5.06&            2.50E-11&   5.26&            5.33E-11&   5.47\\
      \hline
   \end{tabular}
 \end{center}
 \end{table}

  \begin{table}[htbp]
 \begin{center}   \caption{Example \ref{chap03:sin1d}: Numerical error and convergence rate of {\tt SV-cvMSWNO2}-\texttt{SV-cvMSWENO5} schemes at $t =1$ with $M=0.01$.}
 \label{table:chap03sin1d3}
   \begin{tabular}{ccccccccc}
     \hline
    & $N$ & $\delta_N^1$ & $R_N^1$& $\delta_N^2$ & $R_N^{2}$& $\delta_N^\infty$ & $R_N^{\infty}$ & $percent$
     \\
     \hline

{\tt SV-cvMSWNO2} &   5&            3.03E-01&   --- &            3.26E-01&   --- &            3.86E-01&   --- &  30.00\\
&  10&            1.13E-01&   1.42&            1.18E-01&   1.47&            1.40E-01&   1.47&  30.00\\
&  20&            3.35E-02&   1.76&            3.99E-02&   1.56&            5.73E-02&   1.28&  40.00\\
&  40&            9.54E-03&   1.81&            1.31E-02&   1.61&            2.57E-02&   1.16&  20.00\\
&  60&            4.33E-03&   1.95&            6.65E-03&   1.67&            1.59E-02&   1.18&  10.00\\
&  80&            2.41E-03&   2.03&            4.03E-03&   1.74&            1.01E-02&   1.58&   3.33\\
& 100&            1.51E-03&   2.09&            2.72E-03&   1.76&            7.79E-03&   1.17&   7.50\\

{\tt SV-cvMSWNO3}  &   5&            4.34E-02&   --- &            4.54E-02&   --- &            7.25E-02&   --- &  46.67\\
&  10&            4.09E-03&   3.41&            4.53E-03&   3.32&            5.86E-03&   3.63&  46.67\\
&  20&            2.71E-04&   3.91&            2.96E-04&   3.94&            4.65E-04&   3.65&   6.67\\
&  40&            3.10E-05&   3.13&            3.29E-05&   3.17&            4.62E-05&   3.33&  13.33\\
&  60&            9.21E-06&   2.99&            9.70E-06&   3.01&            1.38E-05&   2.97&   6.67\\
&  80&            3.89E-06&   3.00&            4.09E-06&   3.00&            5.86E-06&   2.99&   4.44\\
& 100&            1.99E-06&   3.00&            2.09E-06&   3.00&            3.01E-06&   2.99&   3.33\\

{\tt SV-cvMSWNO4}&   5&            1.92E-02&   --- &            2.24E-02&   --- &            3.96E-02&   --- &  35.00\\
&  10&            6.70E-04&   4.85&            8.86E-04&   4.66&            1.40E-03&   4.82&  35.00\\
&  20&            6.86E-06&   6.61&            9.28E-06&   6.58&            2.69E-05&   5.70&  10.00\\
&  40&            3.54E-07&   4.28&            4.63E-07&   4.33&            1.62E-06&   4.05&  10.00\\
&  60&            7.12E-08&   3.95&            9.74E-08&   3.84&            3.72E-07&   3.63&   5.00\\
&  80&            2.35E-08&   3.85&            3.40E-08&   3.65&            1.63E-07&   2.87&   3.33\\
& 100&            1.01E-08&   3.78&            1.57E-08&   3.47&            8.61E-08&   2.85&   2.50\\

{\tt SV-cvMSWNO5} &   5&            3.75E-02&   --- &            4.05E-02&   --- &            6.51E-02&   --- &  40.00\\
&  10&            1.18E-04&   8.31&            1.64E-04&   7.94&            3.24E-04&   7.65&  40.00\\
&  20&            4.89E-07&   7.91&            9.12E-07&   7.49&            2.33E-06&   7.12&   8.00\\
&  40&            3.19E-09&   7.26&            5.00E-09&   7.51&            1.69E-08&   7.11&   8.00\\
&  60&            3.32E-10&   5.58&            4.34E-10&   6.03&            1.34E-09&   6.26&   4.00\\
&  80&            7.49E-11&   5.17&            9.35E-11&   5.34&            2.68E-10&   5.59&   2.67\\
& 100&            2.37E-11&   5.16&            2.84E-11&   5.34&            8.38E-11&   5.21&   2.00\\

      \hline
   \end{tabular}
 \end{center}
 \end{table}

To test the influence of $\epsilon$ values, we set $\epsilon=10^{-10}$. Table \ref{table:chap03sin1d4} presents the errors and the corresponding rates of convergence for parameter $M=0.01$, along with the proportion of troubled cells. The results indicate that the $l^\infty$ error of the 2nd-order scheme decreases slightly, while the 3rd-order and 4th-order schemes achieve their theoretical convergence rate. Notably, the 5-order scheme exhibits super convergence.

\begin{table}[htbp]
 \begin{center}   \caption{Example \ref{chap03:sin1d}: Numerical error and convergence rate of {\tt SV-cvMSWNO2}-\texttt{SV-cvMSWENO5} schemes at $t =1$ with $M=0.01$ and $\epsilon=10^{-10}$. }
 \label{table:chap03sin1d4}
   \begin{tabular}{ccccccccc}
     \hline
    & $N$ & $\delta_N^1$ & $R_N^1$& $\delta_N^2$ & $R_N^{2}$& $\delta_N^\infty$ & $R_N^{\infty}$ & $percent$
     \\
     \hline

{\tt SV-cvMSWNO2} &   5&            3.03E-01&   --- &            3.26E-01&   --- &            3.86E-01&   --- &  30.00\\
&  10&            1.13E-01&   1.42&            1.18E-01&   1.47&            1.40E-01&   1.47&  30.00\\
&  20&            3.35E-02&   1.76&            3.99E-02&   1.56&            5.74E-02&   1.28&  40.00\\
&  40&            9.56E-03&   1.81&            1.31E-02&   1.61&            2.57E-02&   1.16&  20.00\\
&  60&            4.33E-03&   1.95&            6.65E-03&   1.68&            1.58E-02&   1.20&  12.50\\
&  80&            2.41E-03&   2.04&            4.04E-03&   1.74&            1.01E-02&   1.56&   5.00\\
& 100&            1.52E-03&   2.08&            2.73E-03&   1.75&            7.85E-03&   1.12&   6.25\\

{\tt SV-cvMSWNO3}  &   5&            4.31E-02&   --- &            4.53E-02&   --- &            7.25E-02&   --- &  46.67\\
&  10&            4.25E-03&   3.34&            4.74E-03&   3.26&            6.24E-03&   3.54&  46.67\\
&  20&            2.97E-04&   3.84&            3.24E-04&   3.87&            4.67E-04&   3.74&   6.67\\
&  40&            3.54E-05&   3.07&            3.78E-05&   3.10&            7.56E-05&   2.63&  13.33\\
&  60&            9.21E-06&   3.32&            9.75E-06&   3.35&            1.38E-05&   4.19&   6.67\\
&  80&            3.89E-06&   3.00&            4.09E-06&   3.02&            5.86E-06&   2.99&   4.44\\
& 100&            1.99E-06&   3.00&            2.09E-06&   3.00&            3.01E-06&   2.99&   3.33\\

{\tt SV-cvMSWNO4}&   5&            2.00E-02&   --- &            2.35E-02&   --- &            4.14E-02&   --- &  40.00\\
&  10&            1.28E-03&   3.97&            1.81E-03&   3.69&            3.21E-03&   3.69&  40.00\\
&  20&            2.11E-04&   2.60&            3.75E-04&   2.27&            1.13E-03&   1.51&  10.00\\
&  40&            1.52E-06&   7.12&            2.81E-06&   7.06&            1.09E-05&   6.69&   7.50\\
&  60&            1.07E-07&   6.54&            1.80E-07&   6.77&            7.98E-07&   6.44&   5.00\\
&  80&            2.78E-08&   4.70&            4.57E-08&   4.76&            2.44E-07&   4.12&   3.33\\
& 100&            1.09E-08&   4.18&            1.78E-08&   4.22&            1.00E-07&   3.99&   2.50\\

{\tt SV-cvMSWNO5} &   5&            3.78E-02&   --- &            4.08E-02&   --- &            6.56E-02&   --- &  36.00\\
&  10&            7.36E-04&   5.68&            1.20E-03&   5.09&            3.40E-03&   4.27&  36.00\\
&  20&            4.11E-05&   4.16&            8.70E-05&   3.78&            2.17E-04&   3.97&   8.00\\
&  40&            1.75E-07&   7.87&            4.30E-07&   7.66&            1.51E-06&   7.17&   8.00\\
&  60&            8.93E-09&   7.34&            2.82E-08&   6.72&            1.16E-07&   6.33&   4.00\\
&  80&            8.98E-10&   7.99&            3.05E-09&   7.73&            1.46E-08&   7.20&   2.67\\
& 100&            1.64E-10&   7.63&            5.41E-10&   7.75&            2.70E-09&   7.55&   2.00\\

      \hline
   \end{tabular}
 \end{center}
 \end{table}

\end{example}

\begin{example}[Precision test of the one-dimensional Euler equations]\label{chap03:Eulersin1d}\rm
The problem describes the periodic propagation of a one-dimensional sine wave in the region $\Omega=[0,2]$. The initial condition are given by
\begin{equation*}
 \rho(x,0)=1+0.2\sin(\pi x),\quad U(x,0)=0.7, \quad P(x,0)=1.
\end{equation*}
An exact solution exists in the following form:
\begin{equation*}
\rho(x,t)=1+0.2\sin\big(\pi(x-0.7t)\big),
\quad U(x,t)=0.7, \quad P(x,t)=1.
\end{equation*}

Periodic boundary conditions are employed to compute the result at one period, that is, $t=2$. Table \ref{table:chap03Eulersin1d1} presents the errors and corresponding convergence rates for each order scheme when $M=2$ (without using a limiter). All schemes achieve their theoretical convergence orders. Table \ref{table:chap03Eulersin1d2} displays the complete limiting errors and corresponding convergence rates for each order scheme. All schemes reach their theoretical convergence orders, through the $L^\infty$ error convergence rate slightly decreases for the second-order scheme. Table \ref{table:chap03Eulersin1d3} illustrates the errors, convergence rates and proportions of troubled cells when $M=0.01$. In this scenario, both the fourth and fifth order schemes experience reductions in performance.

\begin{table}[htbp]
 \begin{center}   \caption{Example \ref{chap03:Eulersin1d}: Numerical error and convergence rate of {\tt SV-cvMSWNO2}-\texttt{SV-cvMSWENO5} schemes at $t =2$ with $M=2.0$. }
 \label{table:chap03Eulersin1d1}
   \begin{tabular}{cccccccc}
     \hline
    & $N$ & $\delta_N^1$ & $R_N^1$& $\delta_N^2$ & $R_N^{2}$& $\delta_N^\infty$ & $R_N^{\infty}$
     \\
     \hline

{\tt SV-cvMSWNO2} &   5&            7.85E-02&   --- &            6.24E-02&   --- &            6.39E-02&   --- \\
&  10&            1.77E-02&   2.15&            1.35E-02&   2.21&            1.41E-02&   2.18   \\
&  20&            5.59E-03&   1.66&            4.32E-03&   1.64&            5.15E-03&   1.45   \\
&  40&            1.61E-03&   1.79&            1.26E-03&   1.77&            1.86E-03&   1.47   \\
&  60&            7.56E-04&   1.87&            6.13E-04&   1.78&            1.06E-03&   1.40   \\
&  80&            4.27E-04&   1.98&            3.61E-04&   1.84&            6.85E-04&   1.51   \\
& 100&            2.76E-04&   1.97&            2.38E-04&   1.86&            4.87E-04&   1.53   \\

{\tt SV-cvMSWNO3}  &   5&            5.10E-03&   --- &            4.05E-03&   --- &            4.96E-03&   --- \\
&  10&            1.09E-03&   2.23&            8.68E-04&   2.22&            1.13E-03&   2.14   \\
&  20&            1.81E-04&   2.59&            1.43E-04&   2.60&            1.59E-04&   2.82   \\
&  40&            2.52E-05&   2.84&            1.99E-05&   2.84&            2.24E-05&   2.83   \\
&  60&            7.67E-06&   2.94&            6.05E-06&   2.94&            6.80E-06&   2.94   \\
&  80&            3.27E-06&   2.97&            2.58E-06&   2.97&            2.89E-06&   2.97   \\
& 100&            1.68E-06&   2.98&            1.33E-06&   2.98&            1.49E-06&   2.98   \\

{\tt SV-cvMSWNO4}  &   5&            3.94E-04&   --- &            3.27E-04&   --- &            4.45E-04&   --- \\
&  10&            2.01E-05&   4.29&            1.67E-05&   4.29&            2.53E-05&   4.13   \\
&  20&            1.16E-06&   4.12&            9.41E-07&   4.15&            1.39E-06&   4.19   \\
&  40&            7.31E-08&   3.99&            5.91E-08&   3.99&            8.80E-08&   3.98   \\
&  60&            1.44E-08&   4.00&            1.16E-08&   4.00&            1.73E-08&   4.01   \\
&  80&            4.57E-09&   4.00&            3.68E-09&   4.00&            5.48E-09&   4.00   \\
& 100&            1.87E-09&   4.00&            1.51E-09&   4.00&            2.25E-09&   4.00   \\

{\tt SV-cvMSWNO5}  &    5&            2.36E-05&   --- &            2.06E-05&   --- &            3.67E-05&   --- \\
&  10&            1.11E-06&   4.41&            9.46E-07&   4.44&            1.92E-06&   4.26   \\
&  20&            4.50E-08&   4.62&            3.82E-08&   4.63&            7.54E-08&   4.67   \\
&  40&            1.55E-09&   4.86&            1.30E-09&   4.88&            2.45E-09&   4.94   \\
&  60&            2.05E-10&   4.99&            1.69E-10&   5.03&            2.99E-10&   5.19   \\
&  80&            4.82E-11&   5.03&            3.95E-11&   5.06&            6.99E-11&   5.05   \\
& 100&            1.58E-11&   4.99&            1.32E-11&   4.92&            2.51E-11&   4.59   \\

      \hline
   \end{tabular}
 \end{center}
 \end{table}

\begin{table}[htbp]
 \begin{center}   \caption{Example \ref{chap03:Eulersin1d}:
  Numerical error and convergence rate of {\tt SV-cvMSWNO2}-\texttt{SV-cvMSWENO5} schemes completely restricted by cvMSWEMO limiter at $t=2$.}
 \label{table:chap03Eulersin1d2}
   \begin{tabular}{cccccccc}
     \hline
    & $N$ & $\delta_N^1$ & $R_N^1$& $\delta_N^2$ & $R_N^{2}$& $\delta_N^\infty$ & $R_N^{\infty}$
     \\
     \hline

{\tt SV-cvMSWNO2} &   5&            1.86E-01&   --- &            1.45E-01&   --- &            1.48E-01&   --- \\
&  10&            7.47E-02&   1.32&            6.39E-02&   1.18&            7.25E-02&   1.03 \\
&  20&            2.68E-02&   1.48&            2.34E-02&   1.45&            3.17E-02&   1.19 \\
&  40&            1.03E-02&   1.37&            9.04E-03&   1.37&            1.38E-02&   1.20 \\
&  60&            4.82E-03&   1.89&            4.62E-03&   1.66&            8.31E-03&   1.25 \\
&  80&            2.91E-03&   1.76&            2.90E-03&   1.62&            5.77E-03&   1.27 \\
& 100&            1.89E-03&   1.92&            2.01E-03&   1.65&            4.33E-03&   1.28 \\

{\tt SV-cvMSWNO3} &    5&            1.78E-02&   --- &            1.43E-02&   --- &            1.54E-02&   --- \\
&  10&            2.31E-03&   2.95&            1.82E-03&   2.97&            1.83E-03&   3.07\\
&  20&            2.91E-04&   2.99&            2.29E-04&   2.99&            2.30E-04&   2.99\\
&  40&            3.65E-05&   3.00&            2.87E-05&   3.00&            2.88E-05&   3.00\\
&  60&            1.08E-05&   3.00&            8.50E-06&   3.00&            8.51E-06&   3.00\\
&  80&            4.57E-06&   3.00&            3.59E-06&   3.00&            3.59E-06&   3.00\\
& 100&            2.34E-06&   3.00&            1.84E-06&   3.00&            1.84E-06&   3.00\\

{\tt SV-cvMSWNO4} &    5&            7.72E-04&   --- &            7.96E-04&   --- &            1.34E-03&   ---  \\
 &  10&            2.96E-05&   4.70&            2.52E-05&   4.98&            3.76E-05&   5.16 \\
 &  20&            1.86E-06&   4.00&            1.49E-06&   4.08&            1.86E-06&   4.33 \\
 &  40&            1.16E-07&   4.00&            9.20E-08&   4.02&            1.09E-07&   4.10 \\
 &  60&            2.29E-08&   3.99&            1.82E-08&   4.00&            2.10E-08&   4.06 \\
 &  80&            7.26E-09&   4.00&            5.74E-09&   4.00&            6.56E-09&   4.04 \\
 & 100&            2.97E-09&   4.00&            2.35E-09&   4.00&            2.67E-09&   4.04 \\

{\tt SV-cvMSWNO5} &    5&            1.81E-04&   --- &            1.84E-04&   --- &            3.23E-04&   ---  \\
&  10&            2.45E-06&   6.21&            2.09E-06&   6.46&            3.33E-06&   6.60 \\
&  20&            7.67E-08&   5.00&            6.05E-08&   5.11&            6.71E-08&   5.64 \\
&  40&            2.41E-09&   4.99&            1.89E-09&   5.00&            1.92E-09&   5.13 \\
&  60&            3.17E-10&   5.00&            2.49E-10&   5.00&            2.52E-10&   5.01 \\
&  80&            7.53E-11&   5.00&            5.91E-11&   5.00&            5.97E-11&   5.00 \\
& 100&            2.47E-11&   5.00&            1.94E-11&   5.00&            1.96E-11&   5.00 \\

      \hline
   \end{tabular}
 \end{center}
 \end{table}

\begin{table}[htbp]
 \begin{center}   \caption{Example \ref{chap03:Eulersin1d}: Numerical error and convergence rate of {\tt SV-cvMSWNO2}-\texttt{SV-cvMSWENO5} schemes at $t =2$ with $M=0.01$.}
 \label{table:chap03Eulersin1d3}
   \begin{tabular}{ccccccccc}
     \hline
    & $N$ & $\delta_N^1$ & $R_N^1$& $\delta_N^2$ & $R_N^{2}$& $\delta_N^\infty$ & $R_N^{\infty}$ & $percent$
     \\
     \hline

{\tt SV-cvMSWNO2} &   5&            1.74E-01&   --- &            1.36E-01&   --- &            1.41E-01&   --- &  40.00\\
&  10&            4.76E-02&   1.87&            4.47E-02&   1.61&            5.58E-02&   1.34&  40.00\\
&  20&            2.00E-02&   1.25&            1.69E-02&   1.40&            2.16E-02&   1.37&  25.00\\
&  40&            5.04E-03&   1.99&            4.95E-03&   1.77&            8.04E-03&   1.43&  10.00\\
&  60&            2.20E-03&   2.04&            2.38E-03&   1.81&            4.45E-03&   1.46&   8.33\\
&  80&            1.20E-03&   2.11&            1.40E-03&   1.85&            2.90E-03&   1.48&   6.25\\
& 100&            7.49E-04&   2.12&            9.24E-04&   1.85&            2.20E-03&   1.24&   4.00\\

{\tt SV-cvMSWNO3}  &   5&            8.23E-03&   --- &            6.41E-03&   --- &            8.03E-03&   --- &  33.33\\
&  10&            1.40E-03&   2.56&            1.19E-03&   2.43&            2.28E-03&   1.82&   6.67\\
&  20&            1.90E-04&   2.88&            1.59E-04&   2.91&            2.61E-04&   3.12&  13.33\\
&  40&            2.53E-05&   2.91&            2.06E-05&   2.94&            3.19E-05&   3.03&   6.67\\
&  60&            7.67E-06&   2.94&            6.16E-06&   2.98&            8.91E-06&   3.15&   4.44\\
&  80&            3.27E-06&   2.97&            2.60E-06&   2.99&            3.59E-06&   3.16&   3.33\\
& 100&            1.68E-06&   2.98&            1.34E-06&   2.99&            1.78E-06&   3.14&   2.67\\

{\tt SV-cvMSWNO4}&   5&            1.63E-03&   --- &            1.29E-03&   --- &            1.71E-03&   --- &  25.00\\
&  10&            5.45E-05&   4.90&            4.57E-05&   4.81&            5.45E-05&   4.97&  10.00\\
&  20&            3.43E-06&   3.99&            3.31E-06&   3.79&            6.50E-06&   3.07&  10.00\\
&  40&            2.61E-07&   3.72&            3.02E-07&   3.46&            7.79E-07&   3.06&   5.00\\
&  60&            5.82E-08&   3.70&            7.90E-08&   3.31&            2.38E-07&   2.92&   3.33\\
&  80&            2.15E-08&   3.47&            3.19E-08&   3.15&            1.03E-07&   2.91&   2.50\\
& 100&            1.02E-08&   3.32&            1.64E-08&   2.98&            5.36E-08&   2.93&   2.00\\

{\tt SV-cvMSWNO5} &   5&            5.35E-04&   --- &            5.01E-04&   --- &            7.84E-04&   --- &  16.00\\
&  10&            5.21E-06&   6.68&            5.88E-06&   6.41&            1.00E-05&   6.29&  12.00\\
&  20&            1.06E-07&   5.62&            1.23E-07&   5.58&            3.59E-07&   4.80&   8.00\\
&  40&            2.01E-09&   5.72&            1.99E-09&   5.95&            6.07E-09&   5.89&   4.00\\
&  60&            2.81E-10&   4.86&            2.66E-10&   4.96&            7.40E-10&   5.19&   2.67\\
&  80&            7.83E-11&   4.44&            7.43E-11&   4.43&            2.06E-10&   4.44&   2.00\\
& 100&            3.59E-11&   3.50&            3.70E-11&   3.13&            9.38E-11&   3.54&   1.60\\

      \hline
   \end{tabular}
 \end{center}
 \end{table}

As shown in Table \ref{table:chap03Eulersin1d4} when $M=0.01$ and $\epsilon=10^{-10}$, the error, convergence rate, and proportion of troubled cells for the 4th order scheme are reduced, while the 5th order scheme shows a slight improvement in accuracy.

 \begin{table}[htbp]
 \begin{center}   \caption{Example \ref{chap03:Eulersin1d}: Numerical error and convergence rate of {\tt SV-cvMSWNO2}-\texttt{SV-cvMSWENO5} schemes at $t=2$ with $M=0.01$ and $\epsilon=10^{-10}$.}
 \label{table:chap03Eulersin1d4}
   \begin{tabular}{ccccccccc}
     \hline
    & $N$ & $\delta_N^1$ & $R_N^1$& $\delta_N^2$ & $R_N^{2}$& $\delta_N^\infty$ & $R_N^{\infty}$ & $percent$
     \\
     \hline

{\tt SV-cvMSWNO2} &   5&            1.74E-01&   --- &            1.36E-01&   --- &            1.41E-01&   --- &  40.00\\
&  10&            4.78E-02&   1.86&            4.49E-02&   1.60&            5.60E-02&   1.33&  40.00\\
&  20&            2.00E-02&   1.26&            1.69E-02&   1.41&            2.16E-02&   1.37&  25.00\\
&  40&            5.05E-03&   1.99&            4.96E-03&   1.77&            8.05E-03&   1.43&  12.50\\
&  60&            2.20E-03&   2.05&            2.38E-03&   1.81&            4.45E-03&   1.46&   8.33\\
&  80&            1.20E-03&   2.11&            1.39E-03&   1.85&            2.91E-03&   1.48&   6.25\\
& 100&            7.48E-04&   2.12&            9.21E-04&   1.86&            2.19E-03&   1.27&   4.50\\

{\tt SV-cvMSWNO3}  &    5&            1.07E-02&   --- &            8.26E-03&   --- &            1.08E-02&   --- &  33.33\\
&  10&            1.66E-03&   2.69&            1.41E-03&   2.55&            2.79E-03&   1.96&   6.67\\
&  20&            1.93E-04&   3.10&            1.62E-04&   3.12&            2.68E-04&   3.38&  13.33\\
&  40&            2.53E-05&   2.94&            2.06E-05&   2.97&            3.19E-05&   3.07&   6.67\\
&  60&            7.67E-06&   2.94&            6.16E-06&   2.98&            8.91E-06&   3.15&   4.44\\
&  80&            3.27E-06&   2.97&            2.60E-06&   2.99&            3.59E-06&   3.16&   3.33\\
& 100&            1.68E-06&   2.98&            1.34E-06&   2.99&            1.78E-06&   3.14&   2.67\\

{\tt SV-cvMSWNO4}&   5&            4.15E-03&   --- &            3.39E-03&   --- &            4.64E-03&   --- &  25.00\\
&  10&            3.33E-04&   3.64&            3.34E-04&   3.34&            4.61E-04&   3.33&  10.00\\
&  20&            4.96E-06&   6.07&            5.44E-06&   5.94&            1.22E-05&   5.24&  10.00\\
&  40&            2.68E-07&   4.21&            3.12E-07&   4.12&            7.99E-07&   3.93&   5.00\\
&  60&            5.85E-08&   3.75&            7.95E-08&   3.38&            2.40E-07&   2.97&   3.33\\
&  80&            2.15E-08&   3.48&            3.19E-08&   3.17&            1.03E-07&   2.93&   2.50\\
& 100&            1.02E-08&   3.33&            1.64E-08&   2.99&            5.36E-08&   2.93&   2.00\\

{\tt SV-cvMSWNO5} &   5&            4.35E-03&   --- &            4.10E-03&   --- &            7.37E-03&   --- &  28.00\\
&  10&            2.51E-04&   4.12&            2.90E-04&   3.82&            5.55E-04&   3.73&   8.00\\
&  20&            1.56E-06&   7.33&            2.21E-06&   7.04&            5.89E-06&   6.56&   8.00\\
&  40&            1.06E-08&   7.21&            1.28E-08&   7.43&            3.12E-08&   7.56&   4.00\\
&  60&            6.59E-10&   6.84&            7.04E-10&   7.15&            1.64E-09&   7.27&   2.67\\
&  80&            1.35E-10&   5.51&            1.43E-10&   5.53&            3.49E-10&   5.38&   2.00\\
& 100&            5.30E-11&   4.19&            5.93E-11&   3.96&            1.57E-10&   3.59&   1.60\\

      \hline
   \end{tabular}
 \end{center}
 \end{table}

\end{example}

\subsection{Numerical results for discontinuous solutions}

\begin{example}[One-dimensional Sod shock tube problem]\label{chap03:sod1d}\rm
Considering the one-dimensional Sod shock tube problem, the computational domain is $\Omega=[-5,5]$, with the initial conditions specified as,
\begin{equation*}
(\rho,U,P)(x,0)=\left\{ \begin{array}{lc}
 (1,0,1),& x<0,\\
(0.125,0,0.1),&x>0, \end{array}\right.
\end{equation*}

The time is calculated to $t=2$. With wave propagation, the left sparse wave, contact discontinuity and right shock wave develop. When the parameter $M$ in the TVB detector is small, each order scheme yield stable numerical results. Figure \ref{Fig:sod1d} illustrates that when $M=0.01$, each order scheme employs 100 spectral volumes to obtain the density numerical solution and the corresponding troubled cells. An excessive number of troubled cells leads to significant dissipation, thereby preventing the high-order scheme from demonstrating its high-resolution advantage. To enhance computational accuracy, numerical experiments were conducted with varying values of $M$. Figure \ref{Fig:sod1d2} presents the density numerical solutions and corresponding troubled cells for third to fifth order schemes with $M$ set at 10, 20, and 50, respectively. The troubled cells are consistent with the wave structure, and the number of troubled cells in each order scheme is comparable, with the higher-order schemes showing slightly better resolution.

\begin{figure}[htbp]
    \centering
    \includegraphics[width=0.3\textwidth]{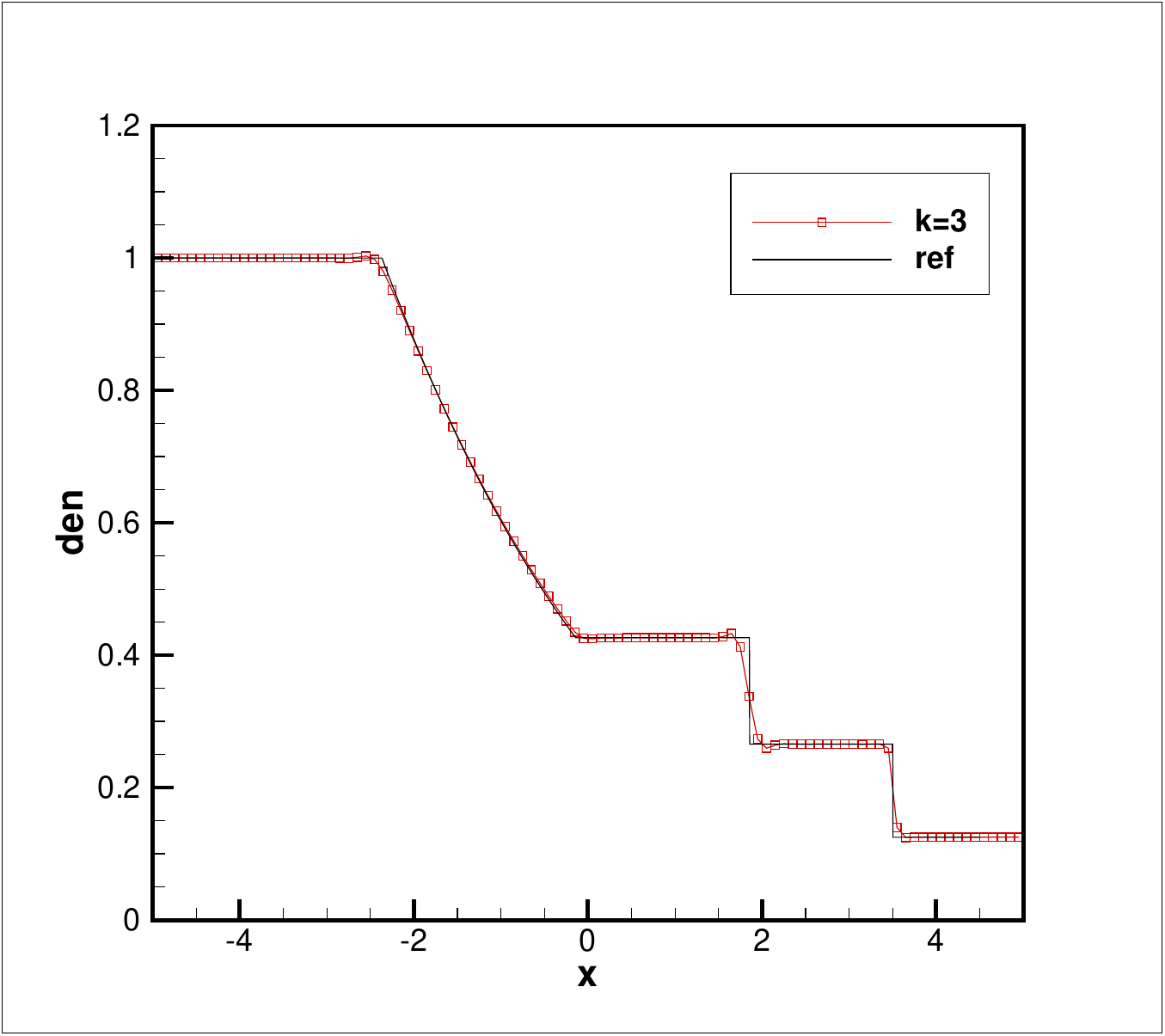}
    \includegraphics[width=0.3\textwidth]{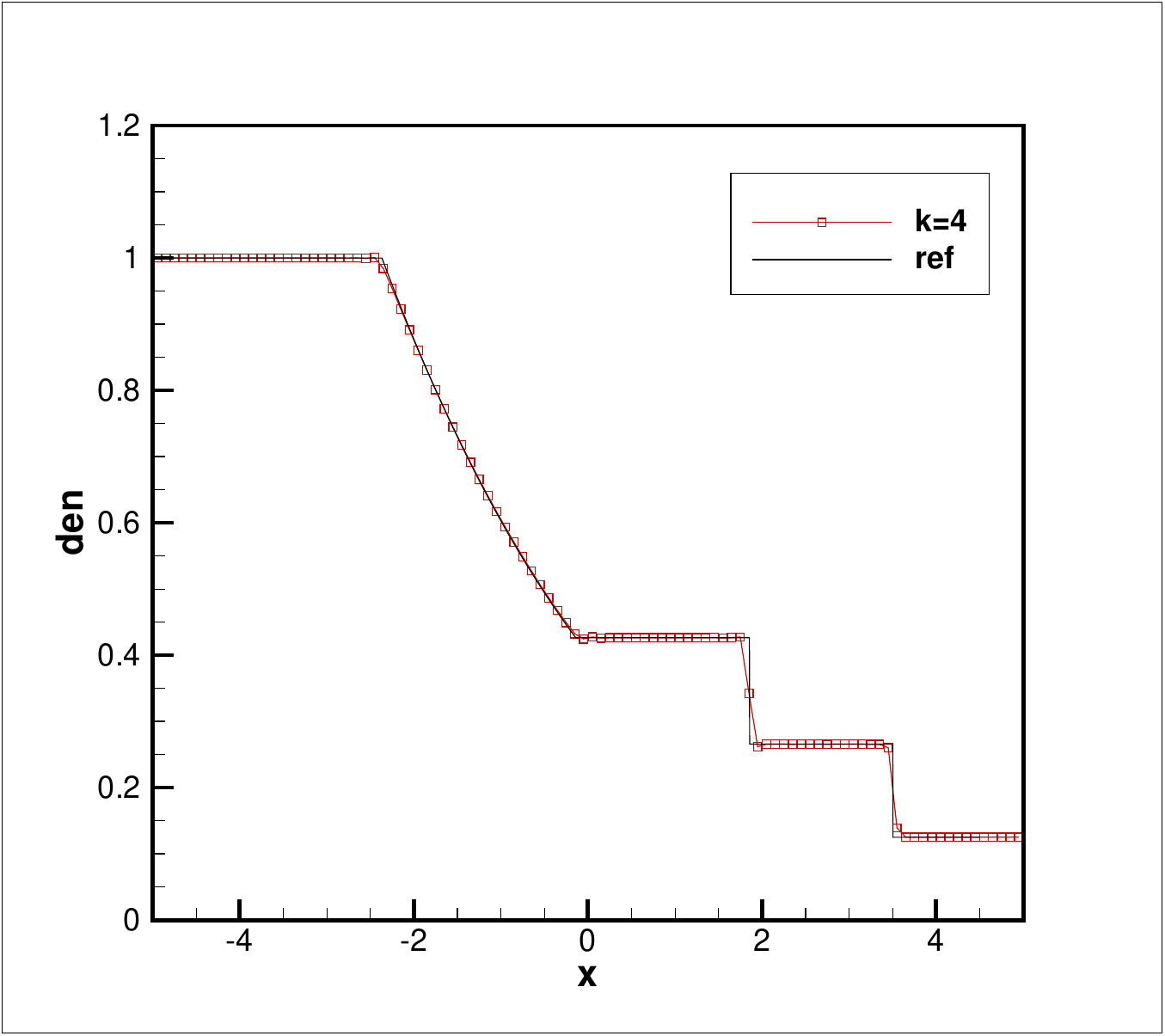}
    \includegraphics[width=0.3\textwidth]{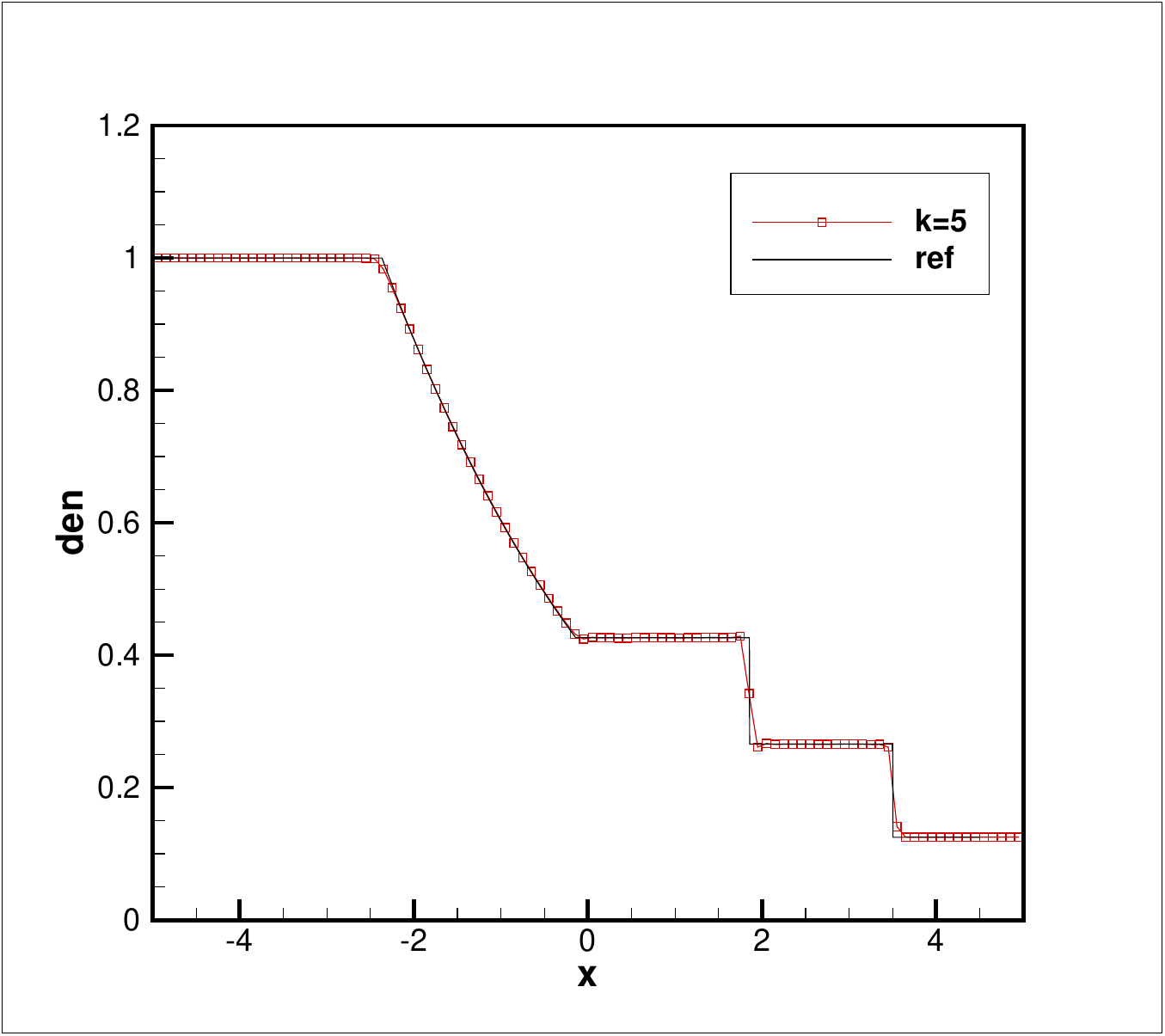}
    \includegraphics[width=0.3\textwidth]{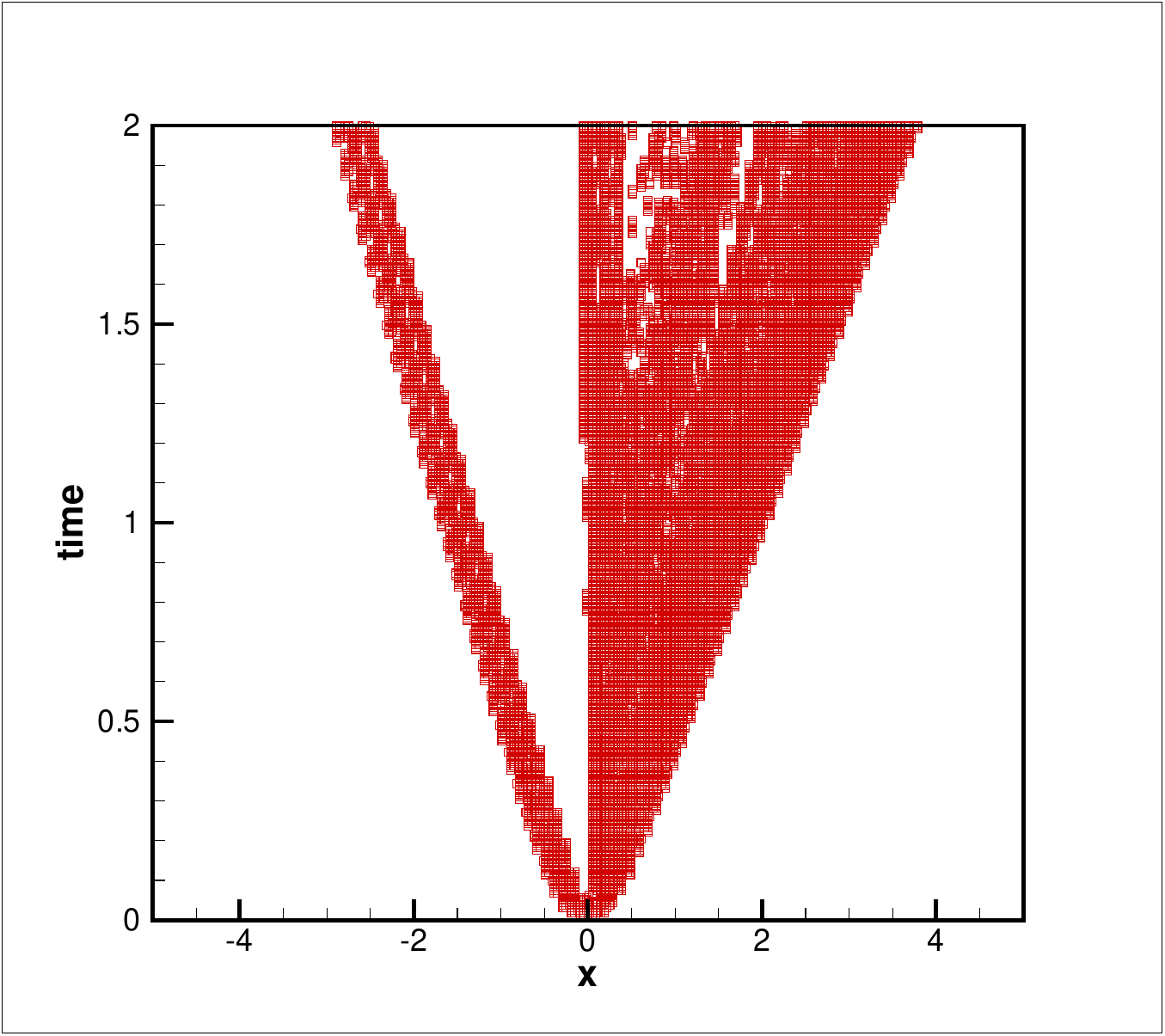}
    \includegraphics[width=0.3\textwidth]{sod1d/sod1d-markcv-m0.01-k4.pdf}
    \includegraphics[width=0.3\textwidth]{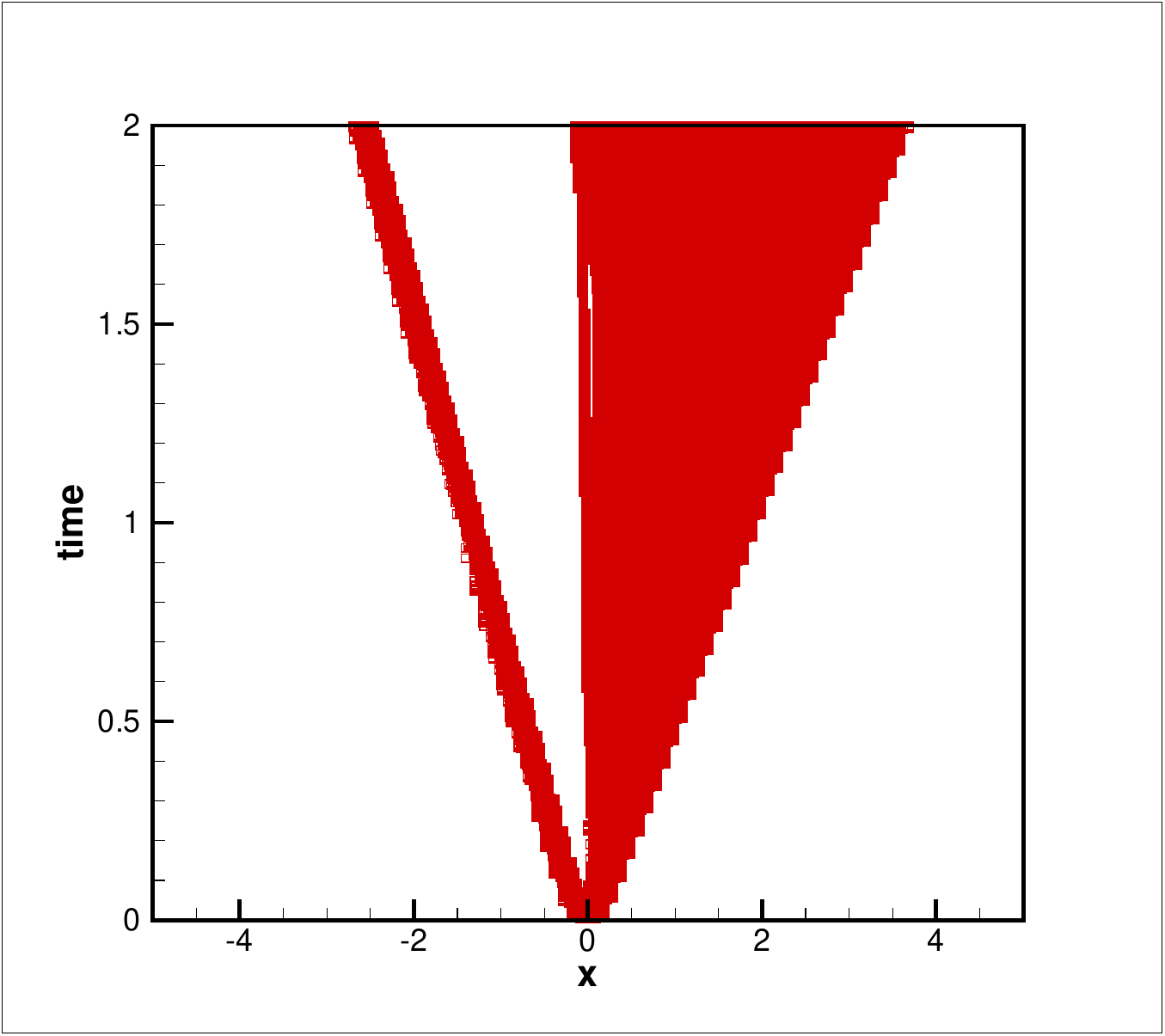}
     \caption{Example \ref{chap03:sod1d}: One-dimensional Sod shock tube problem: When $t=2$ and $M=0.01$,the density numerical solution is obtained using 100 spectral volumes in {\tt SV-cvMSWNO3}-\texttt{SV-cvMSWENO5} scheme. The solid line and symbol ``$\square$'' represent the exact solution and {\tt SV-cvMSWENO} results, respectively.}
    \label{Fig:sod1d}
 \end{figure}

 \begin{figure}[htbp]
    \centering
    \includegraphics[width=0.3\textwidth]{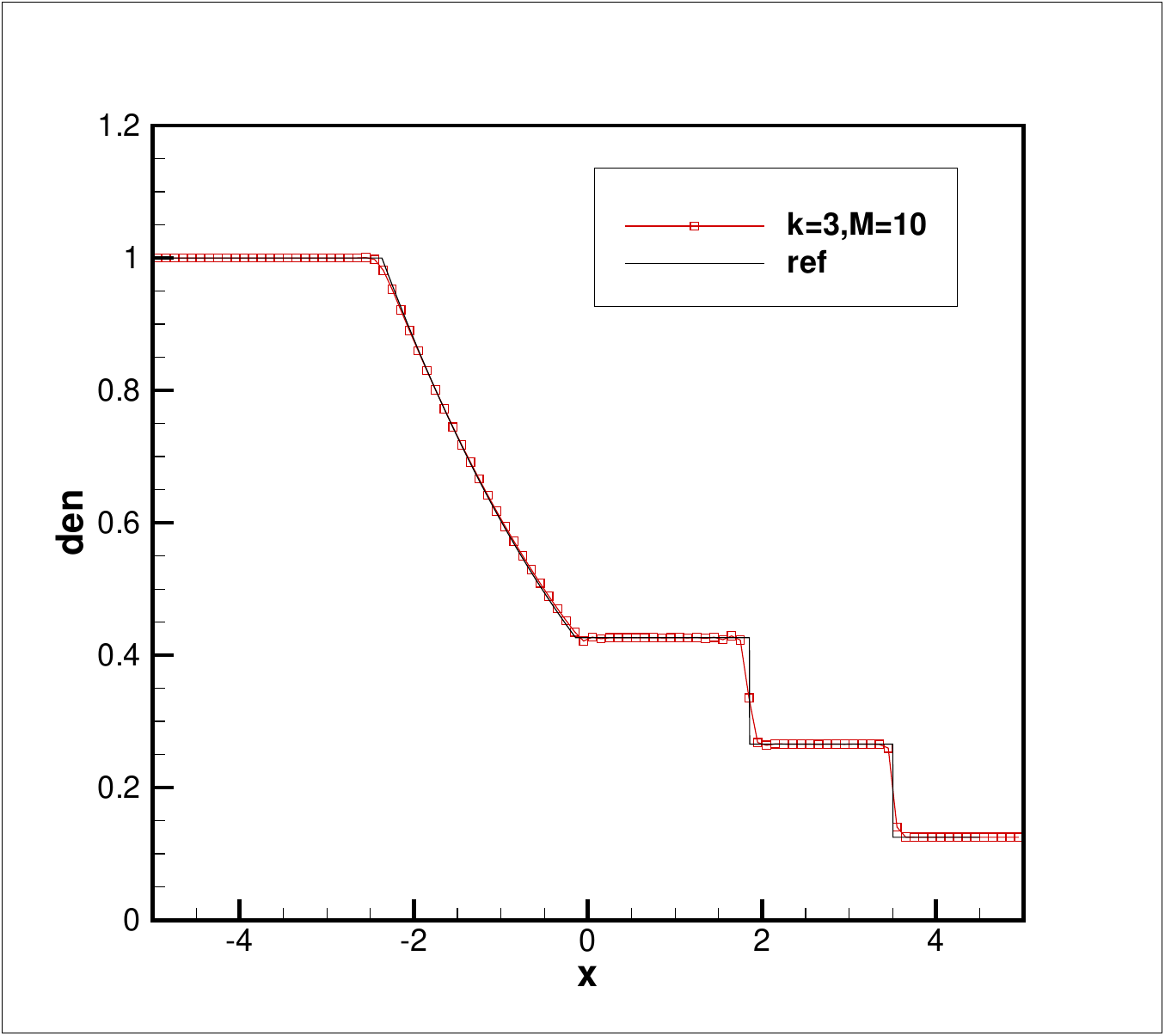}
    \includegraphics[width=0.3\textwidth]{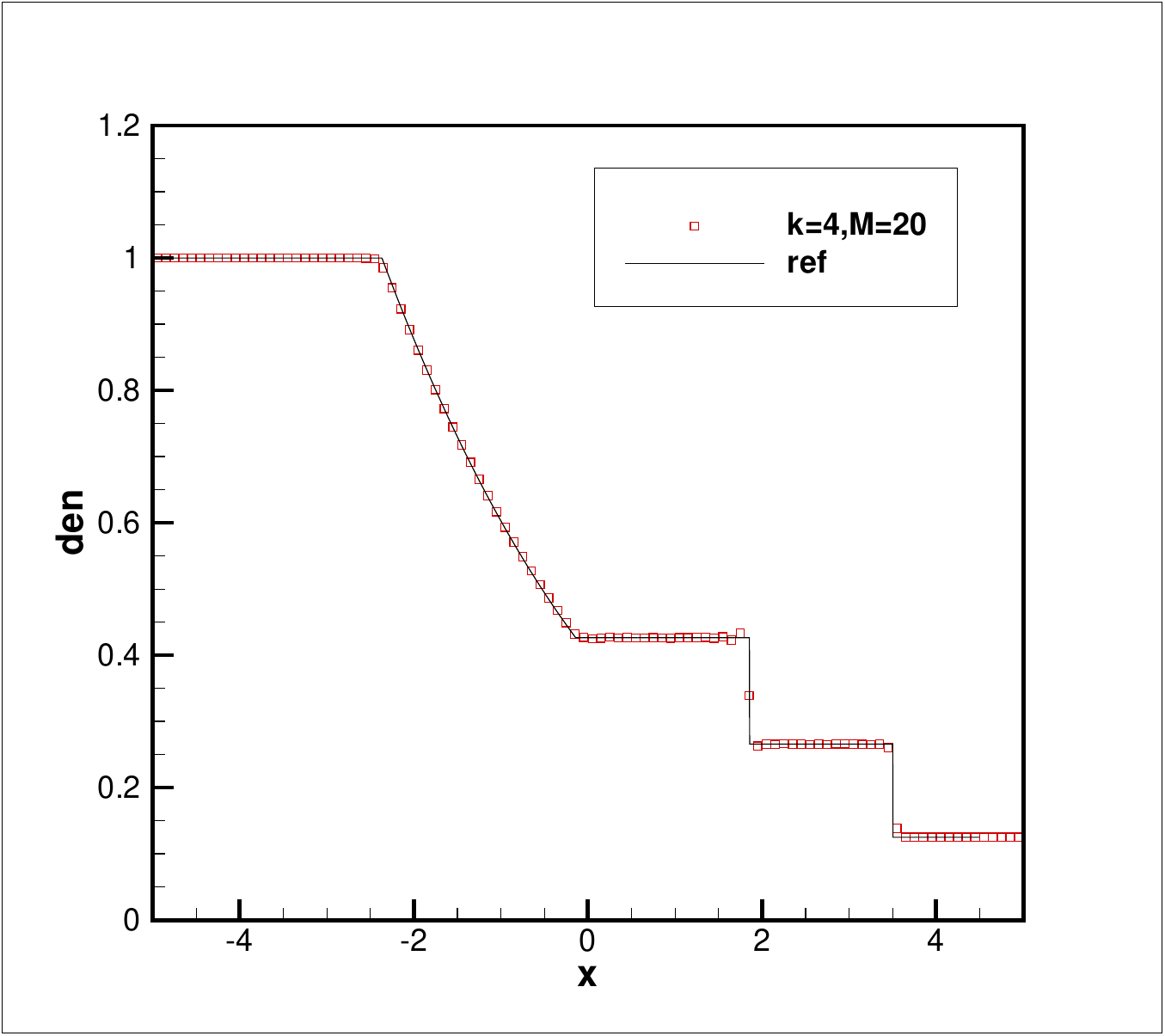}
    \includegraphics[width=0.3\textwidth]{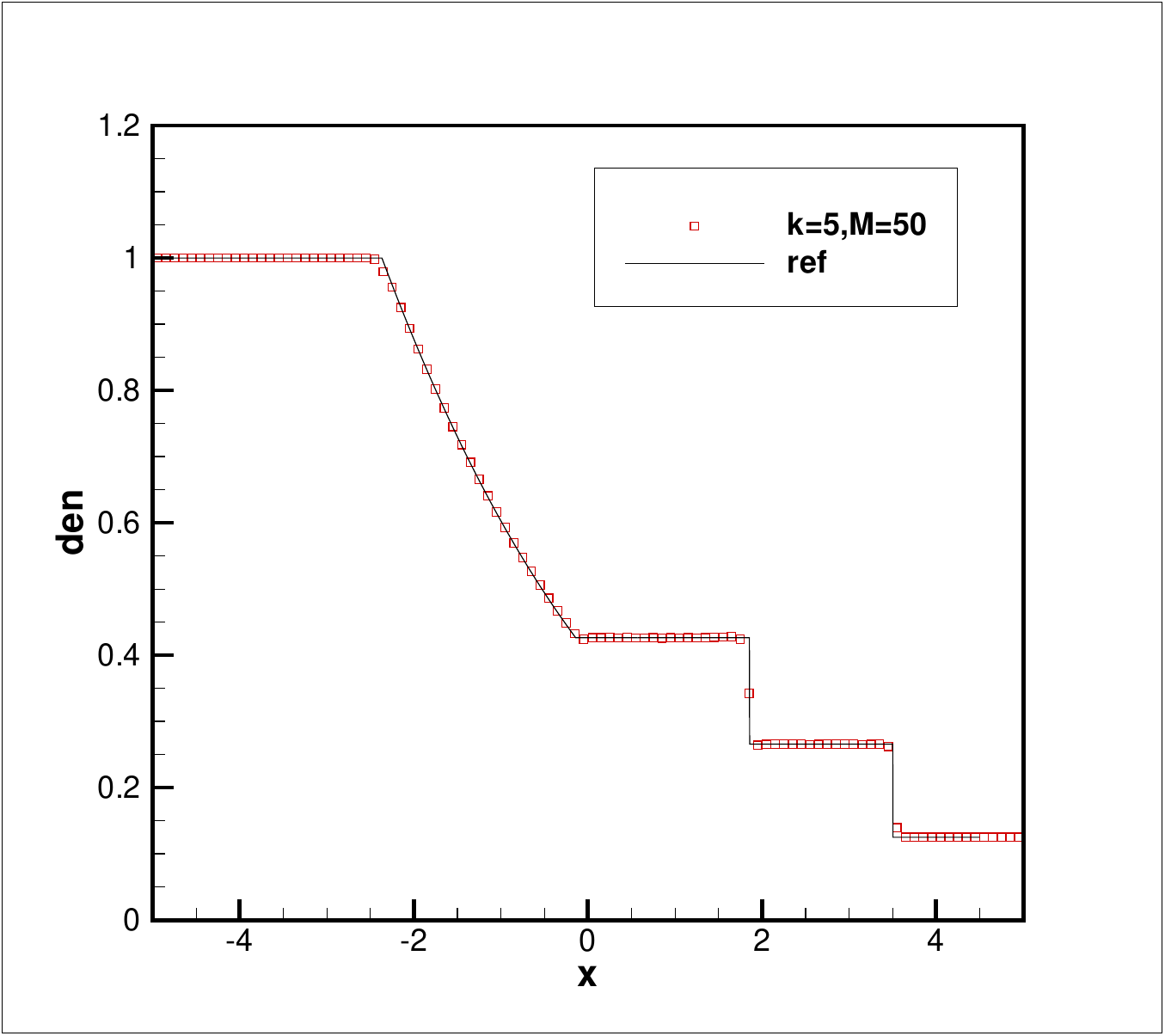}
    \includegraphics[width=0.3\textwidth]{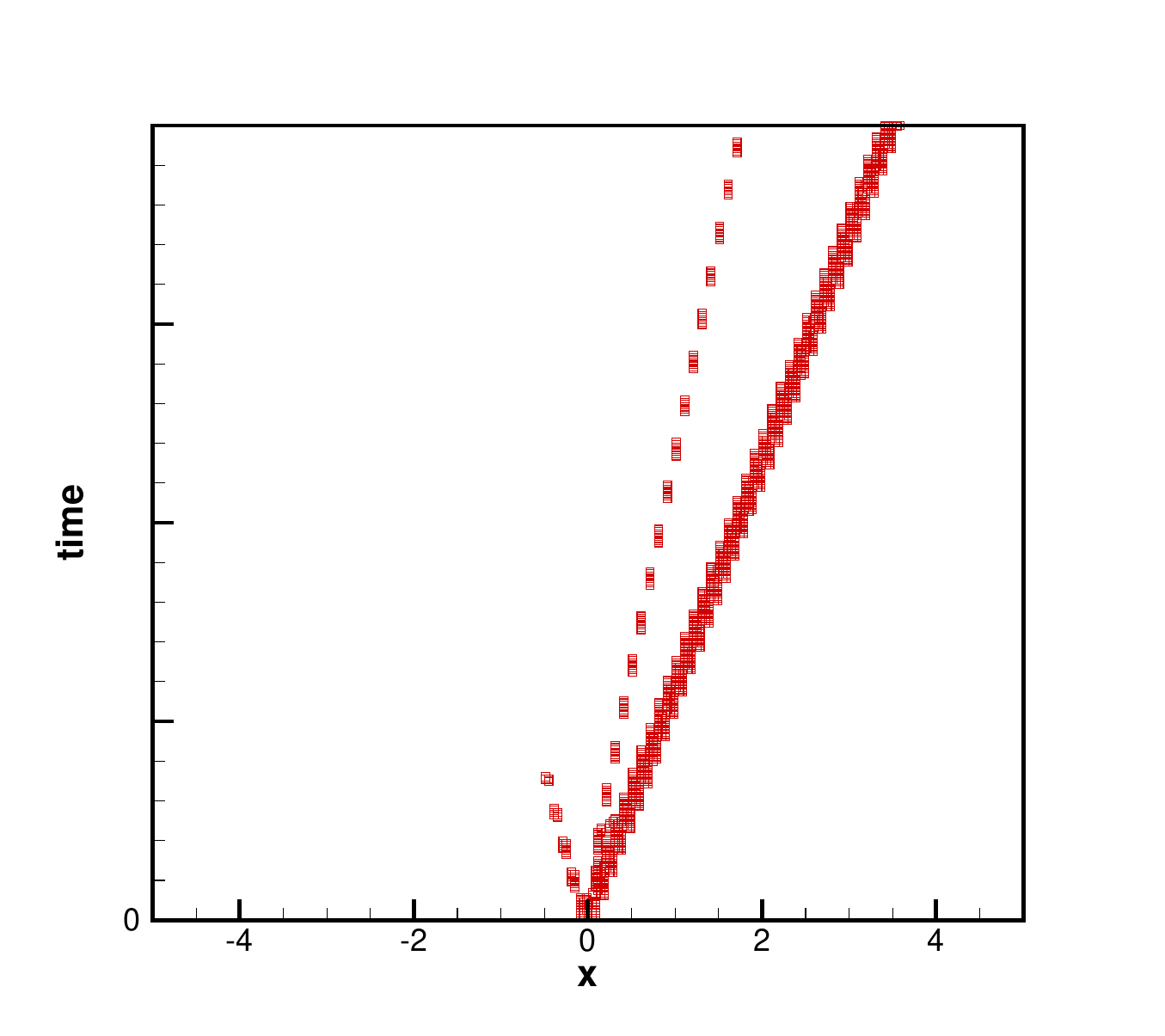}
    \includegraphics[width=0.3\textwidth]{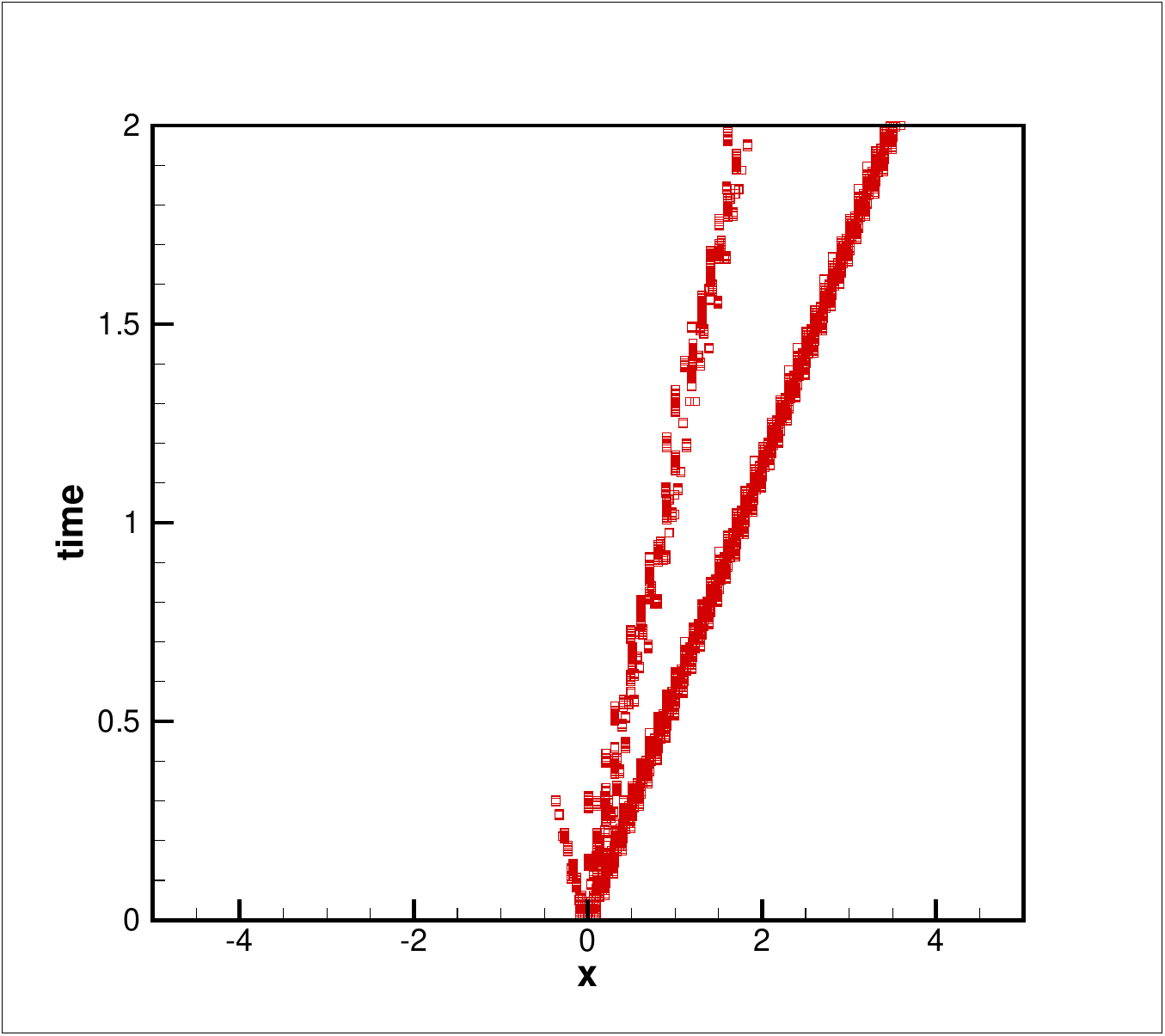}
    \includegraphics[width=0.3\textwidth]{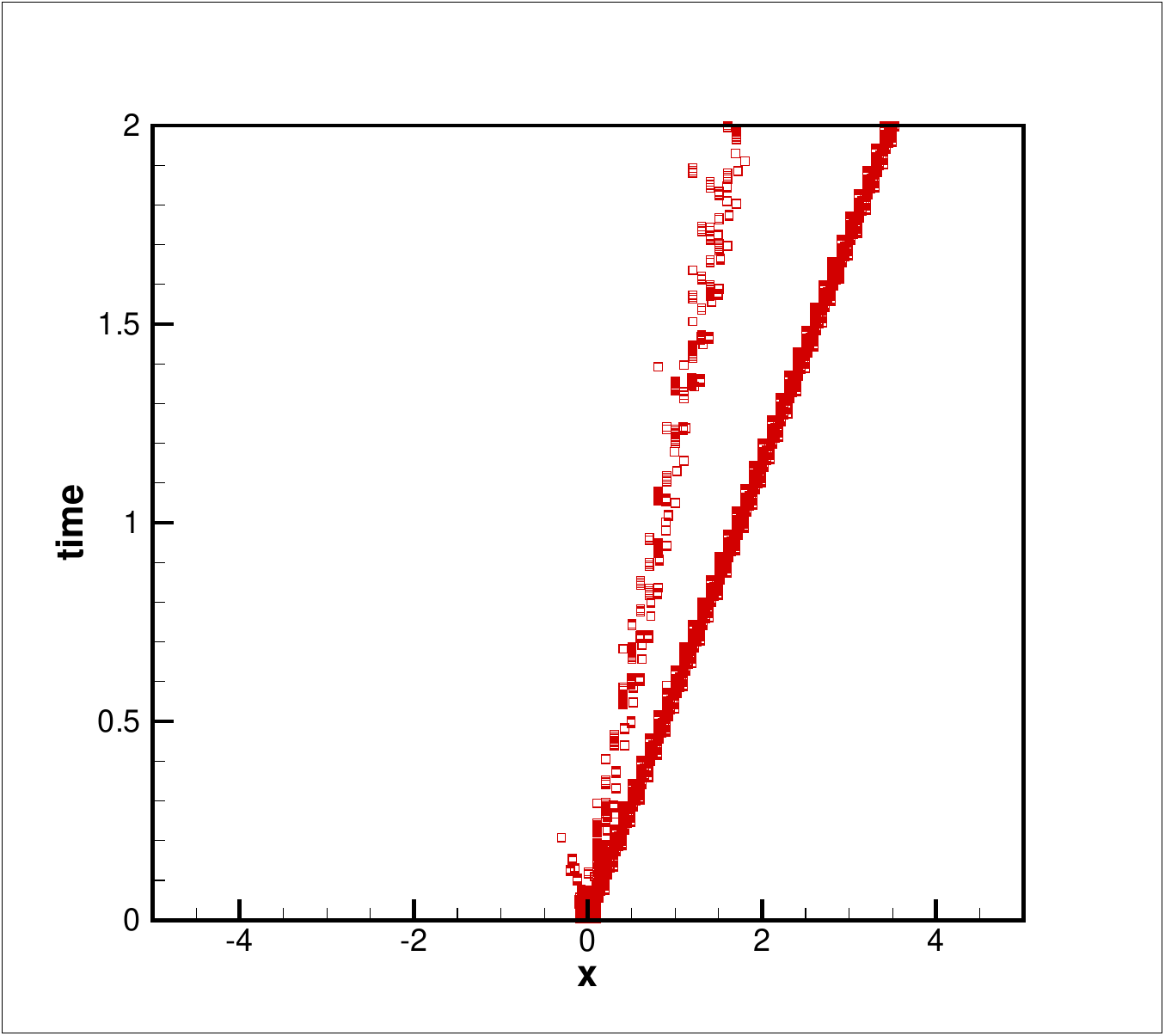}
     \caption{Example \ref{chap03:sod1d}: One-dimensional Sod shock tube problem:
     The density numerical solution obtained using 100 spectral volume elements in {\tt SV-cvMSWNO3}-\texttt{SV-cvMSWENO5} scheme with different $M$ values at $t=2$. The solid line and symbol ``$\square$'' represent the exact solution and {\tt SV-cvMSWENO} results, respectively.}
    \label{Fig:sod1d2}
 \end{figure}
\end{example}

\begin{example}[One-dimensional Lax shock tube problem]\label{chap03:lax1d}\rm
Considering the one-dimensional Lax shock tube problem, the computational domain is $\Omega=[-5,5]$, with the initial conditions specified as,
\begin{equation*}
(\rho,U,P)(x,0)=\left\{ \begin{array}{lc}
(0.445,0.698,3.528),& x<0,\\
(0.5,0,0.571),&x>0, \end{array}\right.
\end{equation*}
The time is calculated to time $t=1.3$. Figure \ref{Fig:lax1d} illustrates the numerical solution of density obtained using 100 spectral volume elements for each order scheme when $M=0.01$, along with the corresponding labeled troubled cells.
Similar to the result of Example \ref{chap03:sod1d}, the high-order schemes fail to demonstrate their high-resolution advantages due to an excessive number of troubled cells.
Additionally, the third order scheme exhibits instability at the contact discontinuity.
We subsequently tested a fully restricted scheme, where all cells were labeled as troubled cells.
As shown in Figure \ref{Fig:lax1d2}, the third order scheme still produces some overshoot and undershoot at the contact discontinuity.

 \begin{figure}[htbp]
    \centering
    \includegraphics[width=0.3\textwidth]{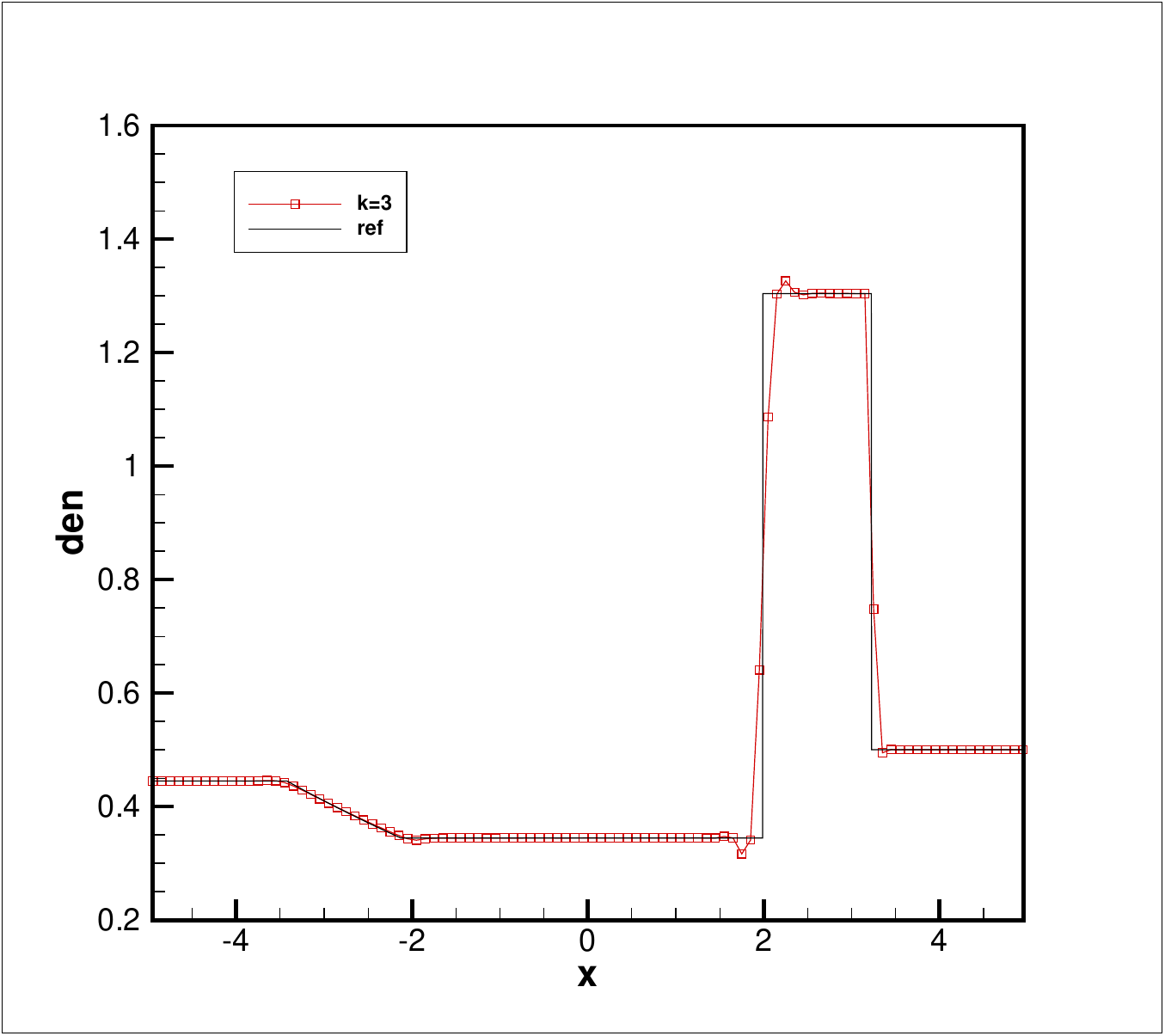}
    \includegraphics[width=0.3\textwidth]{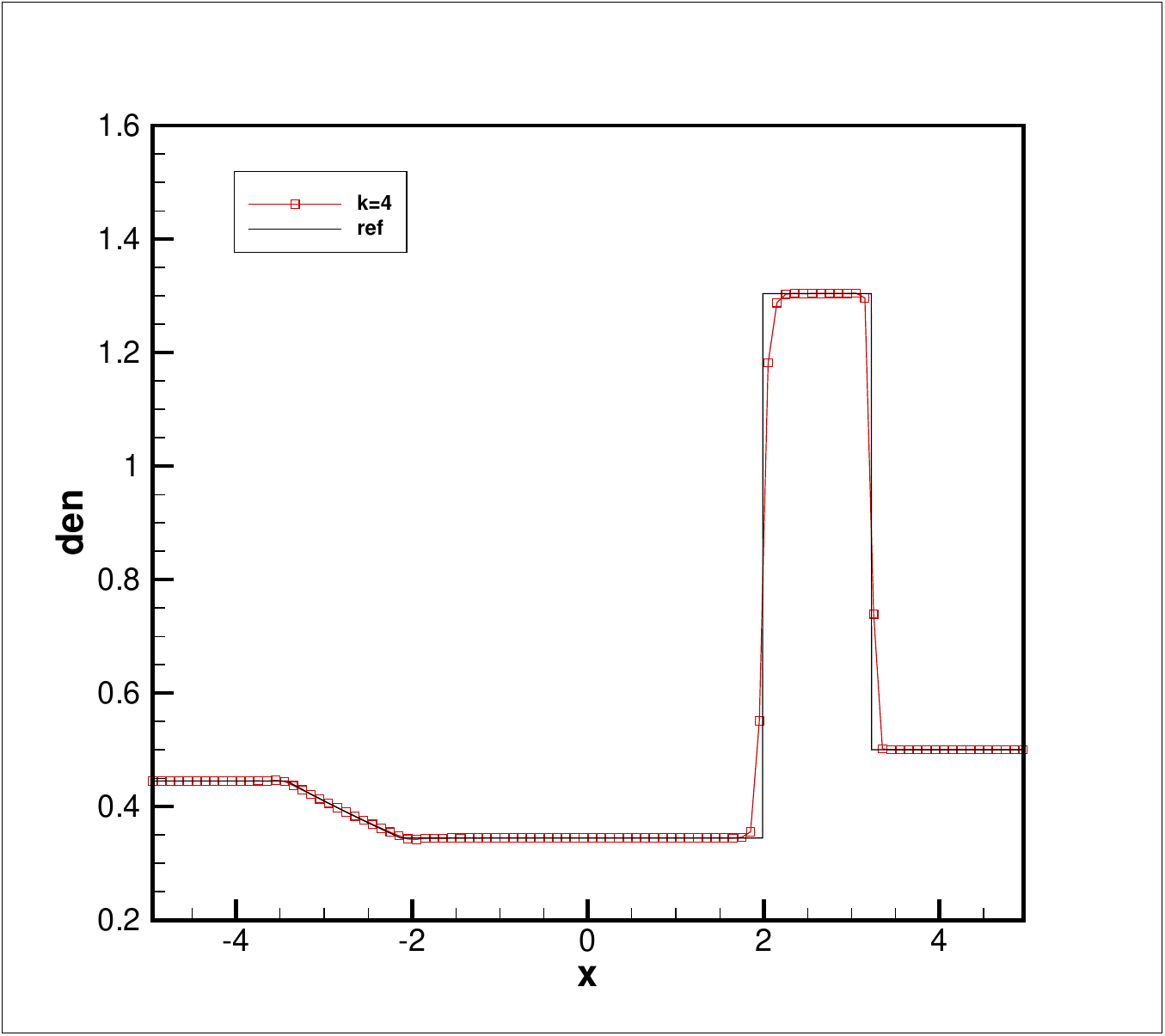}
    \includegraphics[width=0.3\textwidth]{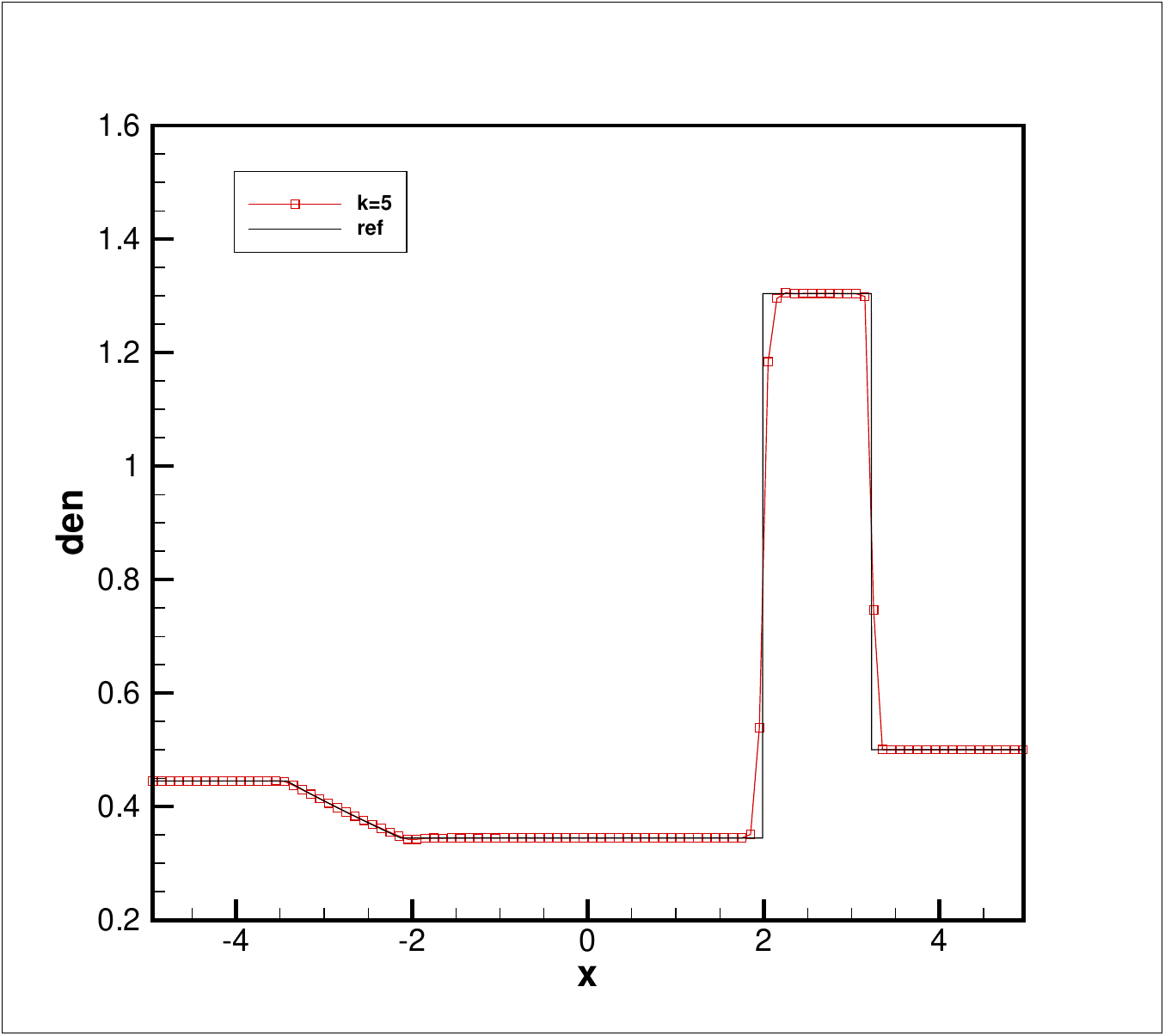}
    \includegraphics[width=0.3\textwidth]{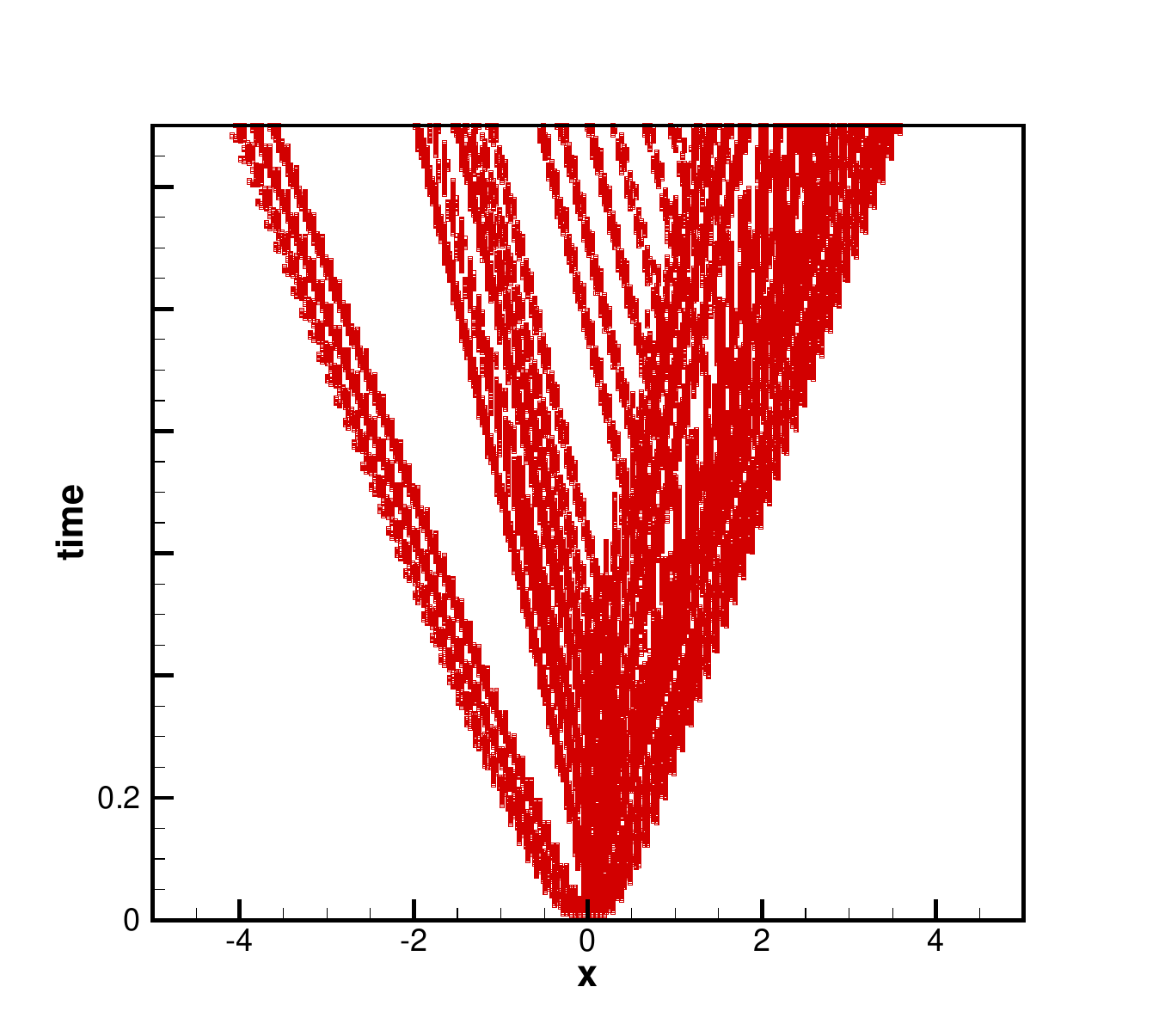}
    \includegraphics[width=0.3\textwidth]{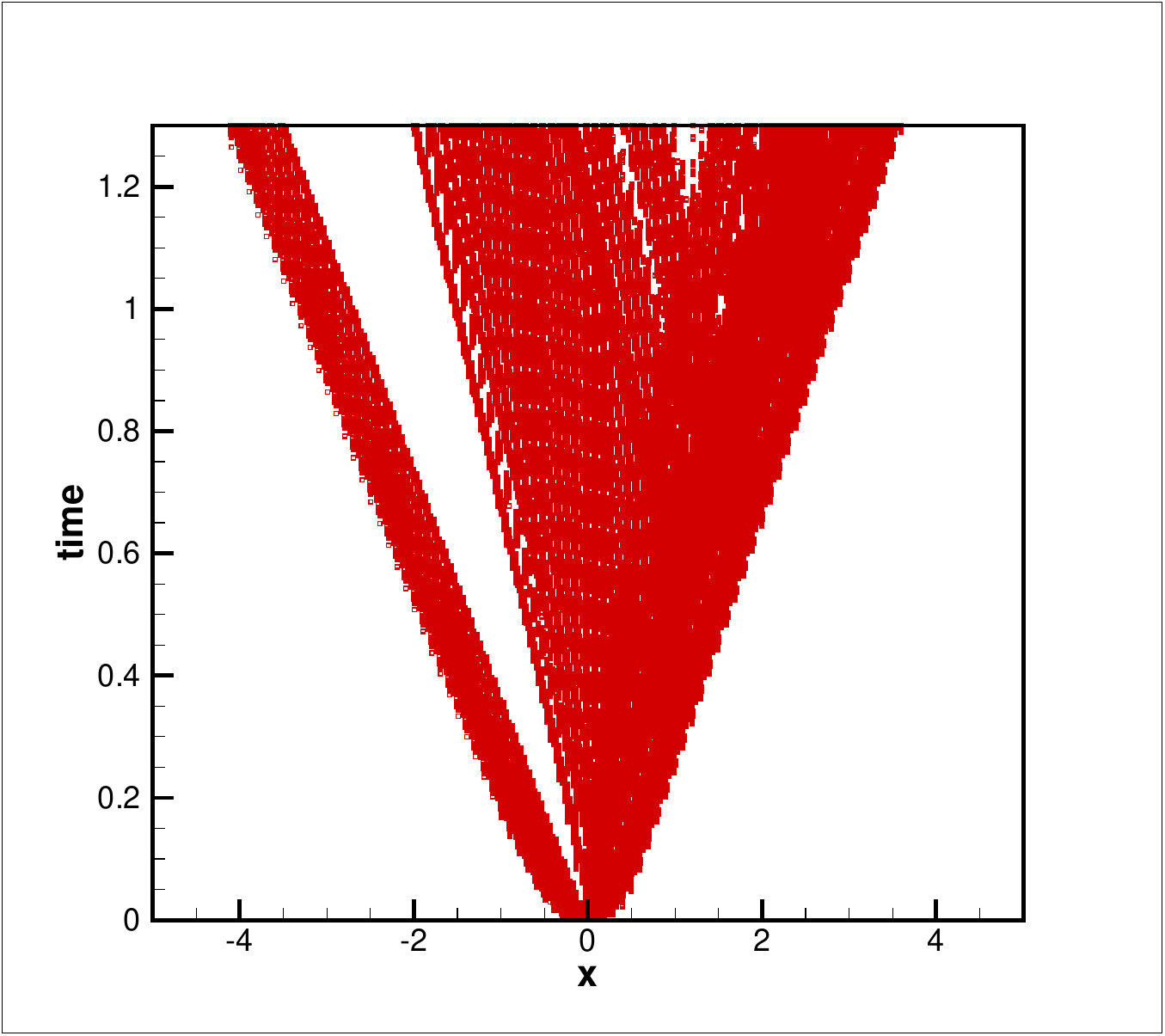}
    \includegraphics[width=0.3\textwidth]{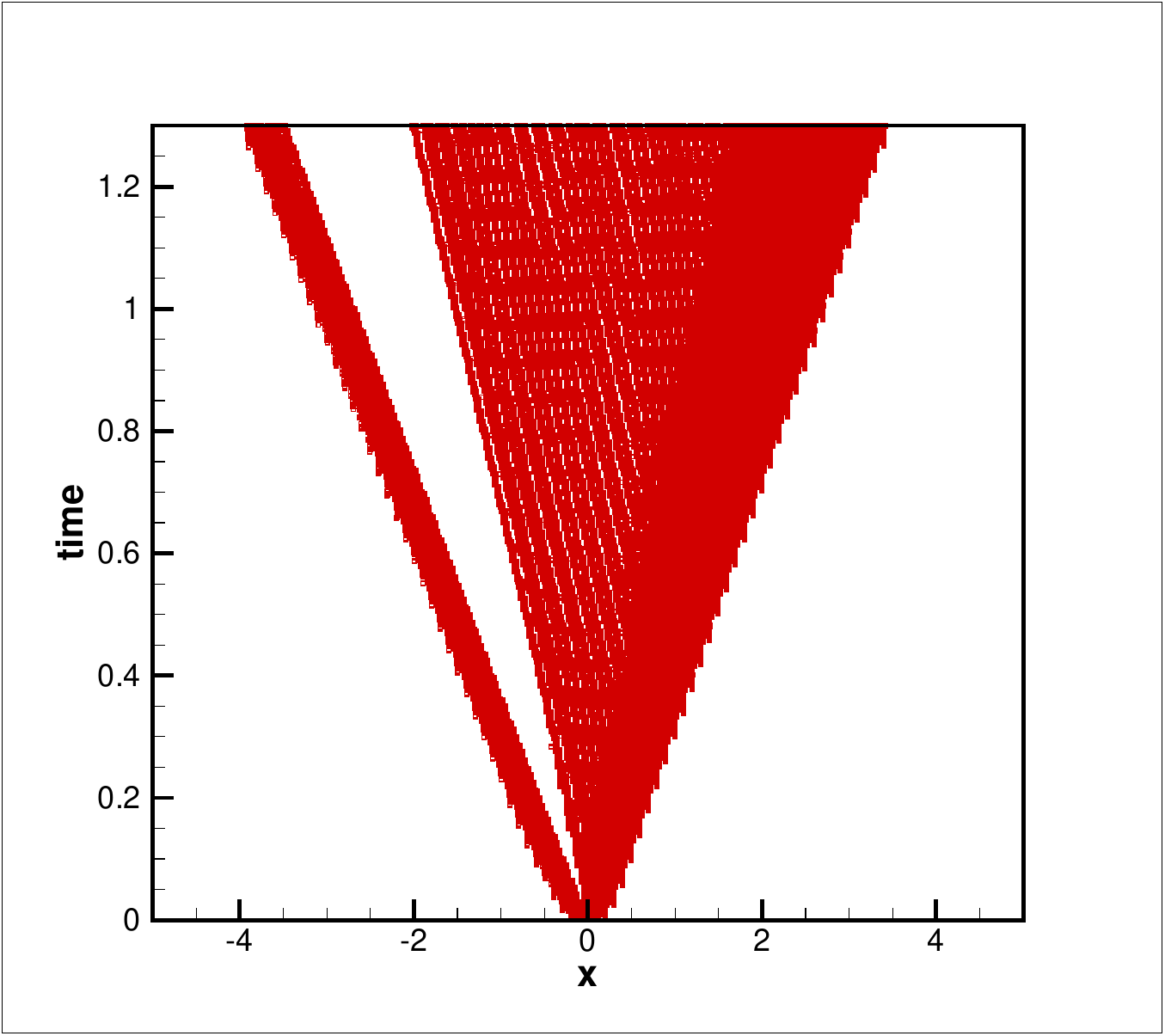}
     \caption{Example \ref{chap03:lax1d}: One-dimensional Lax shock tube problem:
     The density numerical solution obtained using 100 spectral volume elements in {\tt SV-cvMSWNO3}-\texttt{SV-cvMSWENO5} scheme when $M=0.01$ at $t=1.3$.  The solid line and symbol ``$\square$'' represent the exact solution and {\tt SV-cvMSWENO} results, respectively.}
    \label{Fig:lax1d}
 \end{figure}

  \begin{figure}[htbp]
    \centering
    \includegraphics[width=0.3\textwidth]{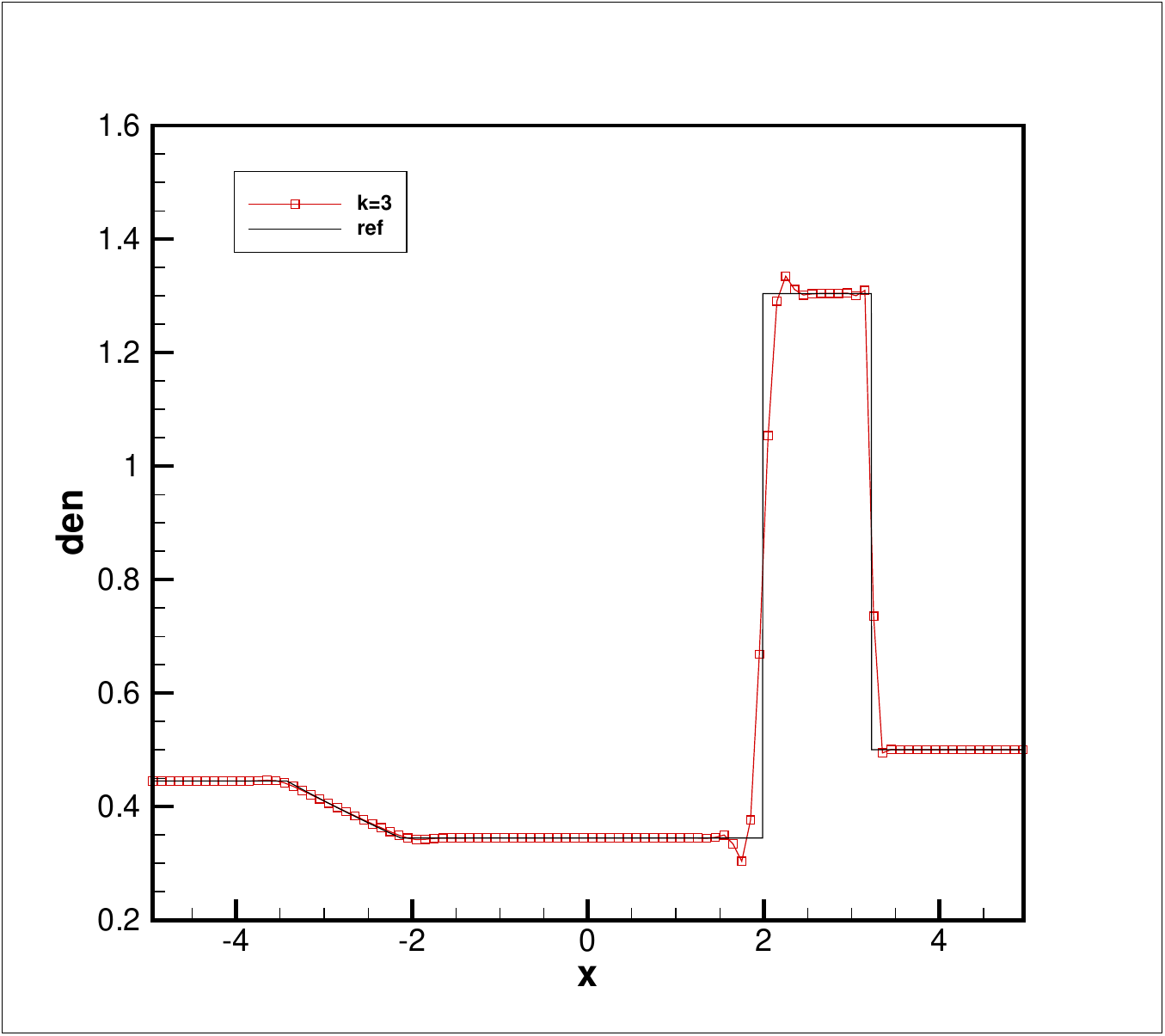}
    \includegraphics[width=0.3\textwidth]{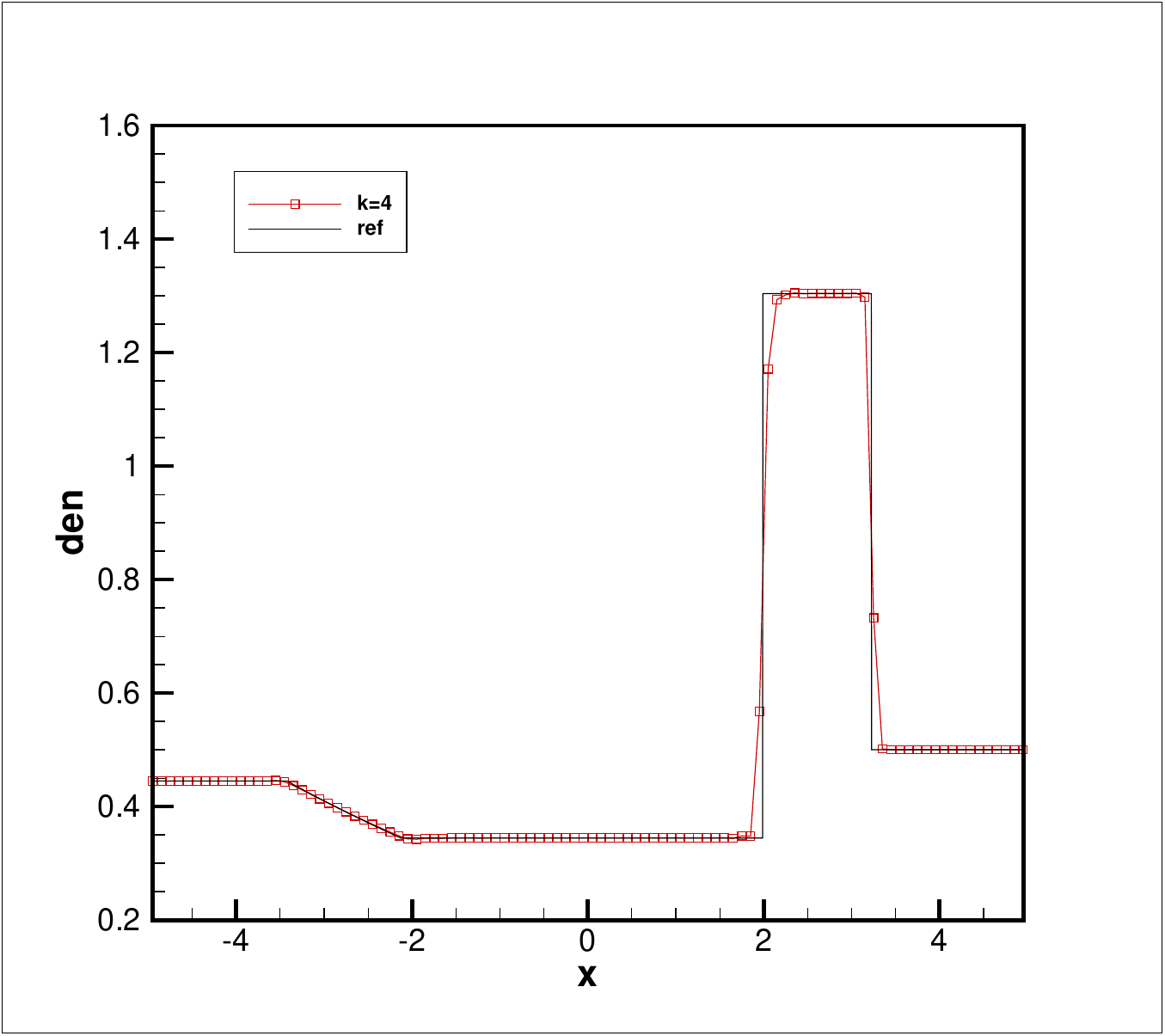}
    \includegraphics[width=0.3\textwidth]{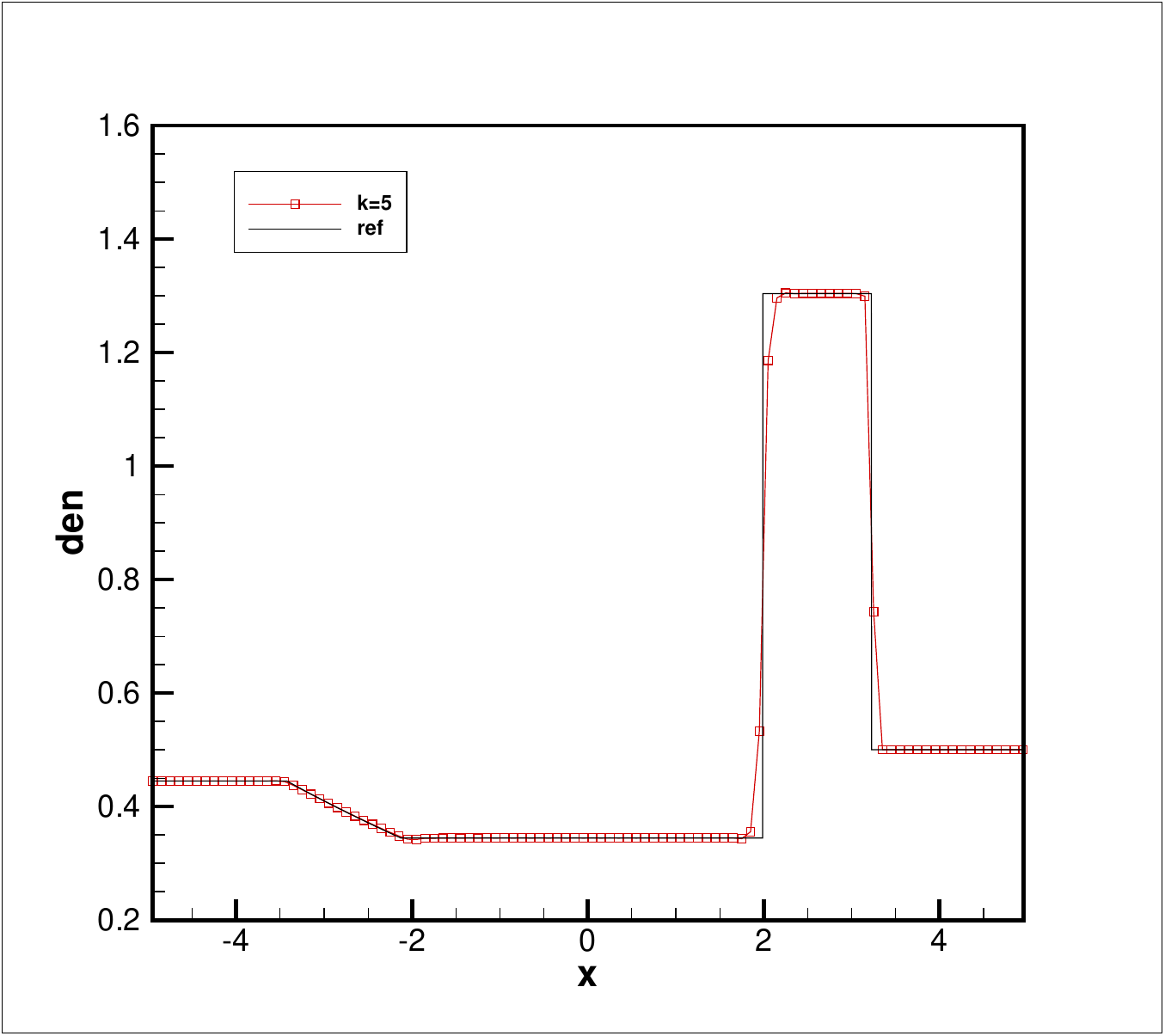}
     \caption{Example \ref{chap03:lax1d}: One-dimensional Lax shock tube problem:
     The density numerical solution obtained using 100 spectral volumes when fully restricted by {\tt SV-cvMSWNO3}-\texttt{SV-cvMSWENO5} scheme at $t=1.3$.  The solid line and symbol ``$\square$'' represent the exact solution and {\tt SV-cvMSWENO} results, respectively.}
    \label{Fig:lax1d2}
 \end{figure}

\end{example}

\begin{example}[Shock and sine wave interaction problem]\label{chap03:shuosher}\rm
This example describes the interaction between a mach 3 shock wave and a sinusoidal density wave. The computational domain is defined as $\Omega=[-5,5]$, with the initial condition specified as follows:
\begin{equation*}
(\rho,U,P)(x,0)=\left\{ \begin{array}{lc}
 (3.857134,2.629369,10.33333),& x<-4,\\
(1+0.2\sin(5x),0,1),&x>-4. \end{array}\right.
\end{equation*}

A right-lateral shock wave with a mach number of 3 propagates into a sinusoidal density basin and interacts with it. Figure \ref{Fig:shuosher1d} shows the numerical solution of density and the corresponding labeled troubled cells obtained using 180 spectral volume elements for each order scheme at $t=1.8$ when $M=0.01$. It is evident that the shock wave's interaction with the sine wave generates a complex yet smooth wave structure on the left side of the shock front. At this moment, the numerical dissipation is significant, leading to a higher number of labeled troubled cells. Figure \ref{Fig:shuosher1d2} illustrates the numerical solution of density along with the corresponding troubled cells when  $M=300$. At this moment, the number of troubled cells is minimal and aligns well with the wave structure, indicating that each order scheme achieves high resolution.

 \begin{figure}[htbp]
    \centering
    \includegraphics[width=0.3\textwidth]{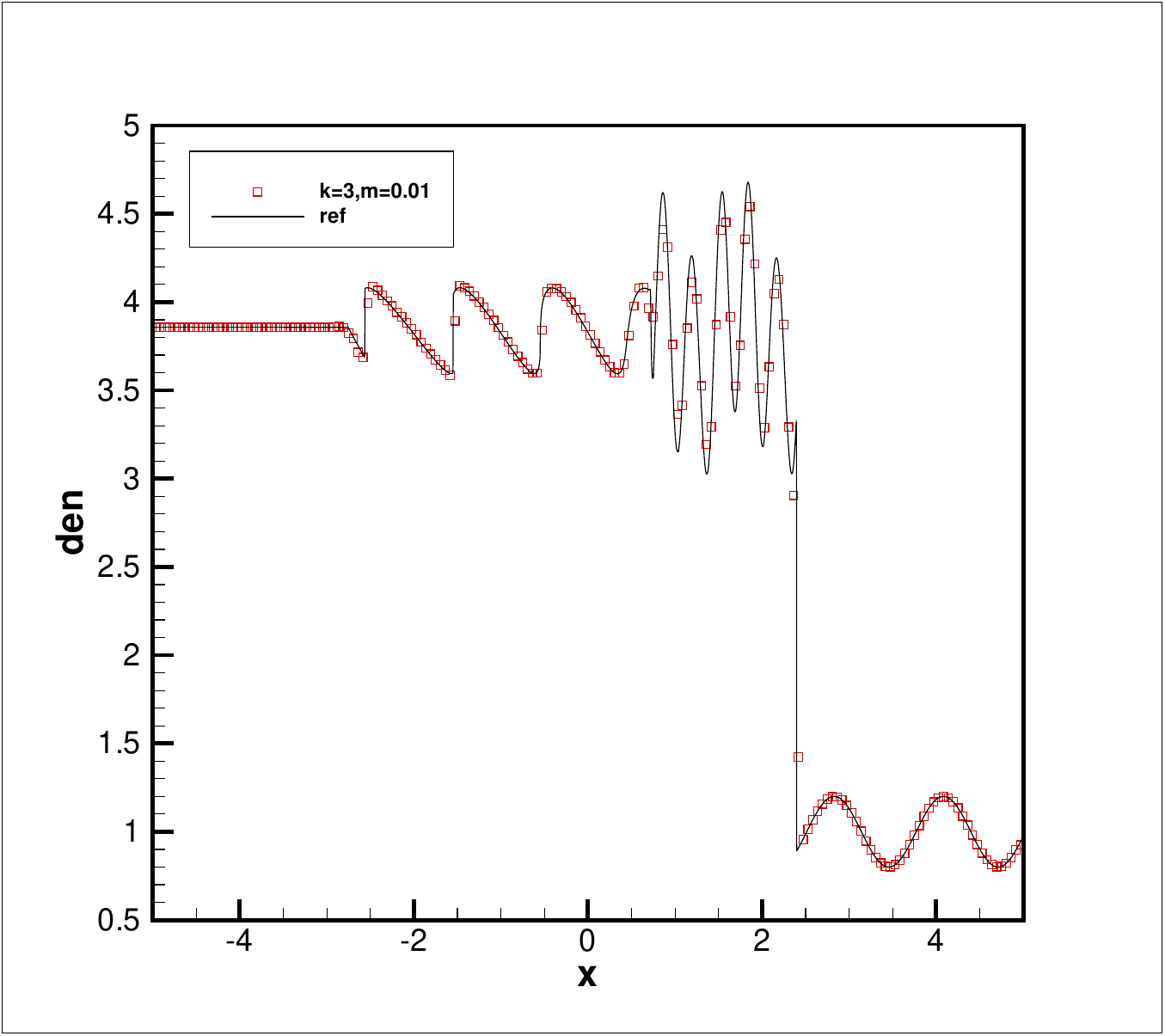}
    \includegraphics[width=0.3\textwidth]{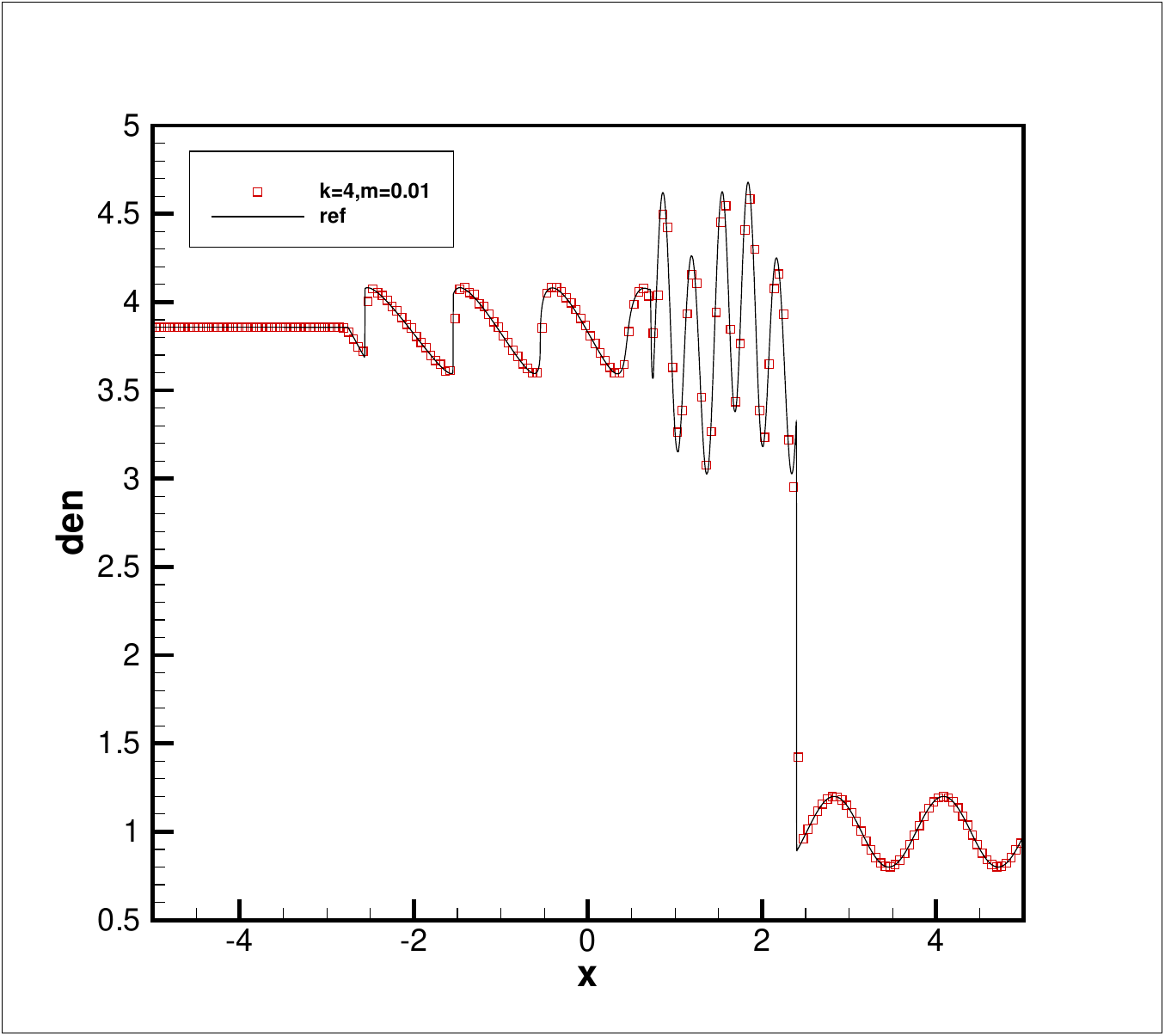}
    \includegraphics[width=0.3\textwidth]{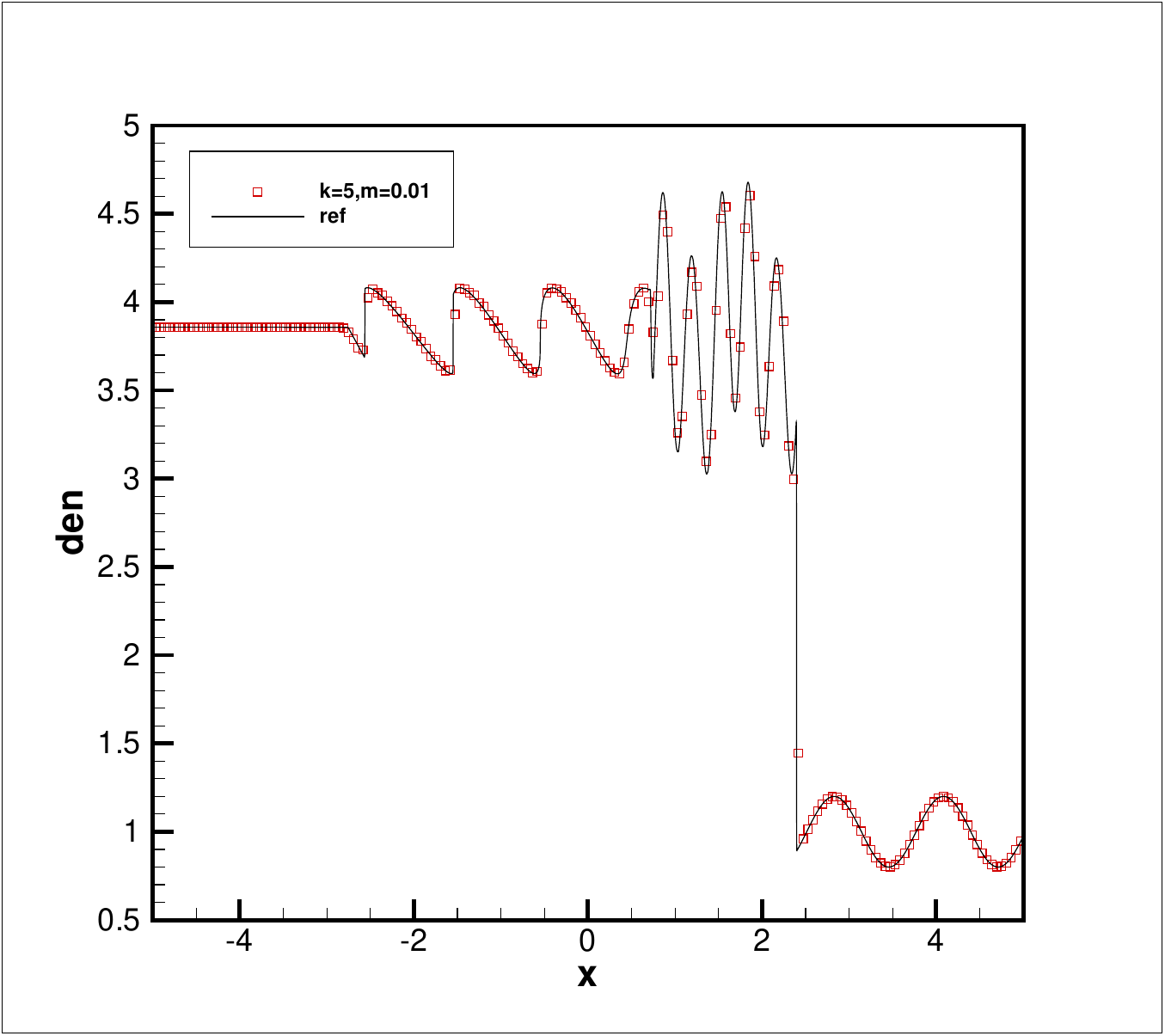}
    \includegraphics[width=0.3\textwidth]{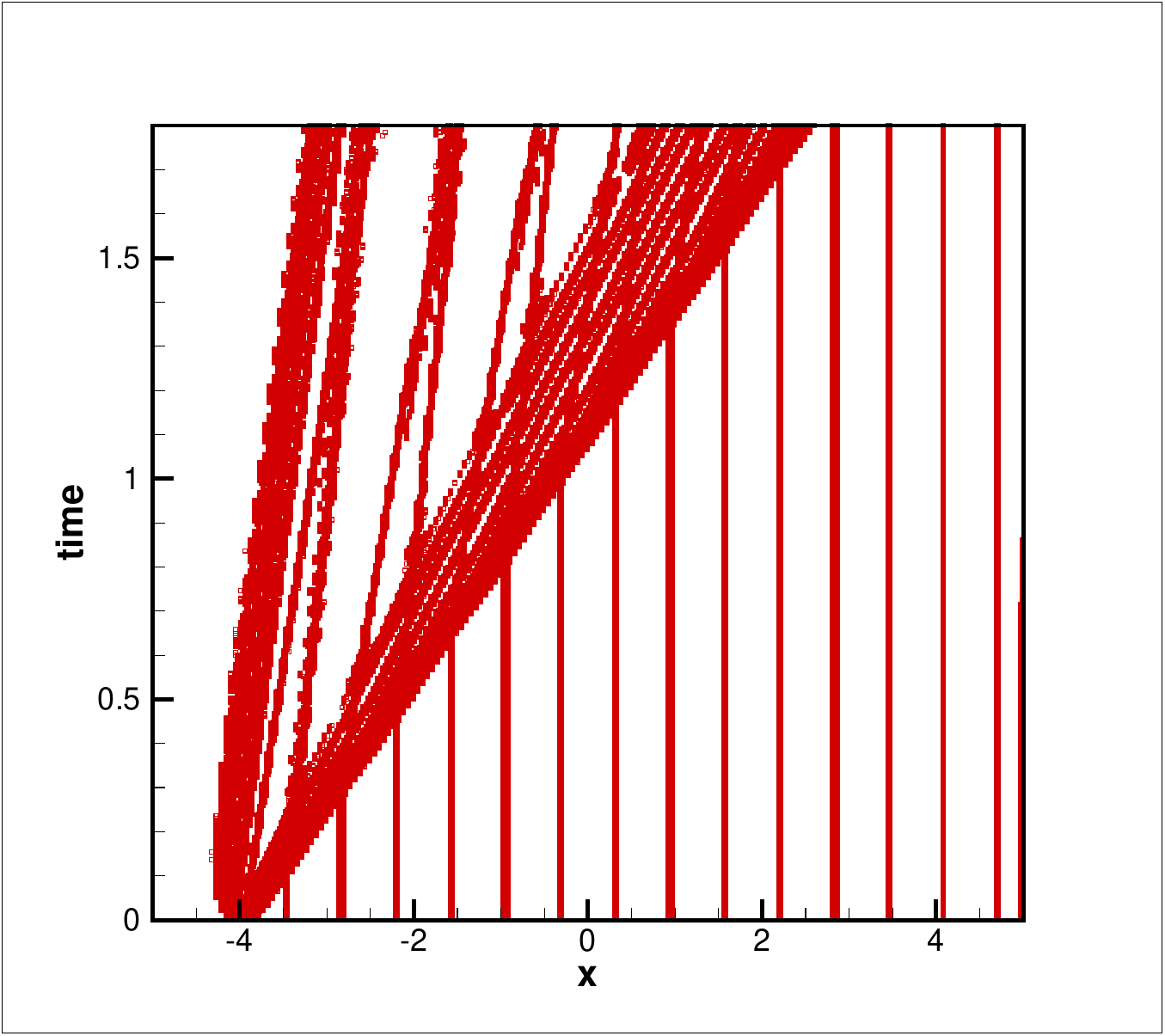}
    \includegraphics[width=0.3\textwidth]{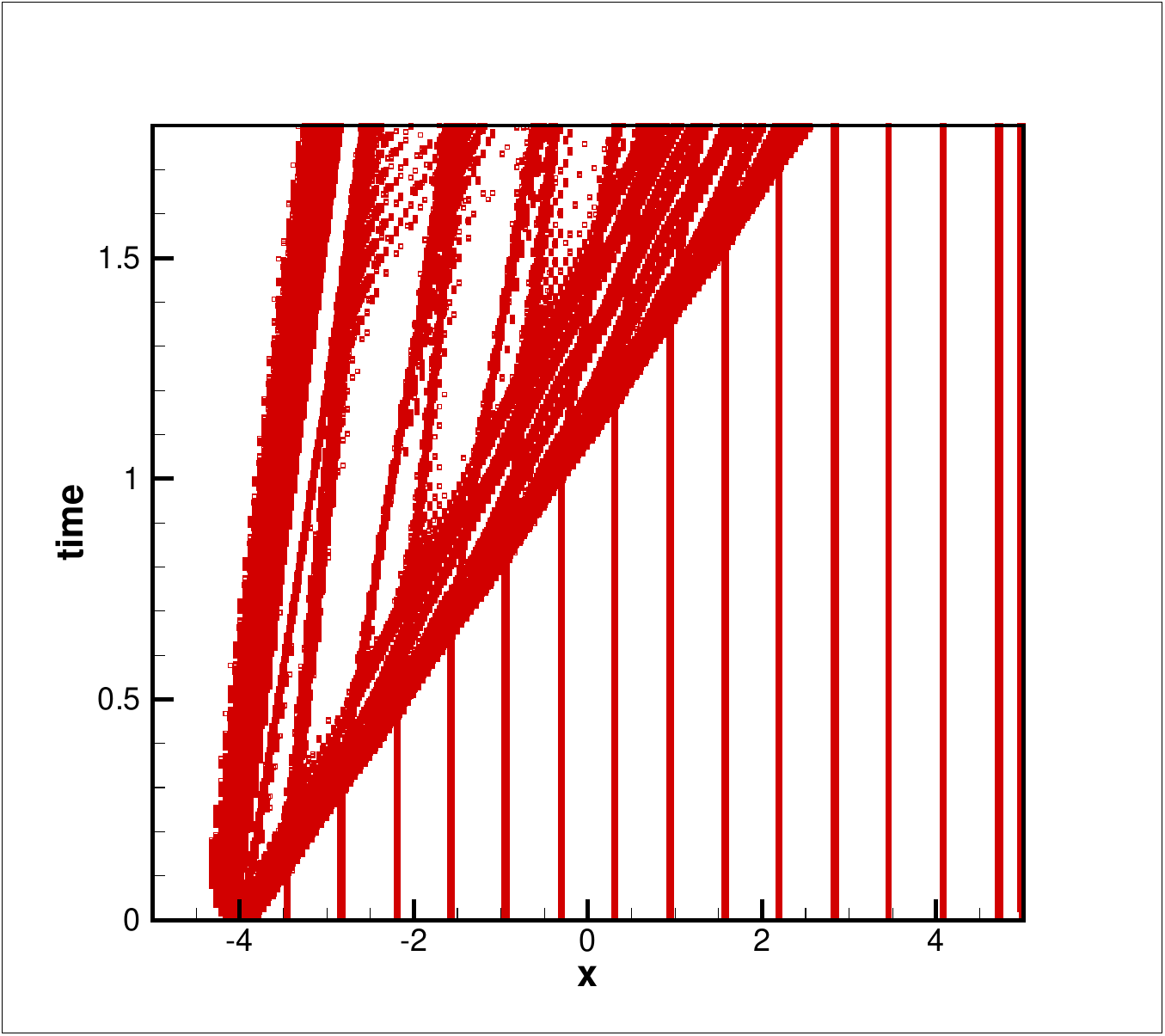}
    \includegraphics[width=0.3\textwidth]{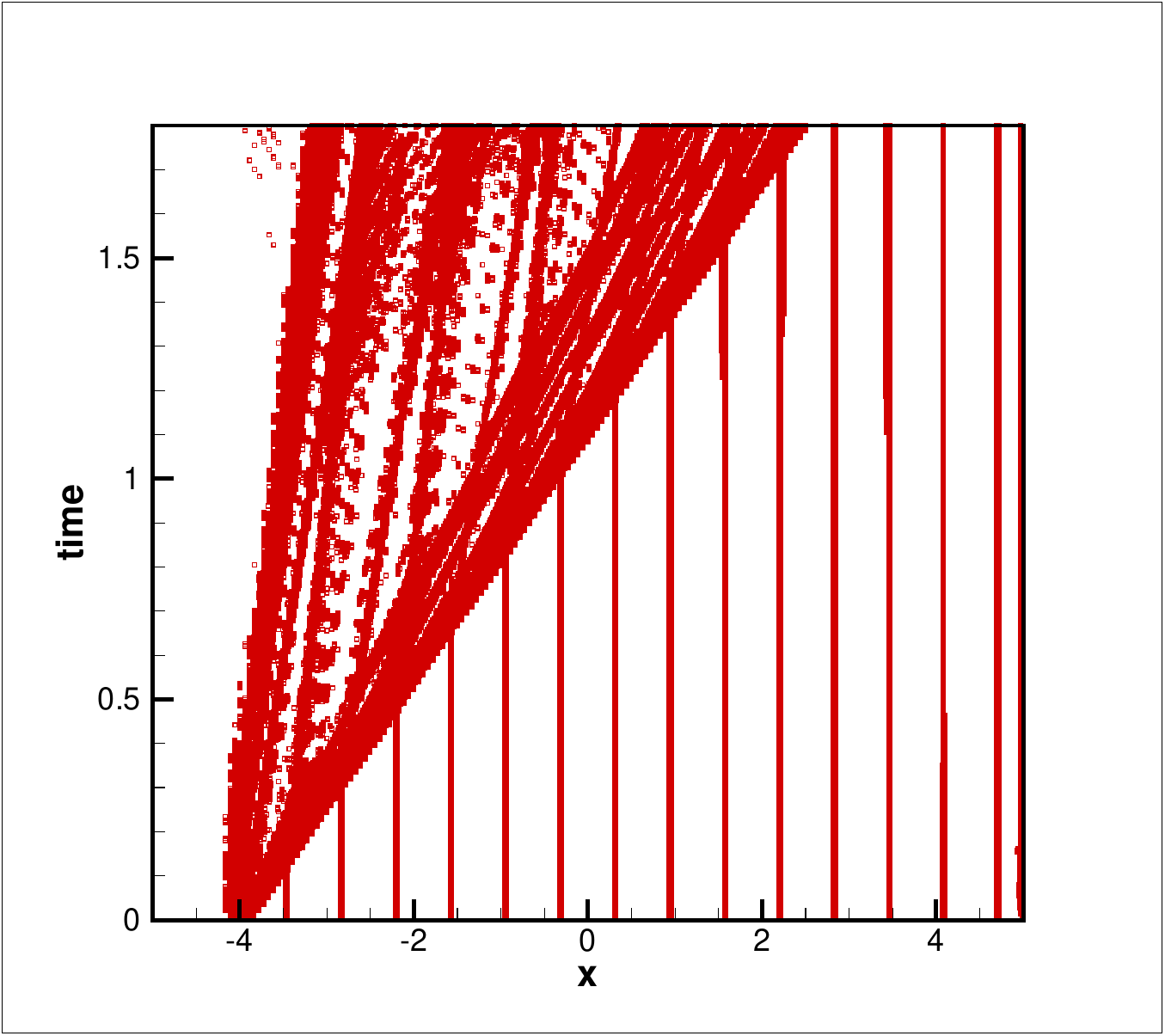}
     \caption{Example \ref{chap03:shuosher}: Shock wave and sine wave interaction problem:
      The density numerical solution obtained using 180 spectral volumes in {\tt SV-cvMSWNO3}-\texttt{SV-cvMSWENO5} scheme when $M=0.01$ at $t=1.8$.  The solid line and symbol ``$\square$'' represent the exact solution and {\tt SV-cvMSWENO} results, respectively.}
    \label{Fig:shuosher1d}
 \end{figure}

 \begin{figure}[htbp]
    \centering
    \includegraphics[width=0.3\textwidth]{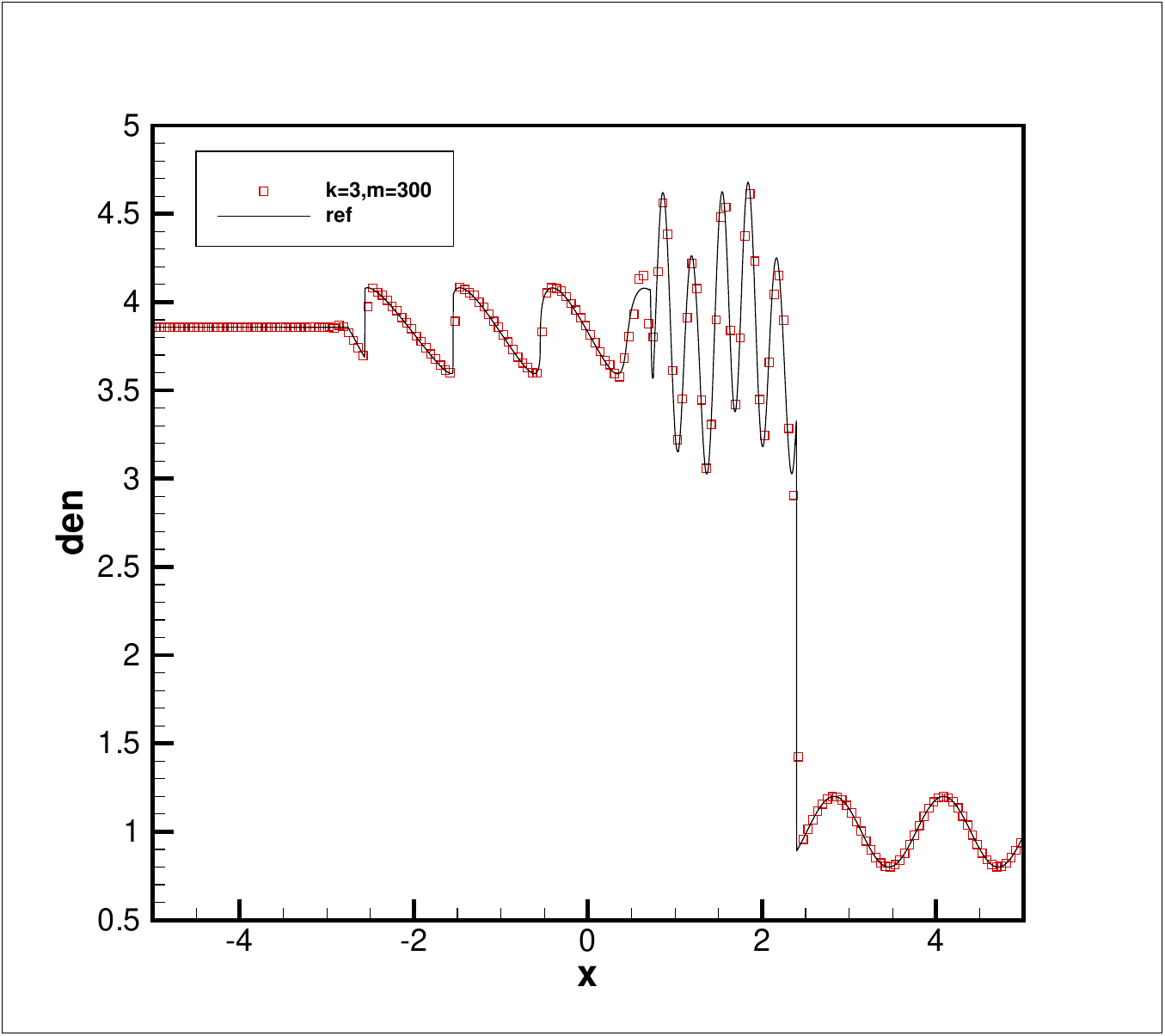}
    \includegraphics[width=0.3\textwidth]{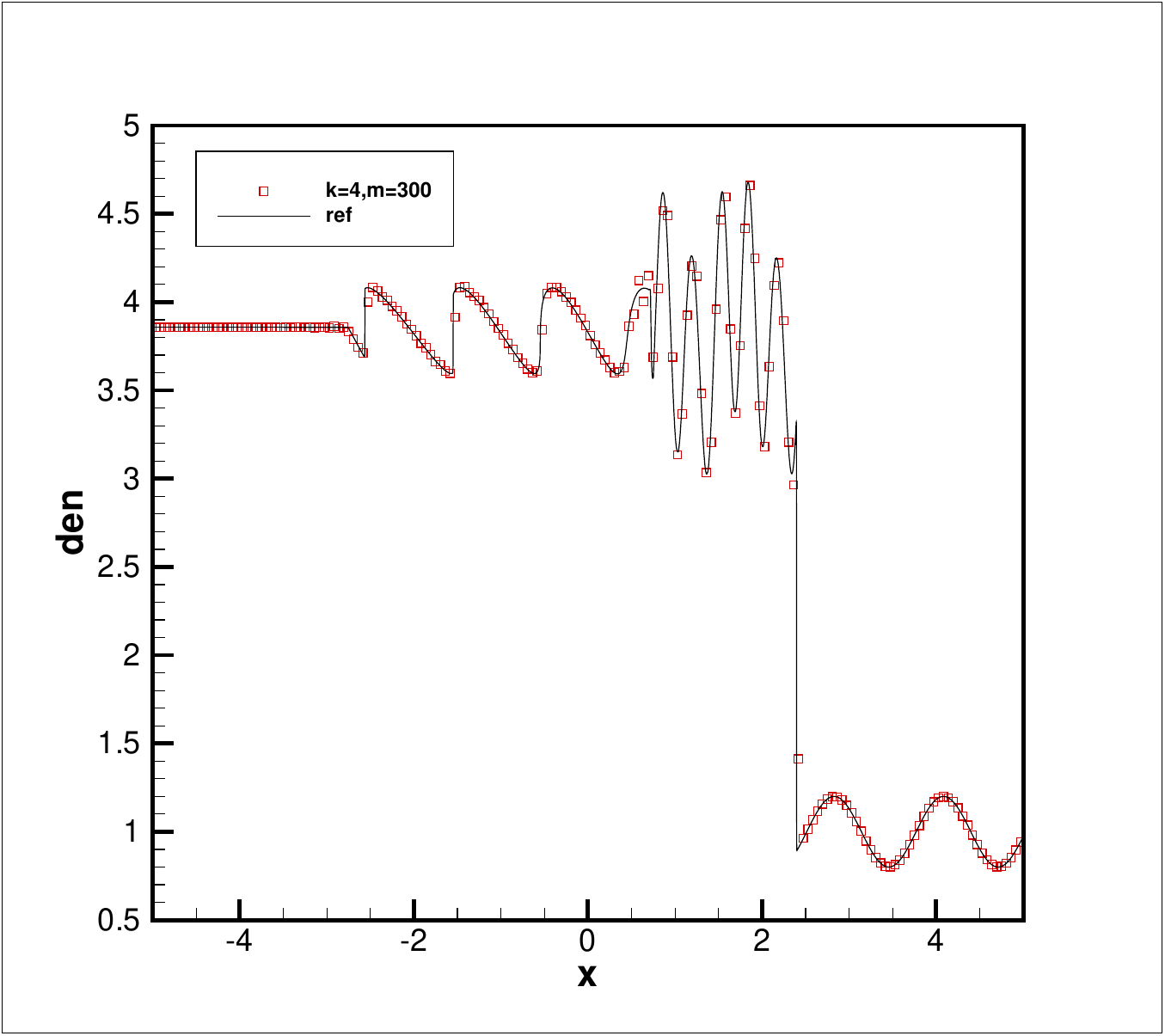}
    \includegraphics[width=0.3\textwidth]{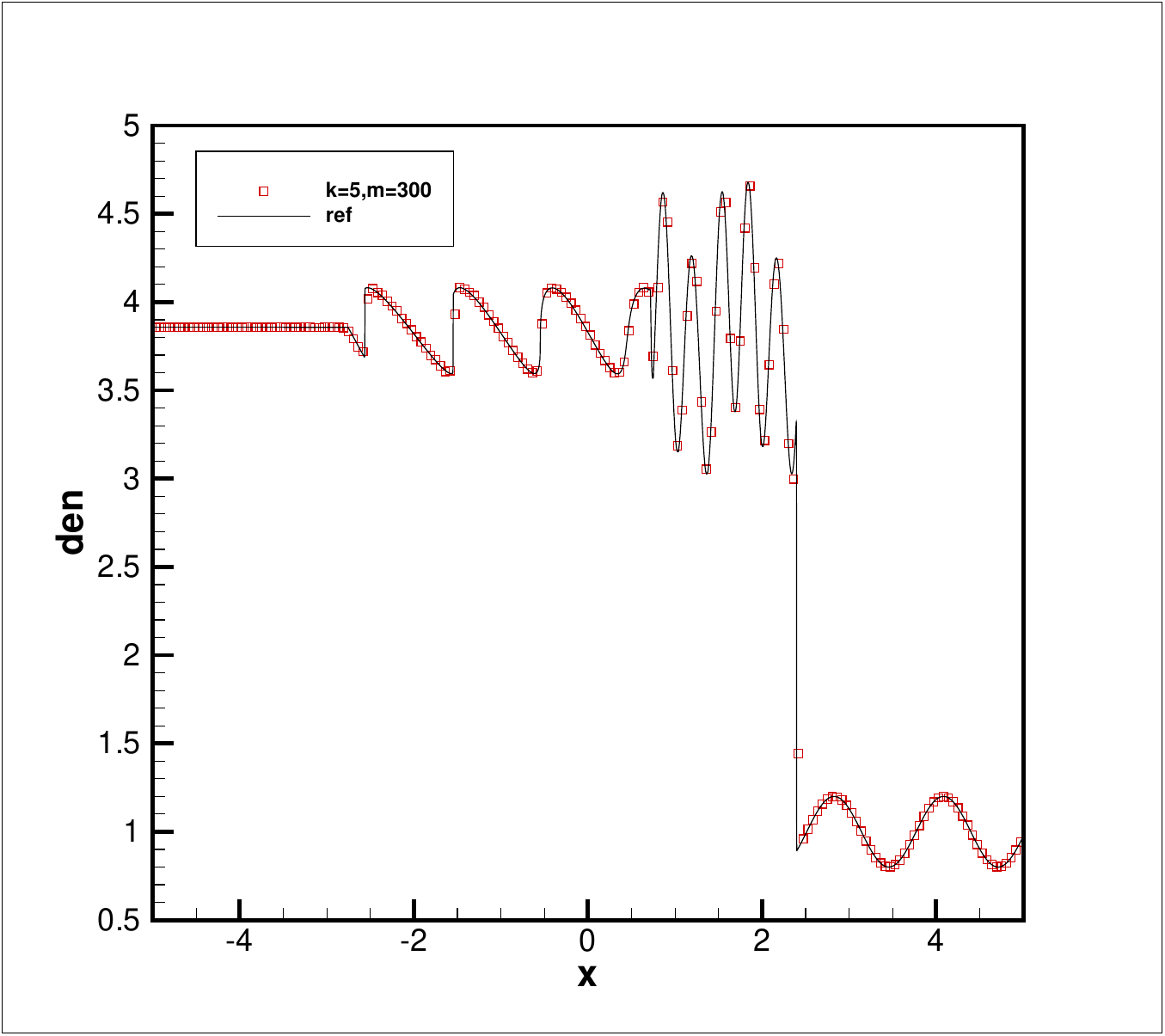}
    \includegraphics[width=0.3\textwidth]{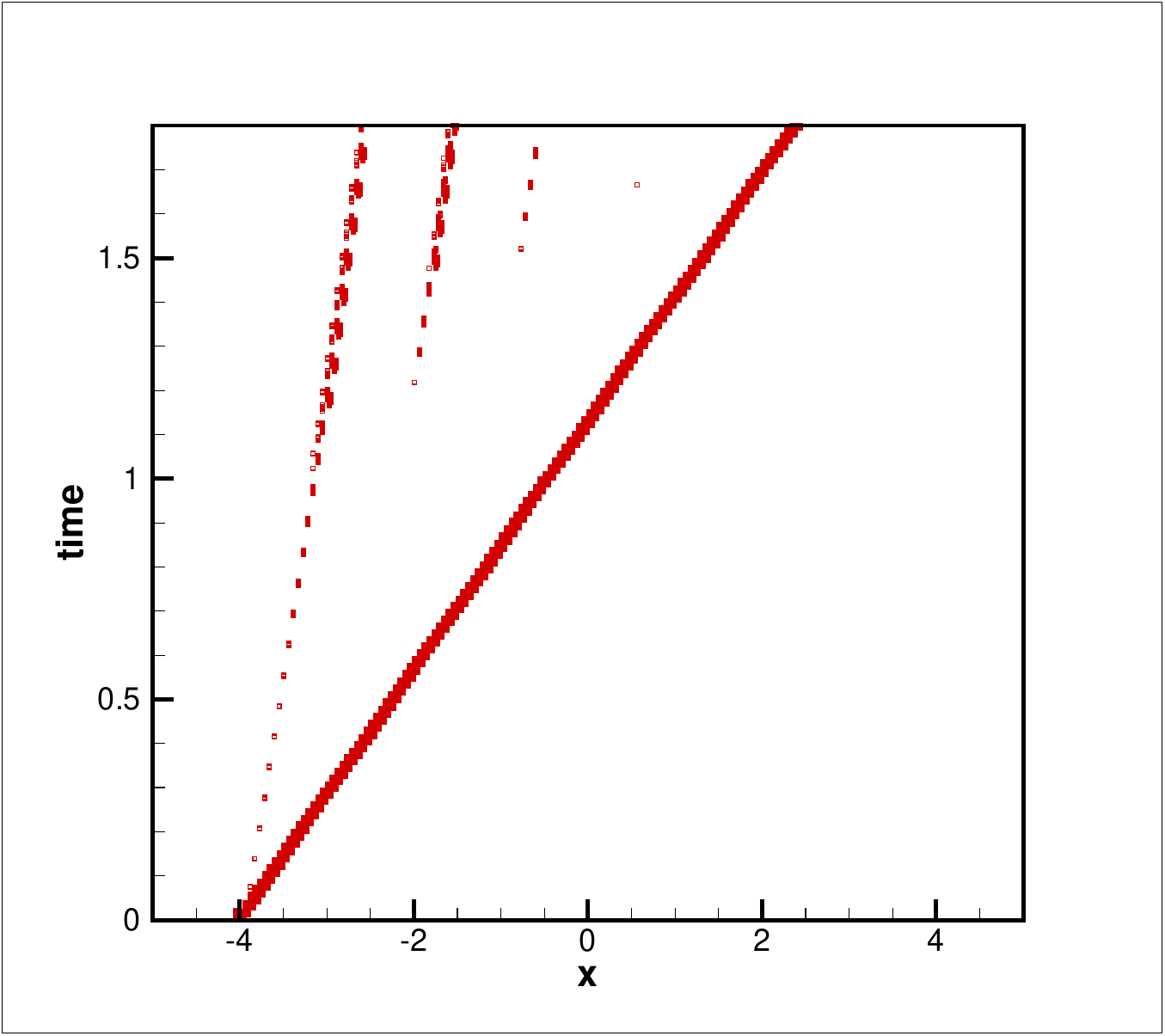}
    \includegraphics[width=0.3\textwidth]{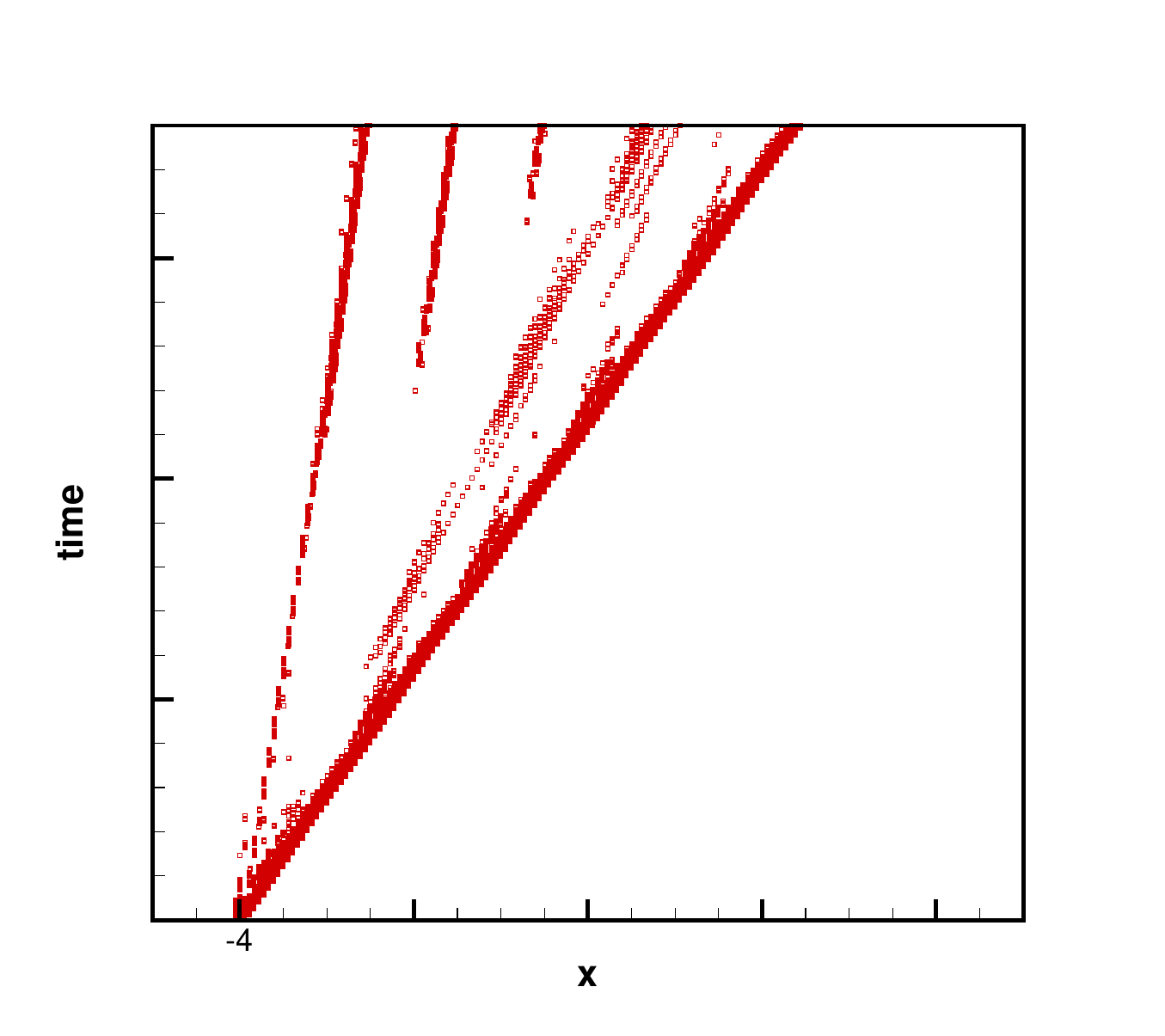}
    \includegraphics[width=0.3\textwidth]{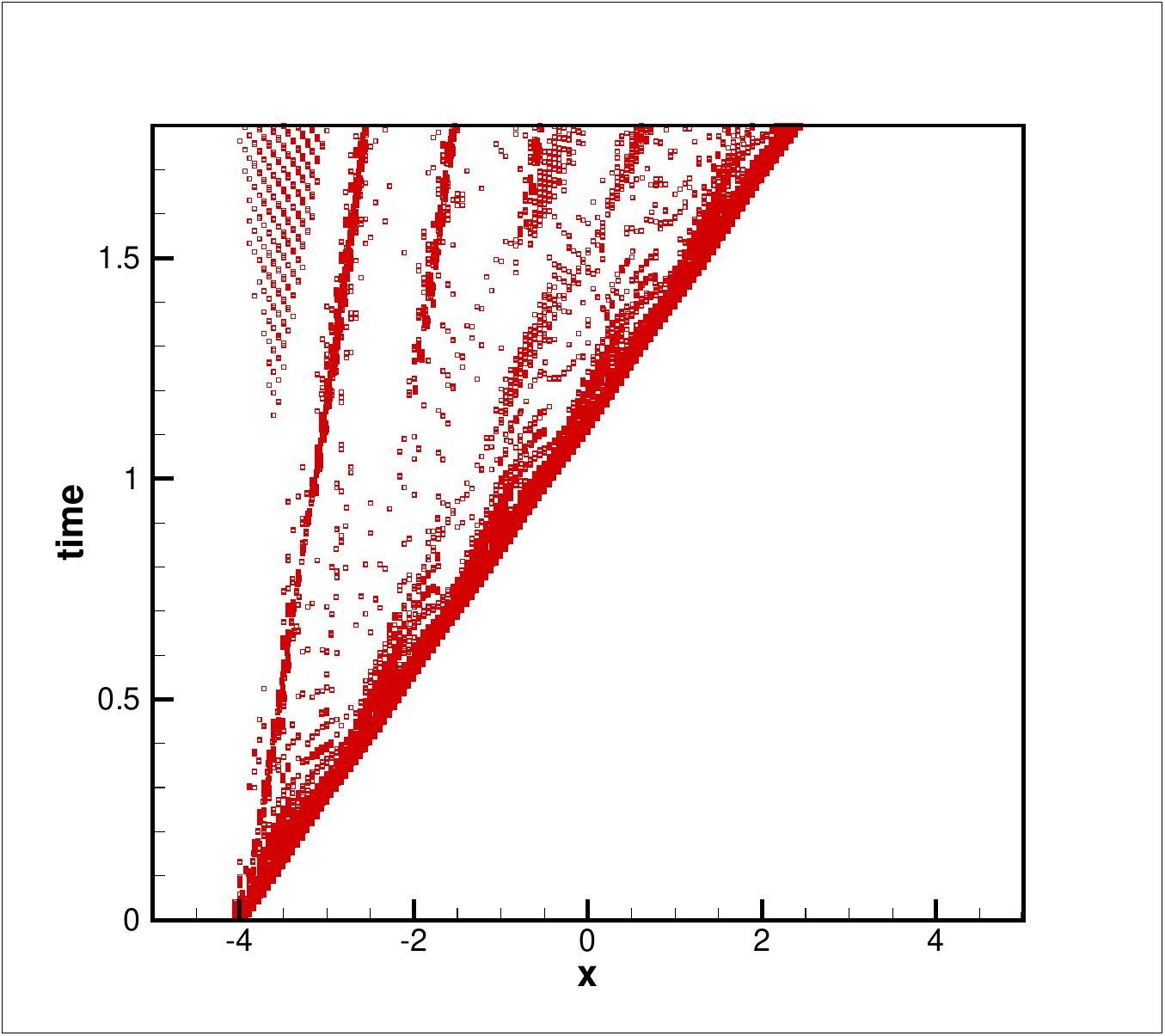}
     \caption{Example \ref{chap03:shuosher}: Shock wave and sine wave interaction problem:
     The density numerical solution obtained using 180 spectral volumes in {\tt SV-cvMSWNO3}-\texttt{SV-cvMSWENO5} scheme when $M=300$ at $t=1.8$.  The solid line and symbol ``$\square$'' represent the exact solution and {\tt SV-cvMSWENO} results, respectively.}
    \label{Fig:shuosher1d2}
 \end{figure}

\end{example}

\begin{example}[Explosive wave problem]\label{chap03:wood}\rm

The explosive wave problem was formulated by Woodward and Colella \cite{woodward1984ppm}. The computational domain $\Omega$ is defined as [0,1] with reflection boundaries at both endpoints of the region. The initial condition is
\begin{equation*}
(\rho,U,P)(x,0)=\left\{ \begin{array}{lc}
 (1,0,1000),& 0\le x\le 0.1,\\
 (1,0,0.01),& 0.1\le x\le 0.9,\\
 (1,0,100),& 0.9\le x\le 1.0. \end{array}\right.
\end{equation*}

Initially, the fluid is characterized by three distinct constant states. The pressure on both sides of the fluid is significantly high. Due to the high pressure environment, a left-traveling explosion wave and a right-traveling explosion wave are generated respectively. And the velocity of the right-traveling strong shock wave is notably faster. Over time, these two shock waves eventually meet and interact. Figure \ref{Fig:woodward1d} shows the numerical solution of density and the corresponding troubled cells obtained using 400 spectral volumes for each order scheme at $t=0.038$ when $M=0.01$. From the numerical results, it is evident that after the labeling and limiting in the reconstruction process, there are no significant numerical oscillations for any order of schemes. Additionally, the higher-order scheme exhibits higher resolution and a greater number of troubled cells. Figure \ref{Fig:woodward1d2} illustrates the results when $M=100$, indicating that there is no significant change in the number of troubled cells or resolution compared to when $M=0.01$.

 \begin{figure}[htbp]
    \centering
    \includegraphics[width=0.3\textwidth]{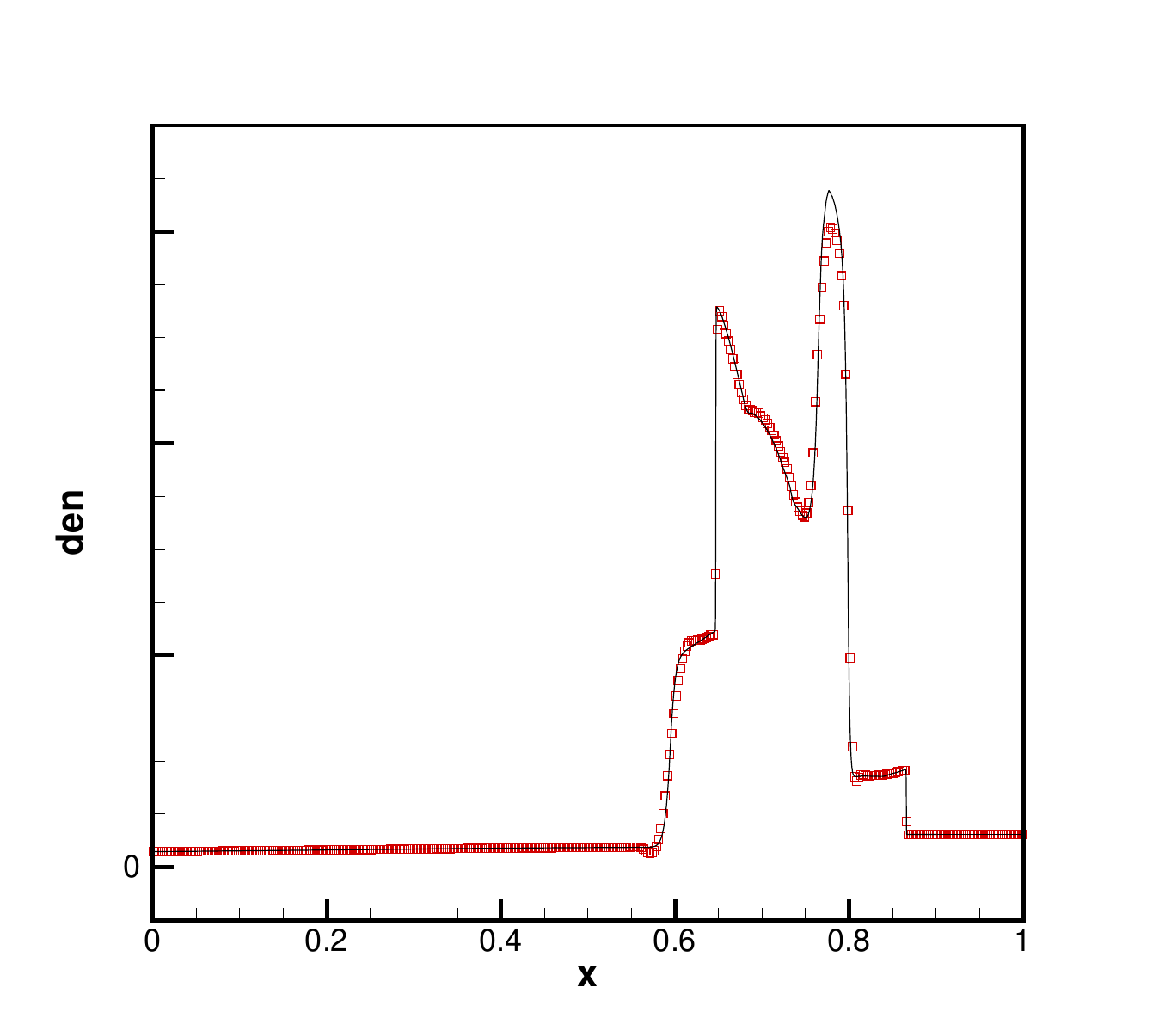}
    \includegraphics[width=0.3\textwidth]{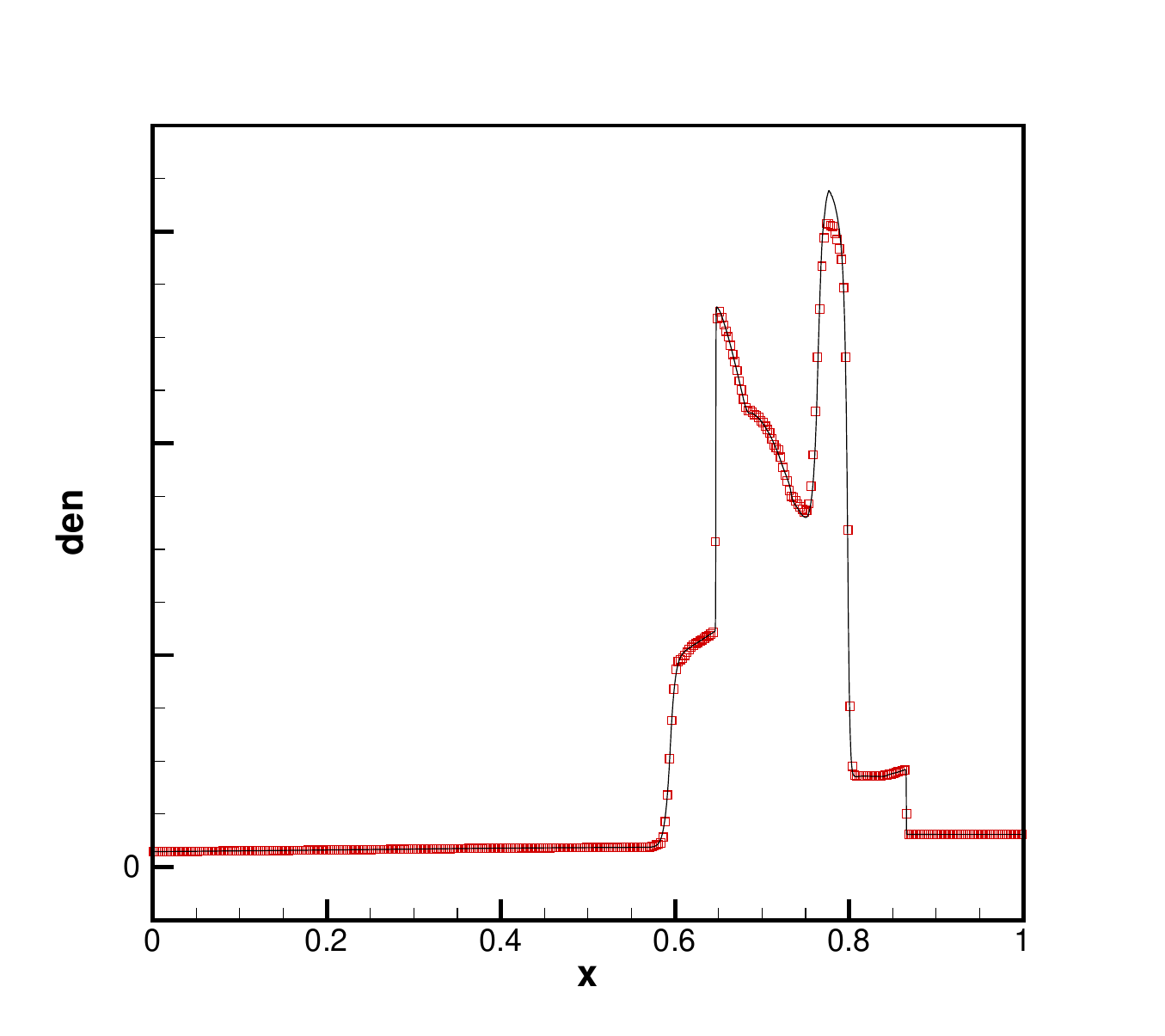}
    \includegraphics[width=0.3\textwidth]{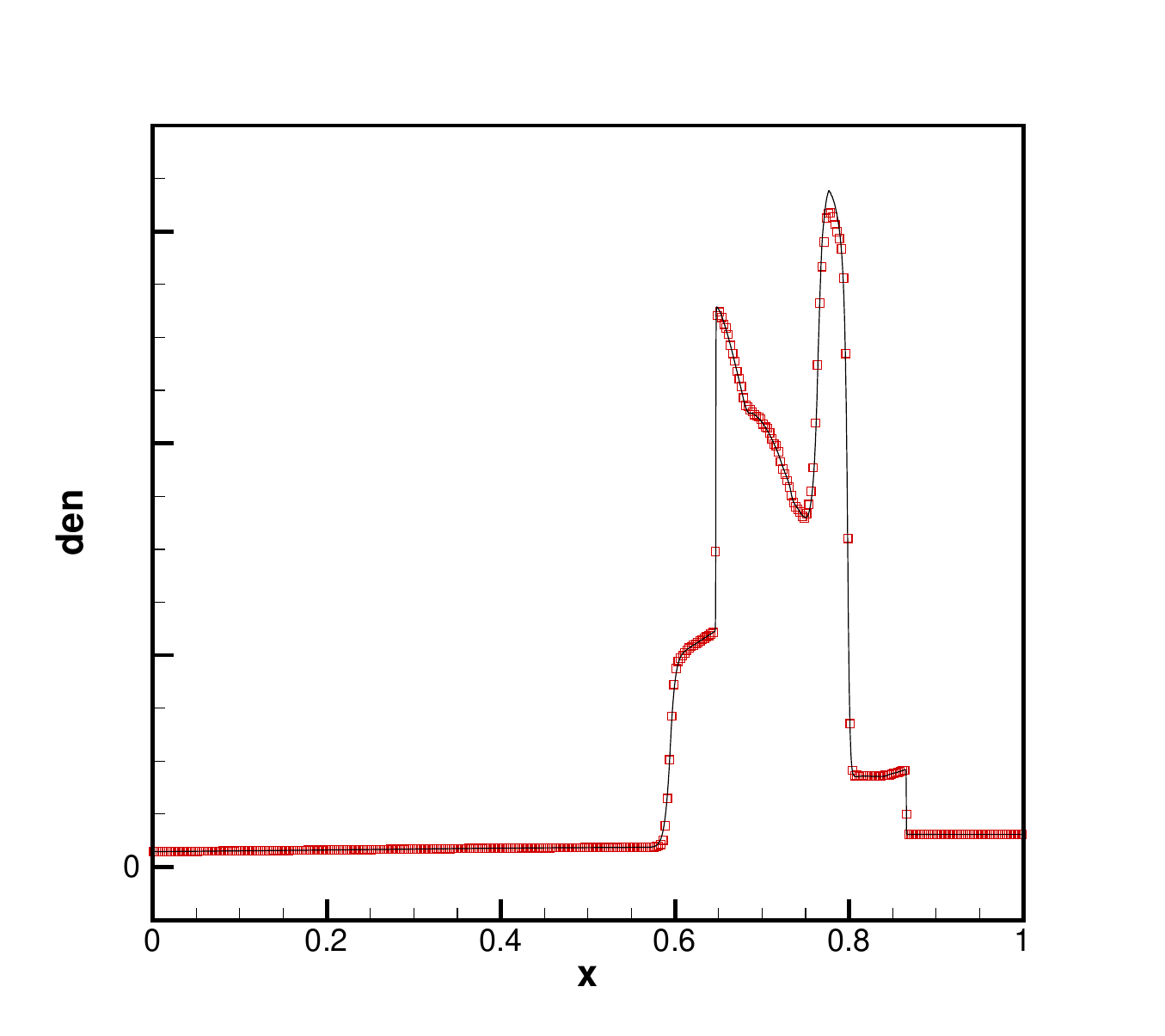}
    \includegraphics[width=0.3\textwidth]{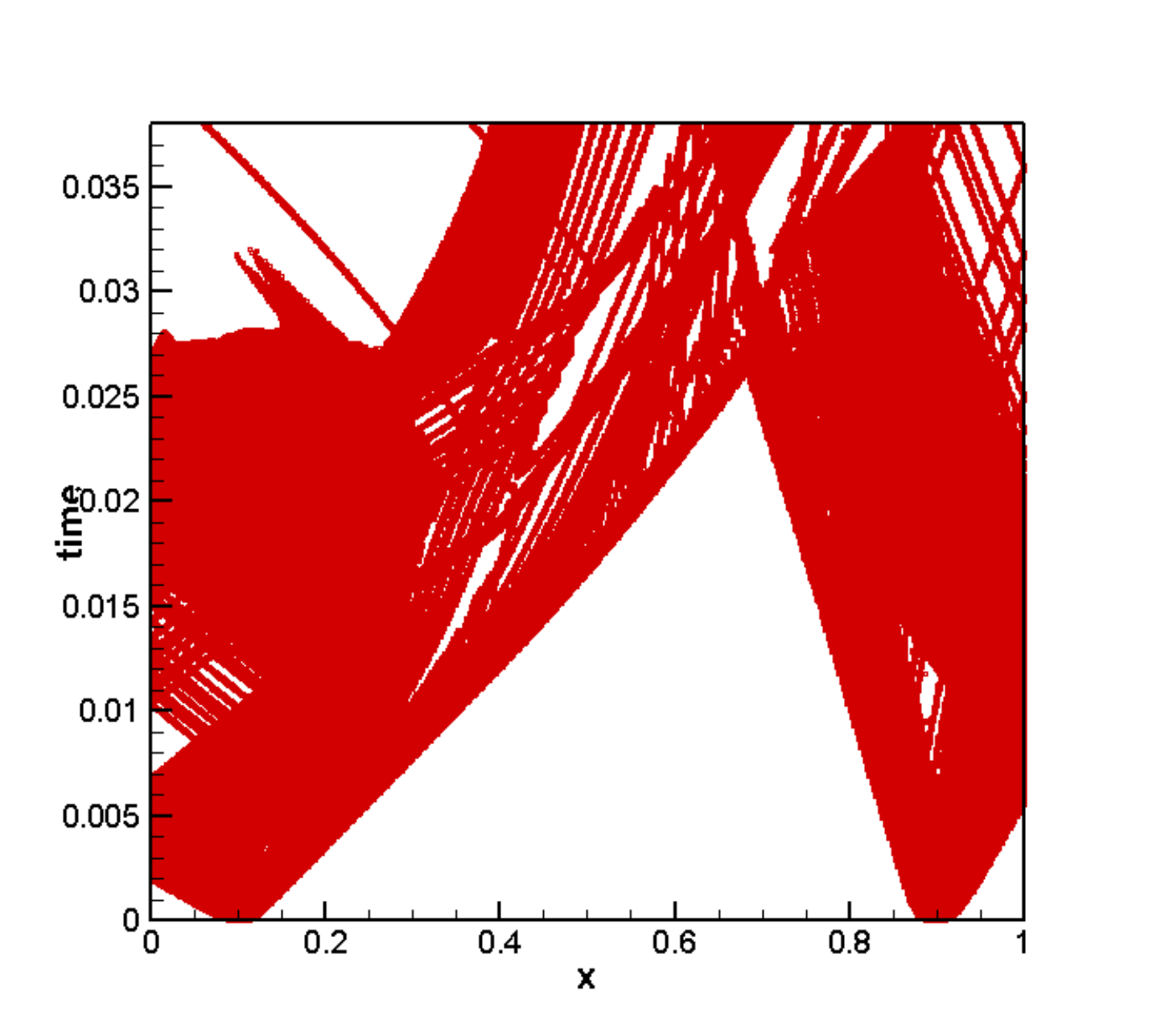}
    \includegraphics[width=0.3\textwidth]{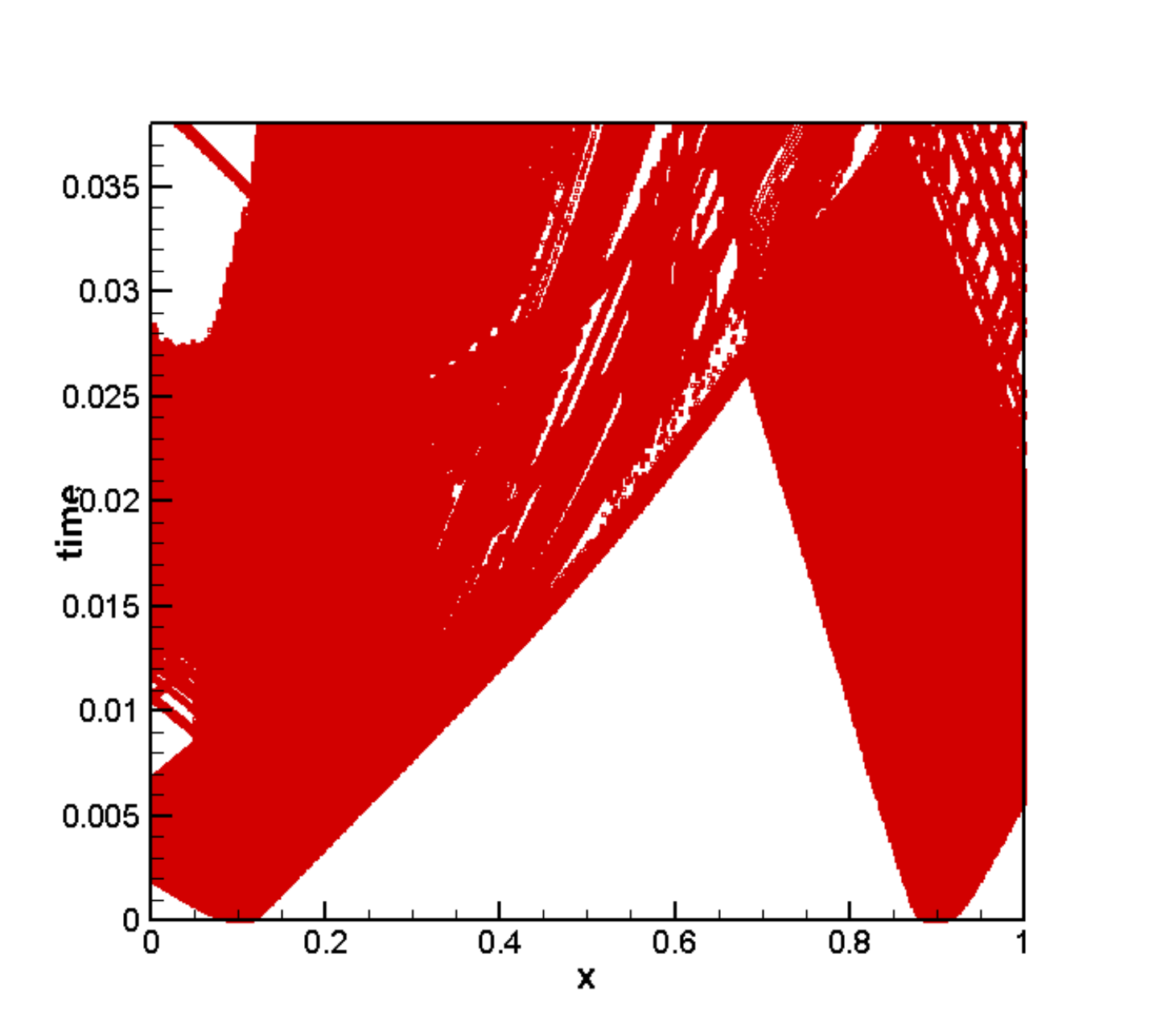}
    \includegraphics[width=0.3\textwidth]{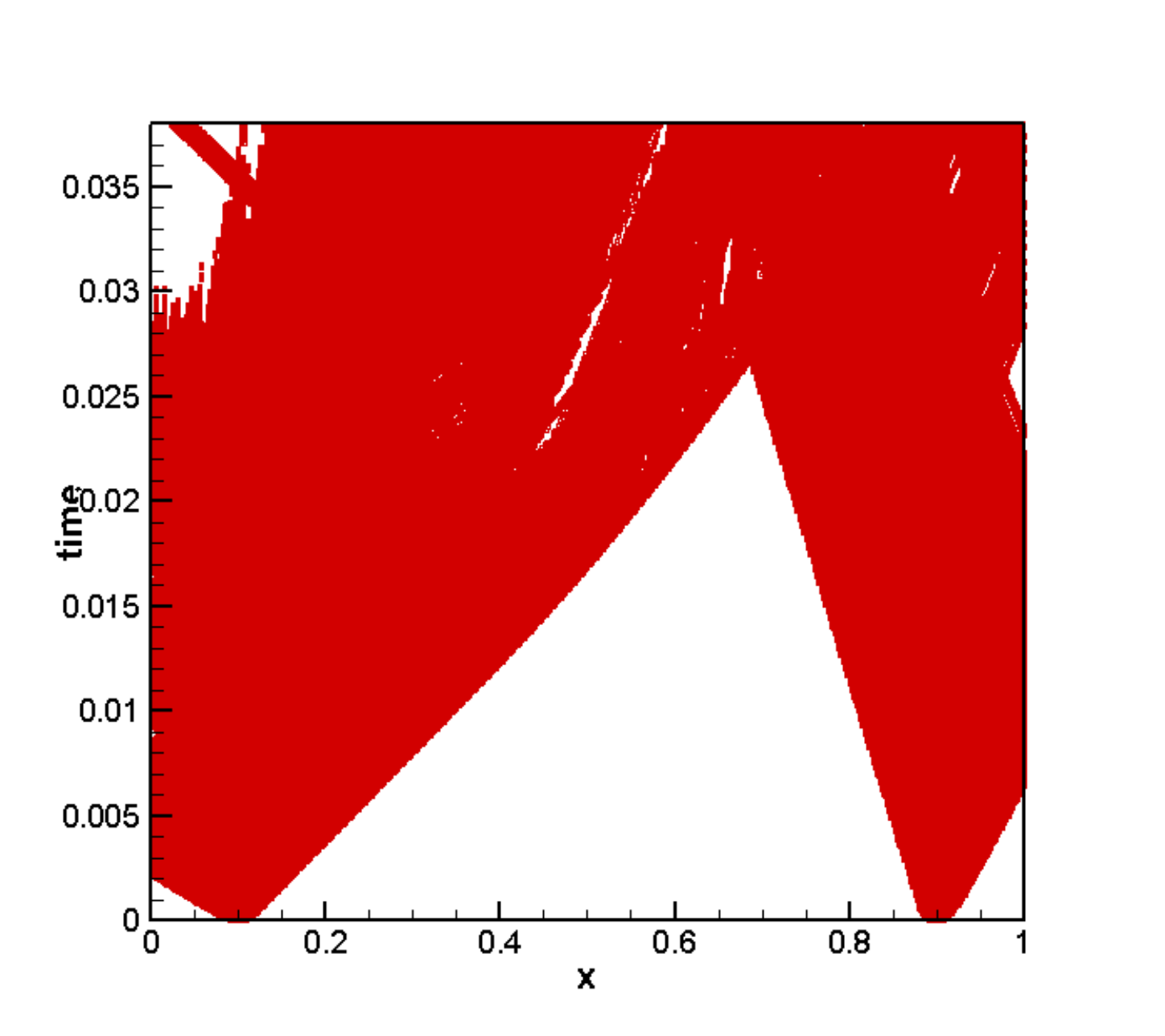}
     \caption{Example \ref{chap03:wood}: Explosion wave problem:
     The density numerical solution is obtained using 400 spectral volumes in {\tt SV-cvMSWNO3}-\texttt{SV-cvMSWENO5} scheme when $M=0.01$ at $t=0.038$.  The solid line and symbol ``$\square$'' represent the exact solution and {\tt SV-cvMSWENO} results, respectively.}
    \label{Fig:woodward1d}
 \end{figure}

 \begin{figure}[htbp]
    \centering
    \includegraphics[width=0.3\textwidth]{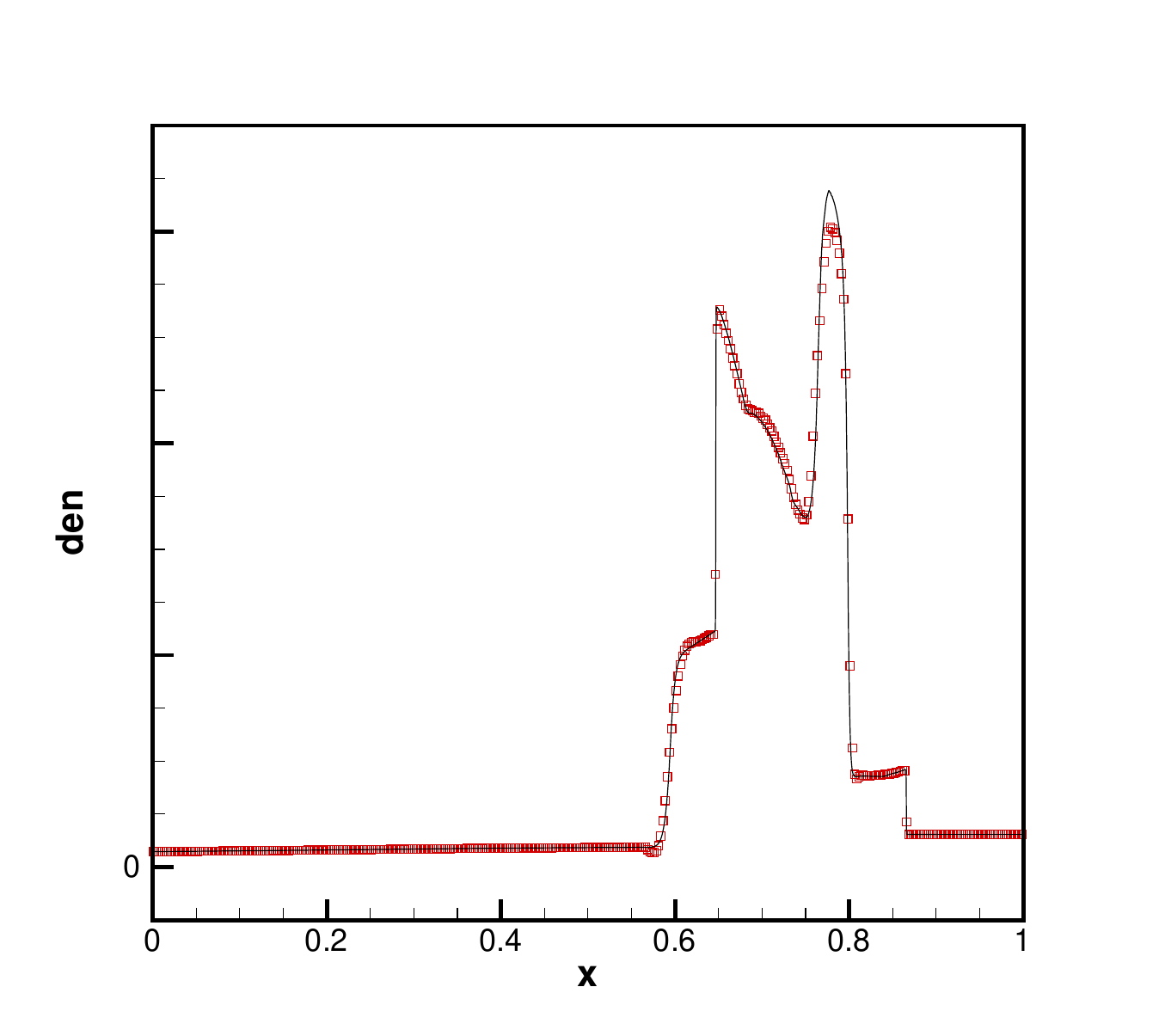}
    \includegraphics[width=0.3\textwidth]{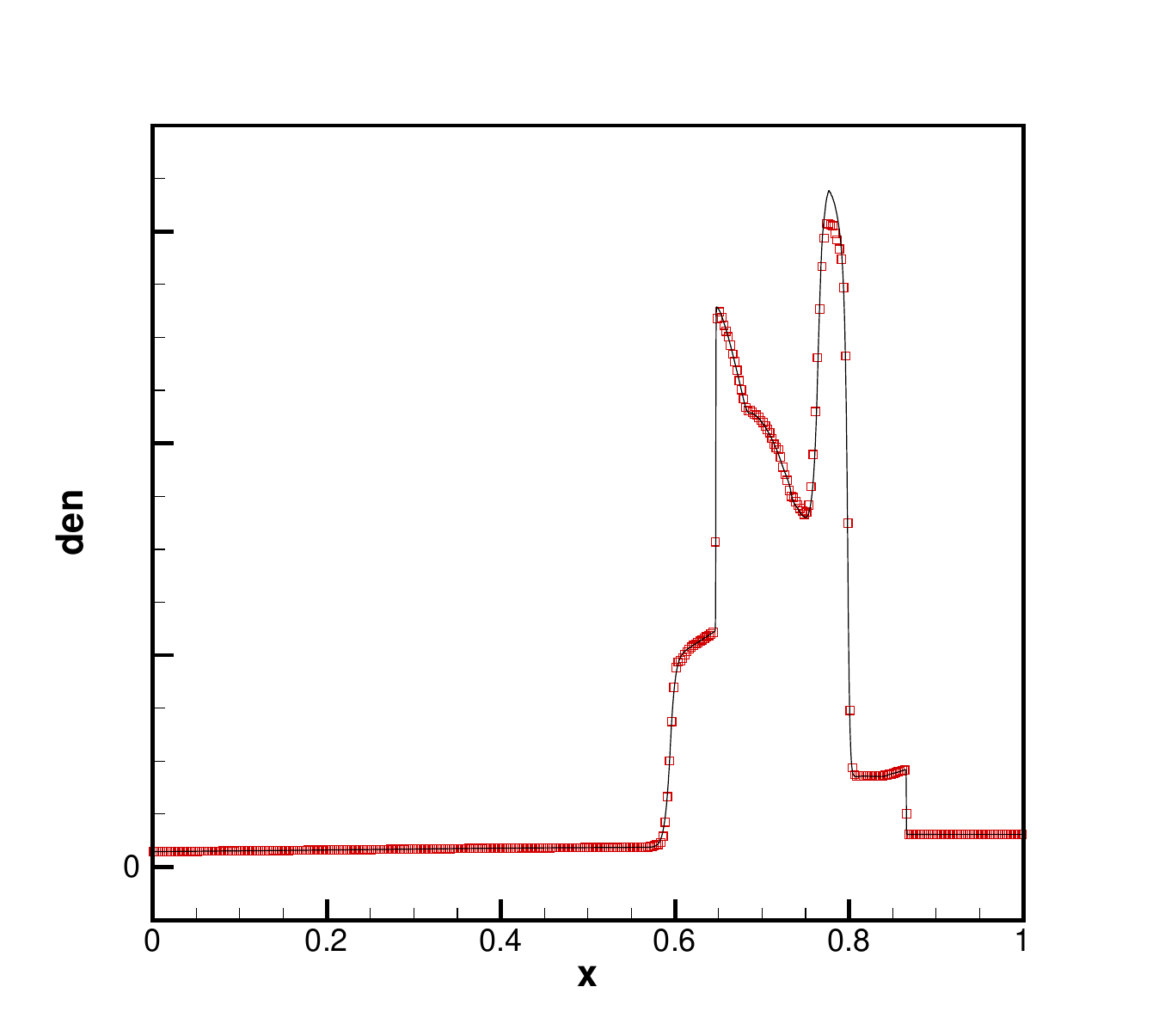}
    \includegraphics[width=0.3\textwidth]{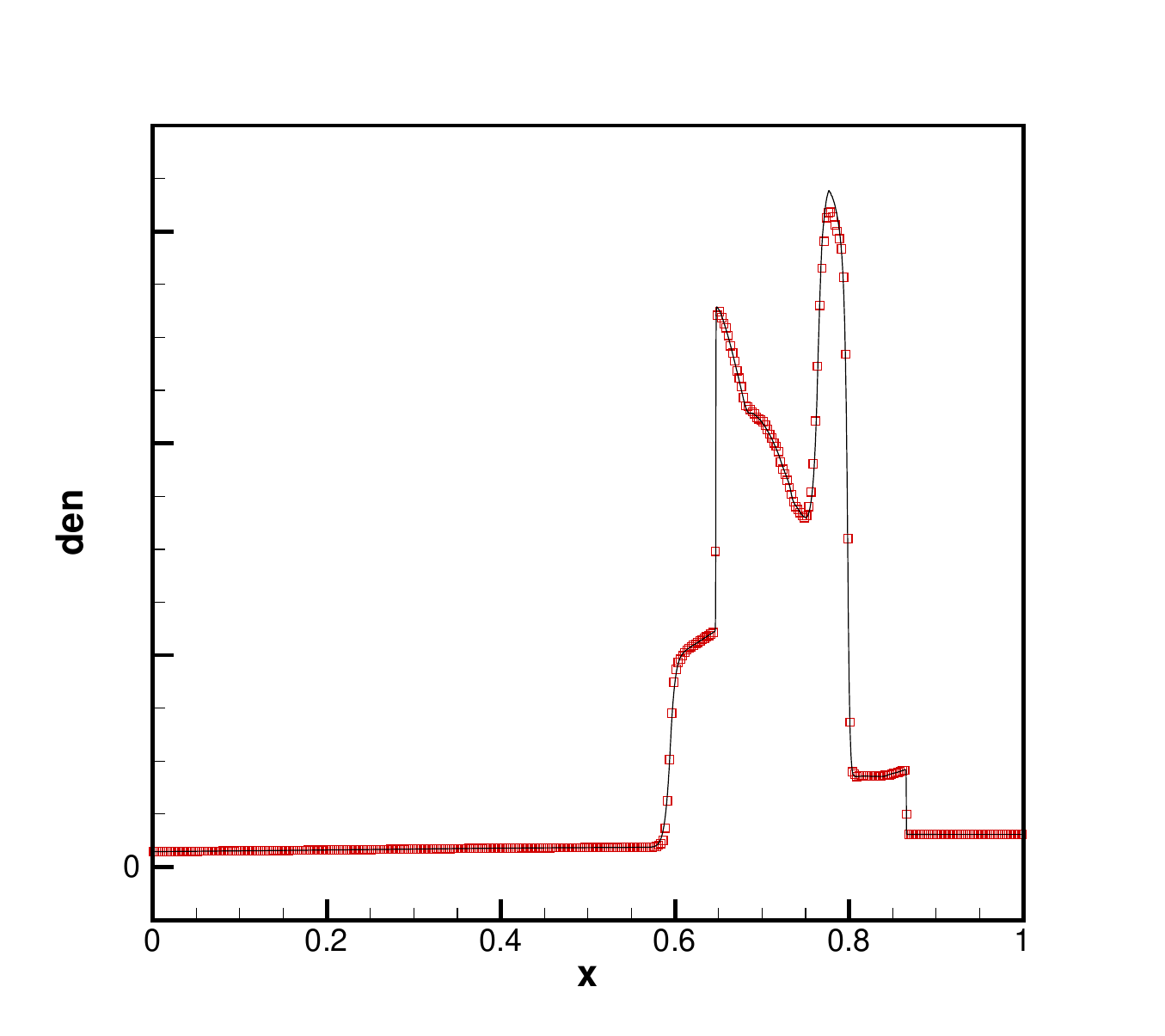}
    \includegraphics[width=0.3\textwidth]{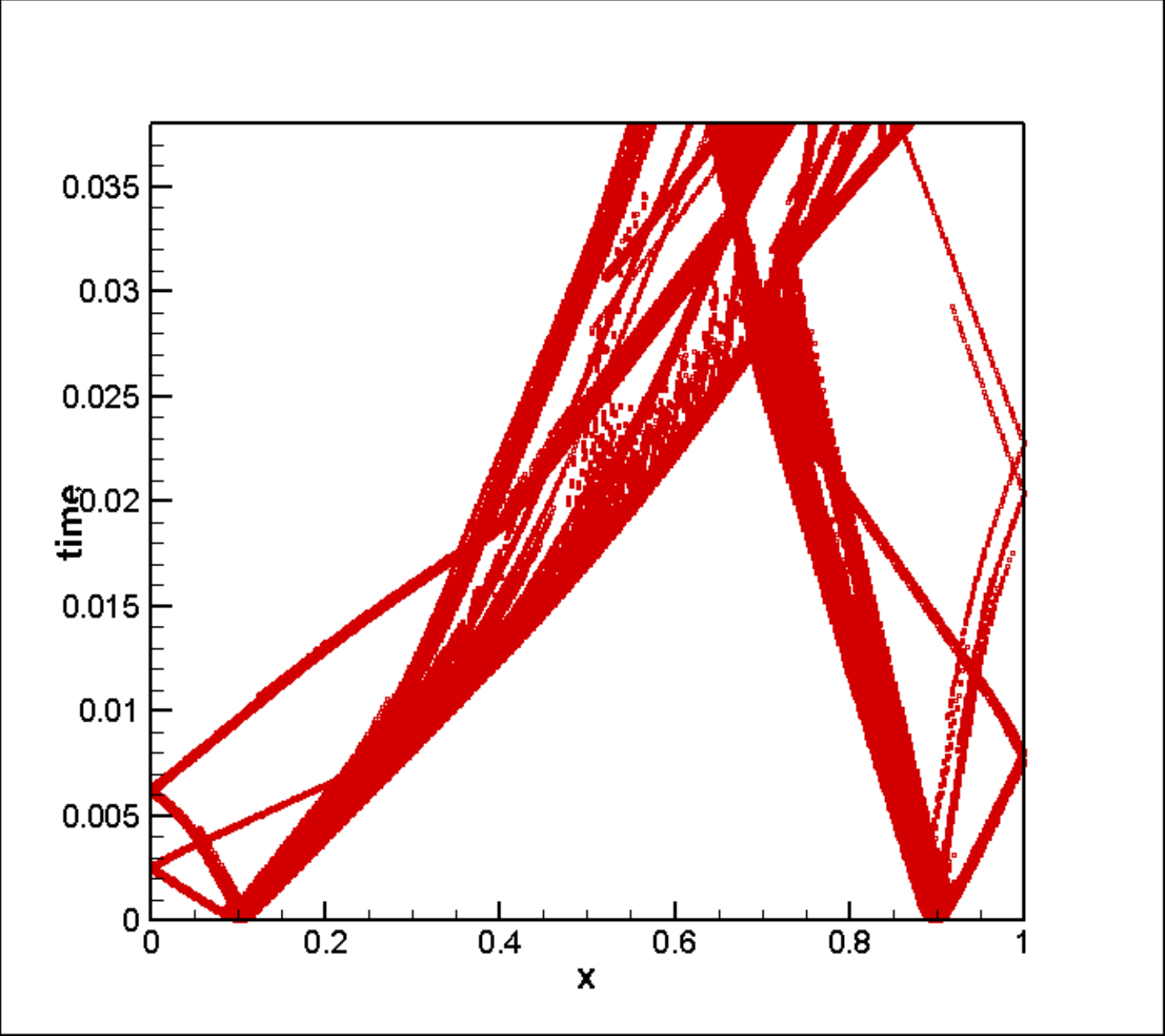}
    \includegraphics[width=0.3\textwidth]{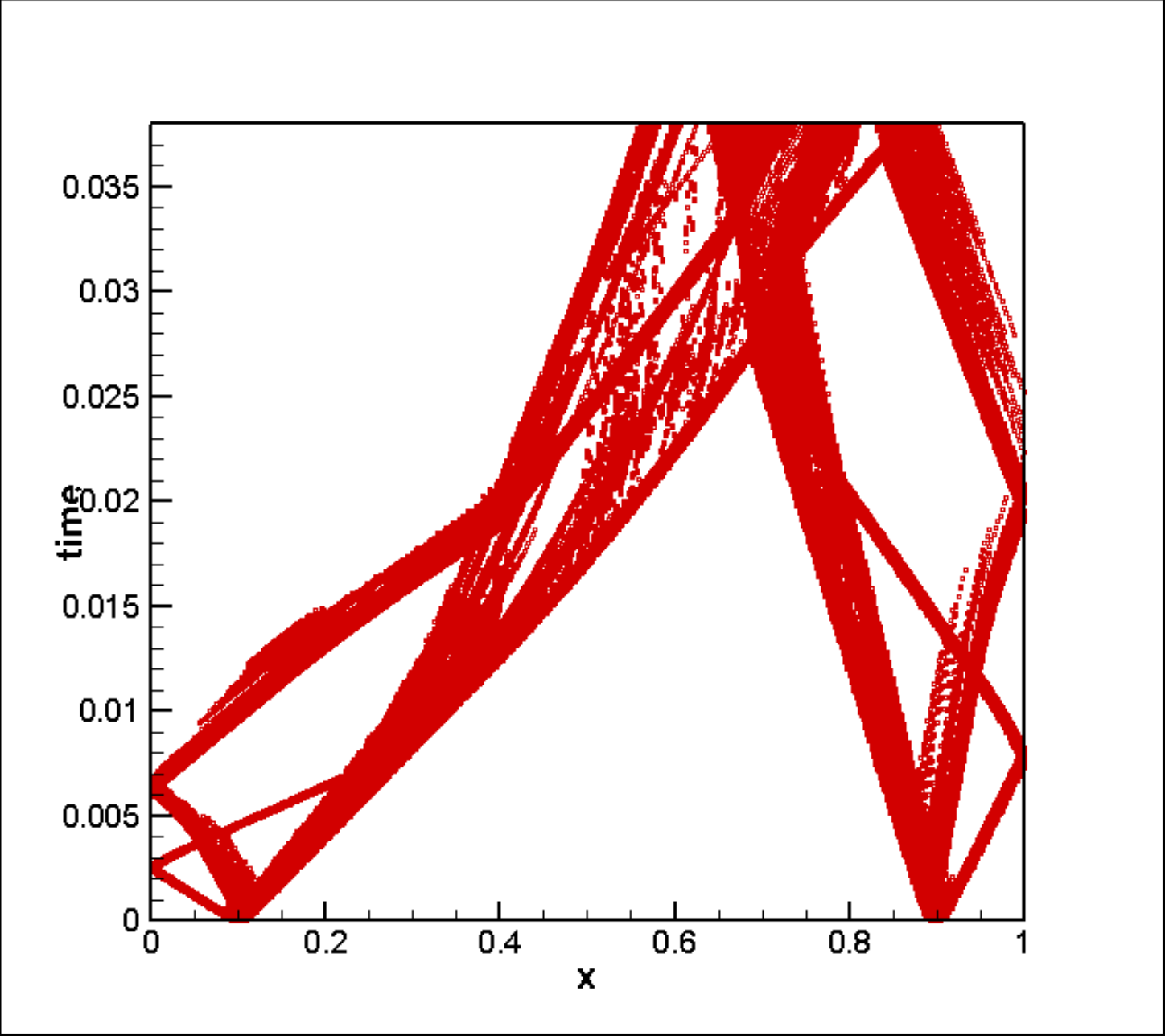}
    \includegraphics[width=0.3\textwidth]{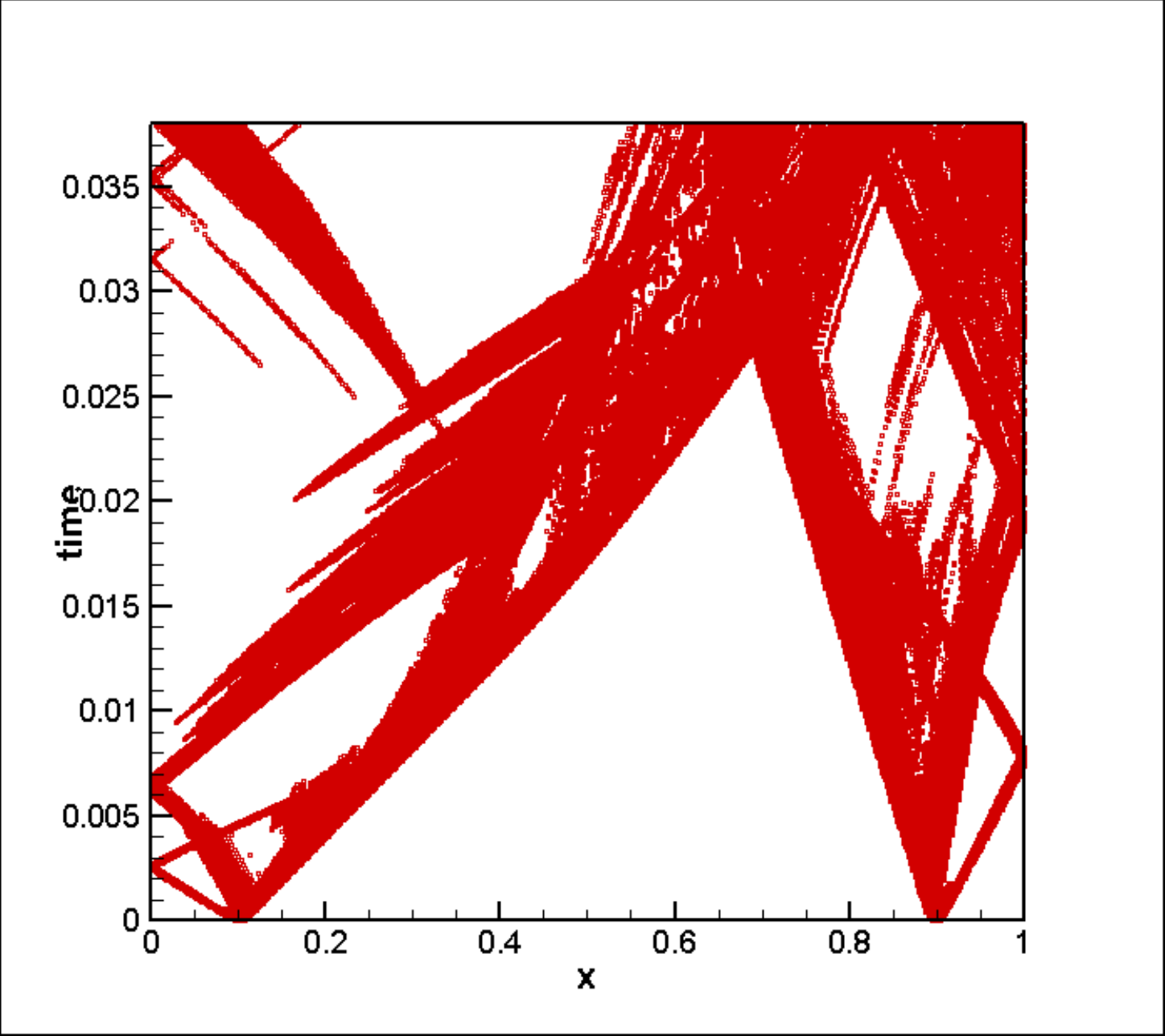}
     \caption{Example \ref{chap03:wood}: Explosion wave problem:
     The density numerical solution is obtained using 400 spectral volumes in {\tt SV-cvMSWNO3}-\texttt{SV-cvMSWENO5} scheme when $M=100$ at $t=0.038$.  The solid line and symbol ``$\square$'' represent the exact solution and {\tt SV-cvMSWENO} results, respectively.}
    \label{Fig:woodward1d2}
 \end{figure}

\end{example}

\begin{example}[2D Riemannproblem  \uppercase\expandafter{\romannumeral1}]\label{chap03:Riemann2d1}\rm
Next, we consider some two dimensional Riemann problems\cite{lax1998}, where the initial values were in four different constant states in the
four quadrants. The computational domain $\Omega$ is defined as $[0, 1] \times [0, 1]$ with outflow boundaries at both endpoints of the region.
Due to the space limitations, we only present the results of two of the Riemann problems here.

The initial condition of the first Riemann problem is
 \begin{align*}
(\rho,U,V,P)(x,y,0)=\begin{cases}
 (0.5313,         0,     0,  0.4),  &   x>0.5,  y>0.5,\\
 (1.0, 0.7276,     0,  1), &   x<0.5,  y>0.5,\\
 (0.8,    0, 0, 1),  &   x<0.5,  y<0.5,\\
 (1.0,      0, 0.7276, 1), &   x>0.5,  y<0.5.
\end{cases}\end{align*}
It contains a rightward fast shockwave and an upward fast shockwave, which intersect in the first quadrant to form a Mach stem and generate a jet along the
diagonal direction. The fluid in the lower left third quadrant is a static high pressure zone, forming a horizontal slip line with the fluid in the second
quadrant and a vertical contact discontinuity with the fluid in the fourth quadrant.

Figure \ref{Fig:riemann2d-1-k3} presents the contour map of the density numerical solution and the troubled cells calculated at $t=0.25$ by the \texttt{SV-cvMSWENO3} scheme using $100\times100$ uniform spectral volume grids, from top to bottom, followed by $M=0.01$(top),$M=100$(middle),$M=200$(bottom).
Figure \ref{Fig:riemann2d-1-k4} gives the results by \texttt{SV-cvMSWENO4} scheme. Here, 30 contour lines are taken, with the range uniformly set from 0.6 to
1.6. The results show that correct results can be obtained under all three values of M. And as M increases, the number of detected troubled cells decreases
gradually, while the resolution improves. The resolution of the \texttt{SV-cvMSWENO4} scheme is better than that of the \texttt{SV-cvMSWENO3} scheme with the same M.

\begin{figure}[htbp]
    \centering
    \includegraphics[width=0.48\textwidth]{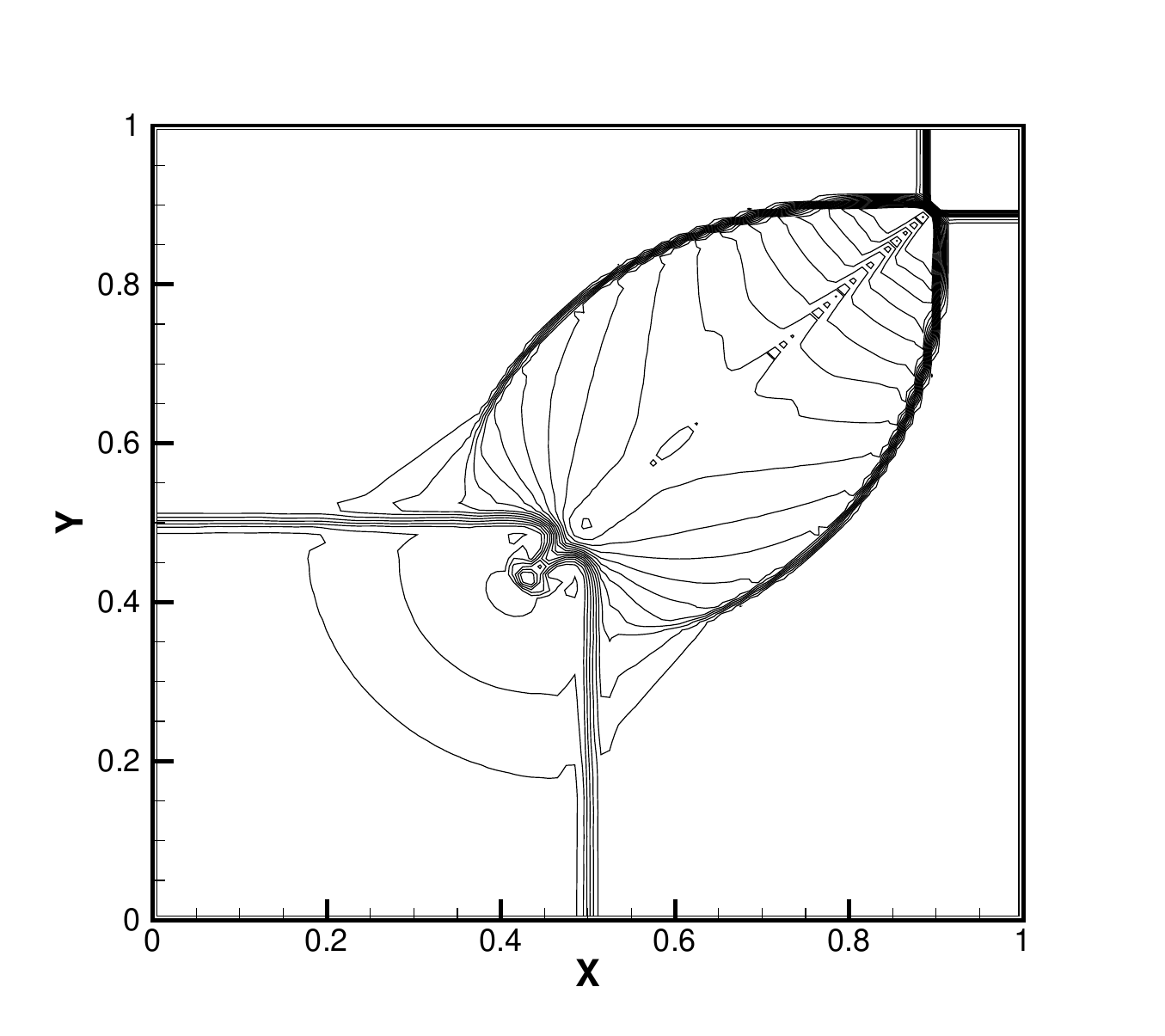}
    \includegraphics[width=0.48\textwidth]{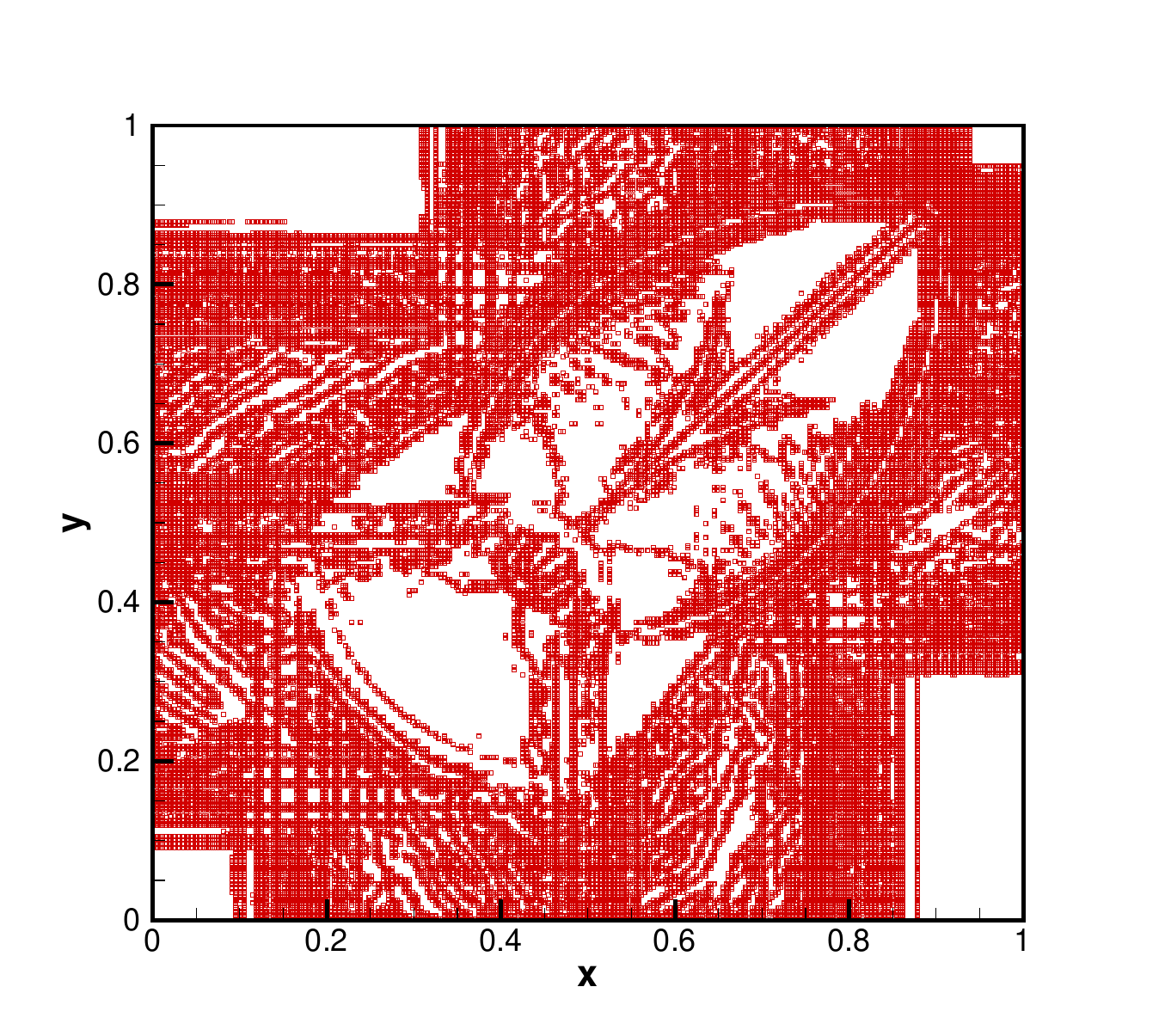}
    \includegraphics[width=0.48\textwidth]{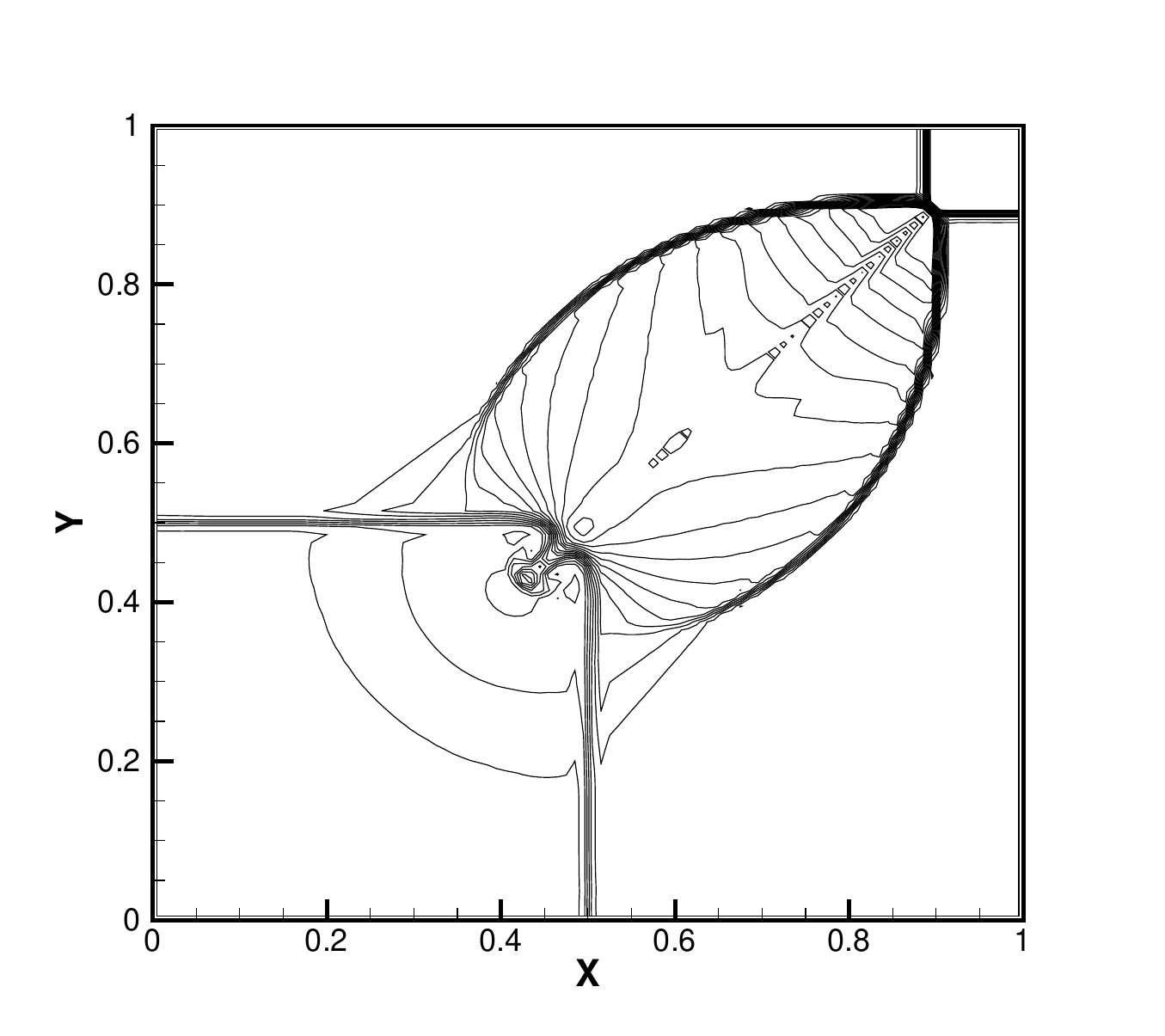}
    \includegraphics[width=0.48\textwidth]{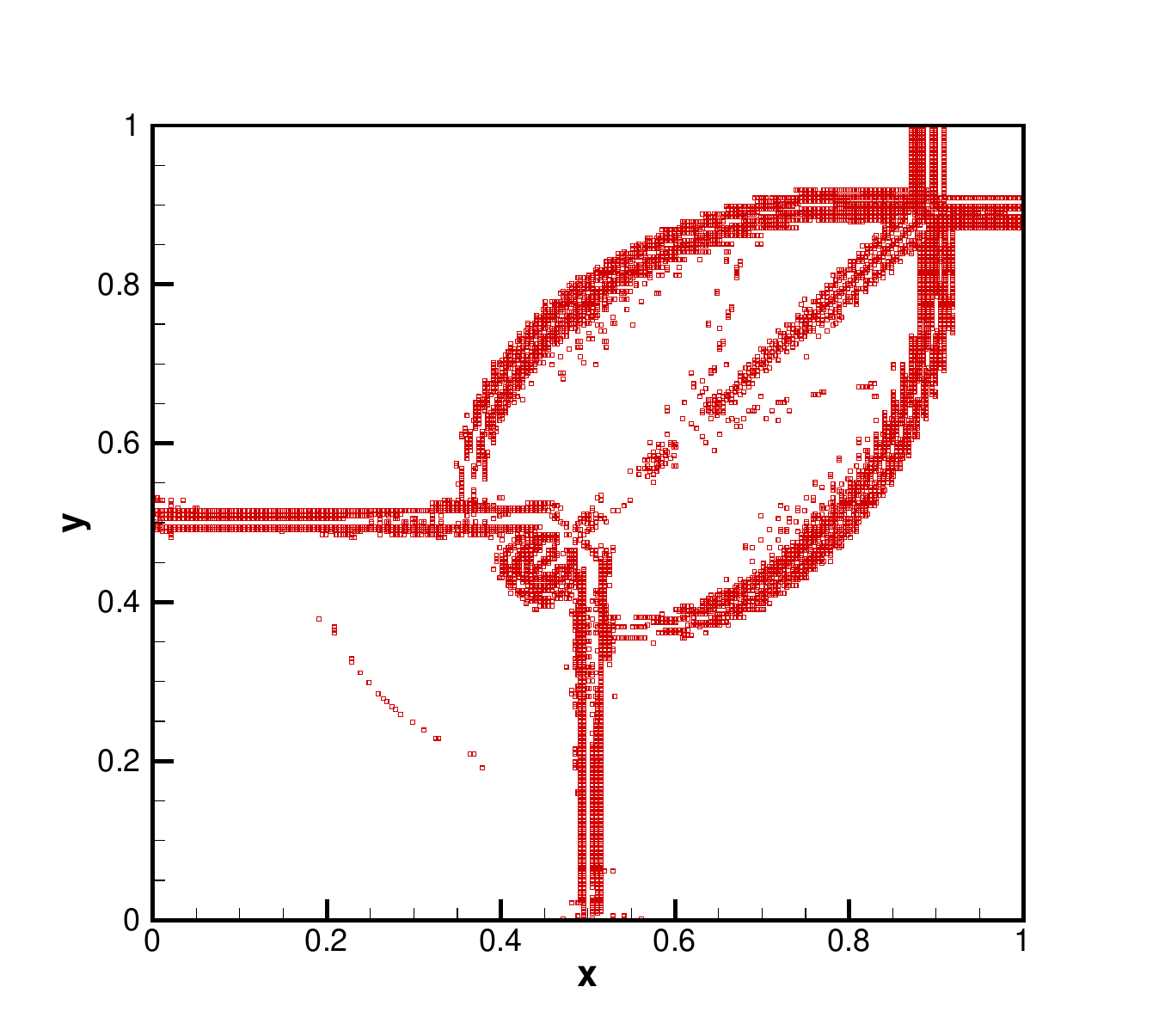}
    \includegraphics[width=0.48\textwidth]{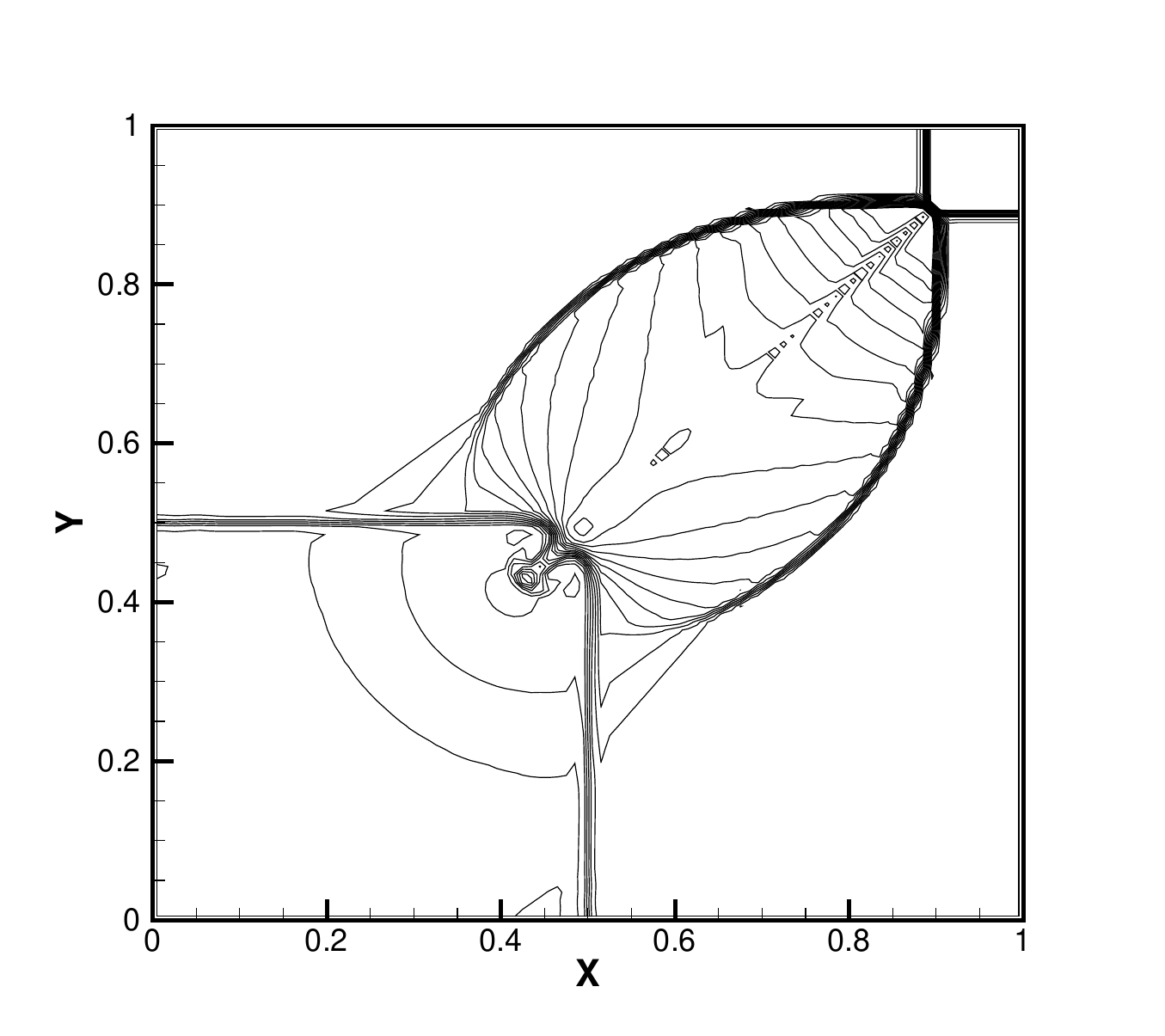}
    \includegraphics[width=0.48\textwidth]{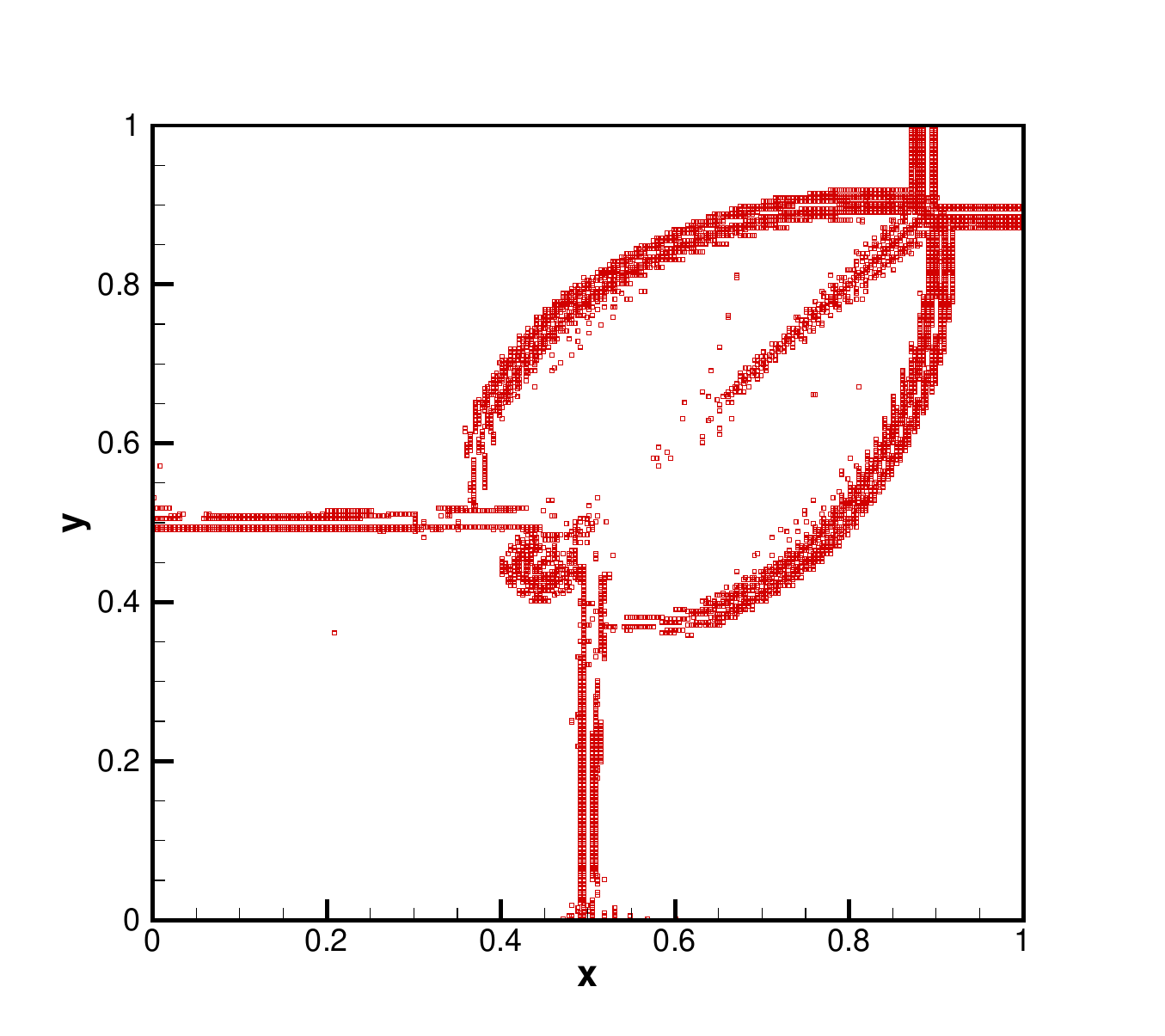}
     \caption{Example \ref{chap03:Riemann2d1}: 2D Riemannproblem \uppercase\expandafter{\romannumeral1}:
     30 equally density contour lines(left) from 0.6 to 1.6 and troubled cells(right) of \texttt{SV-cvMSWENO3} scheme when $M=0.01$(top),$M=100$(middle),$M=200$(bottom).}
    \label{Fig:riemann2d-1-k3}
 \end{figure}

 \begin{figure}[htbp]
    \centering
    \includegraphics[width=0.48\textwidth]{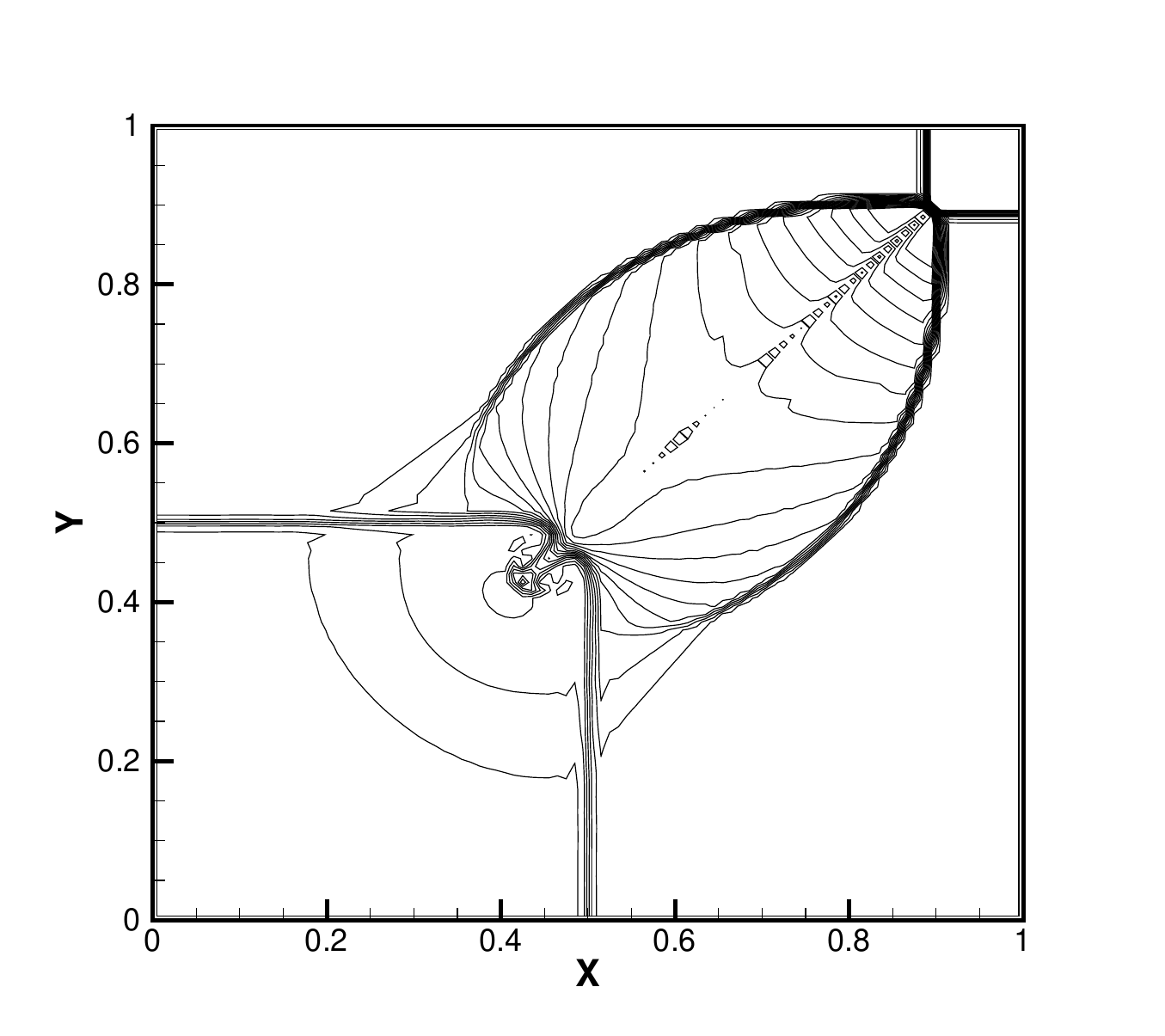}
    \includegraphics[width=0.48\textwidth]{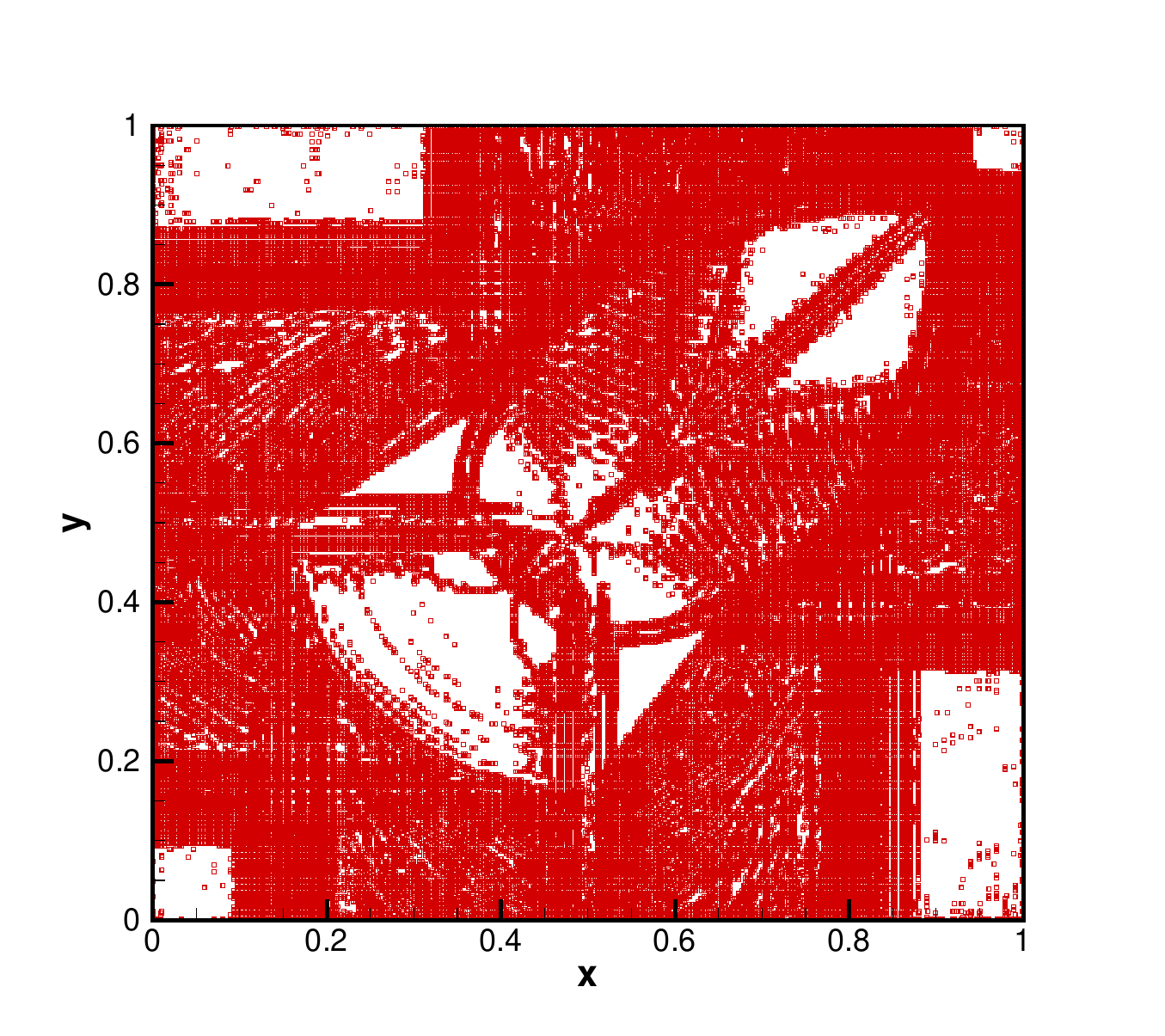}
    \includegraphics[width=0.48\textwidth]{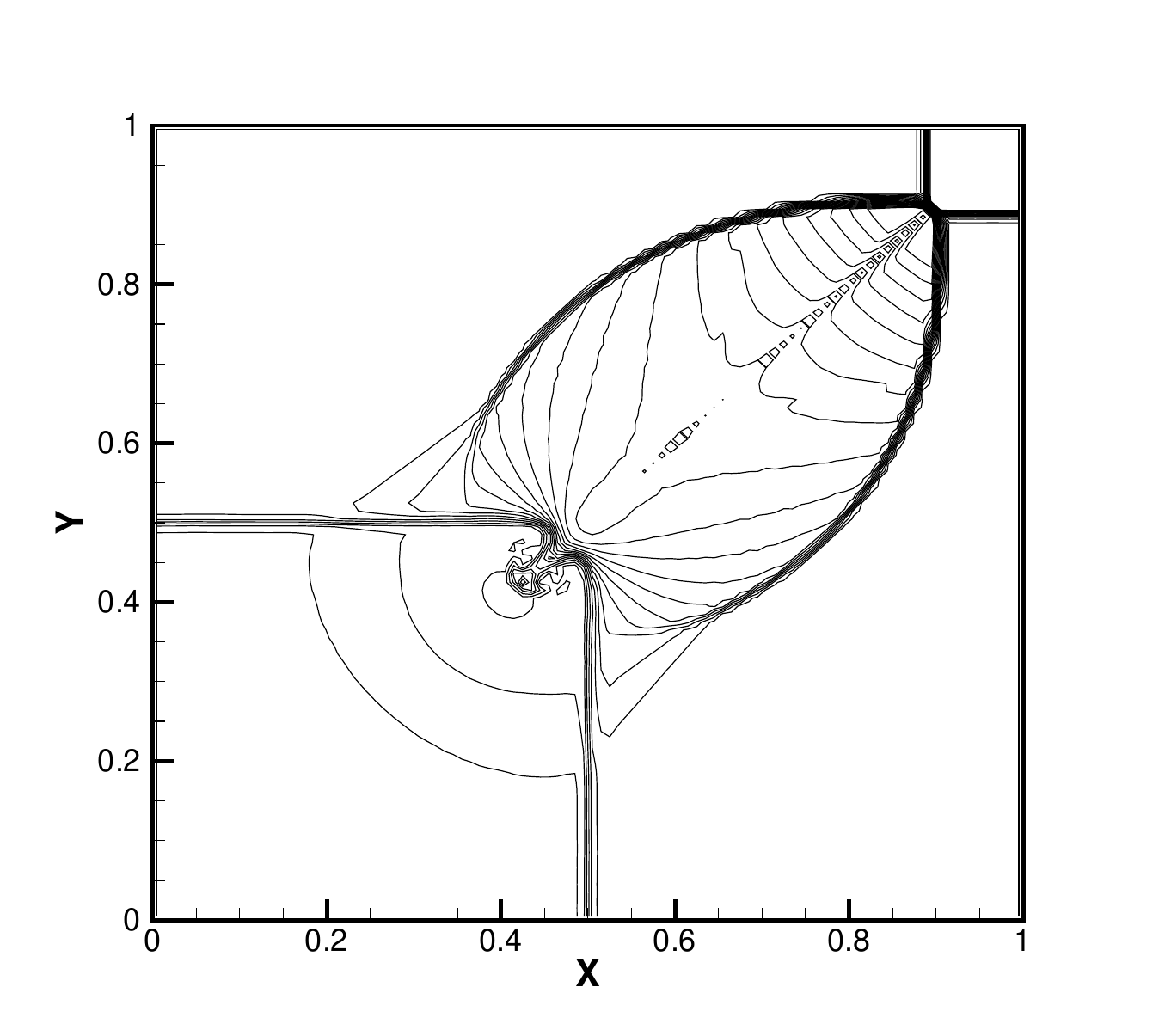}
    \includegraphics[width=0.48\textwidth]{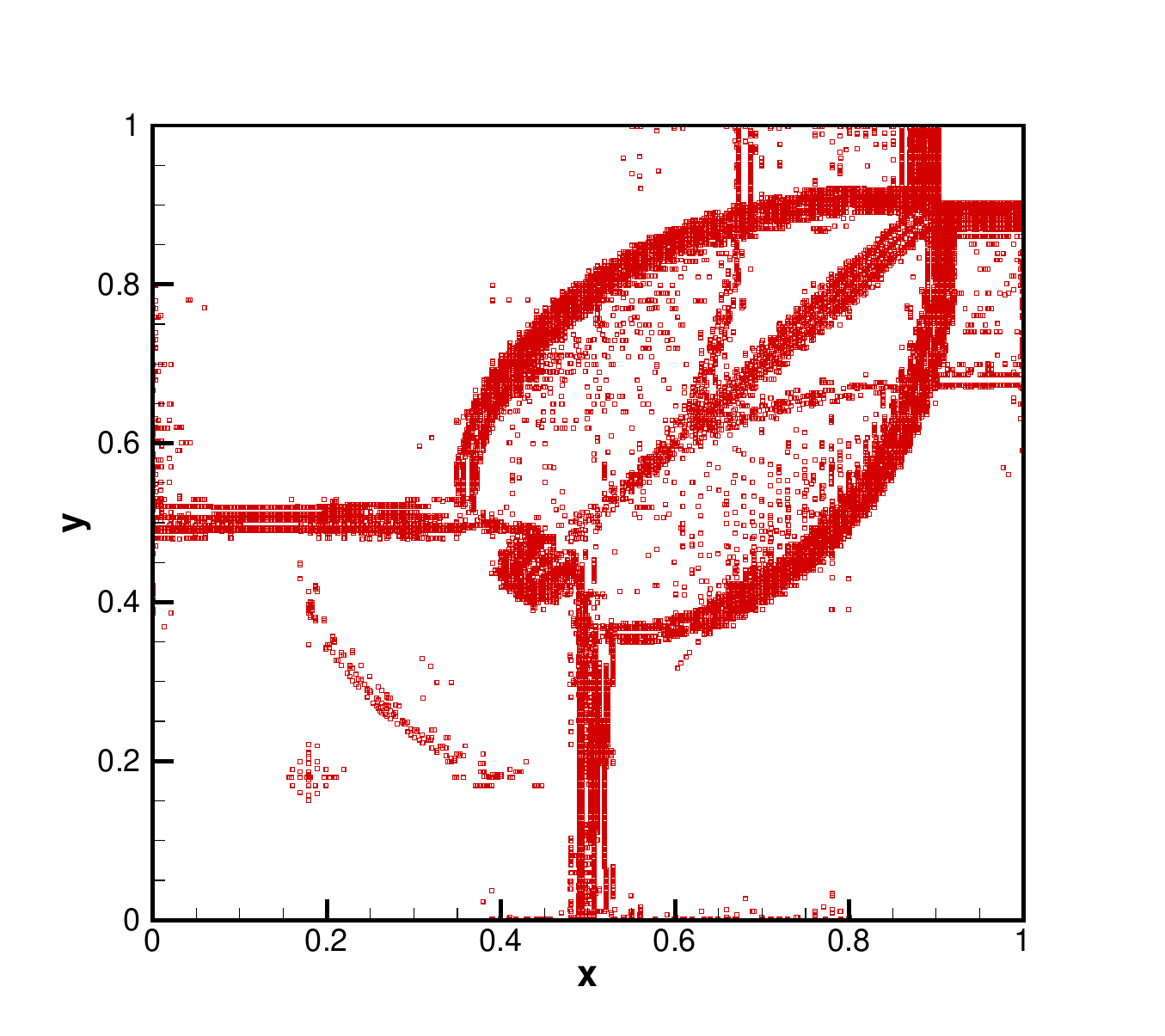}
    \includegraphics[width=0.48\textwidth]{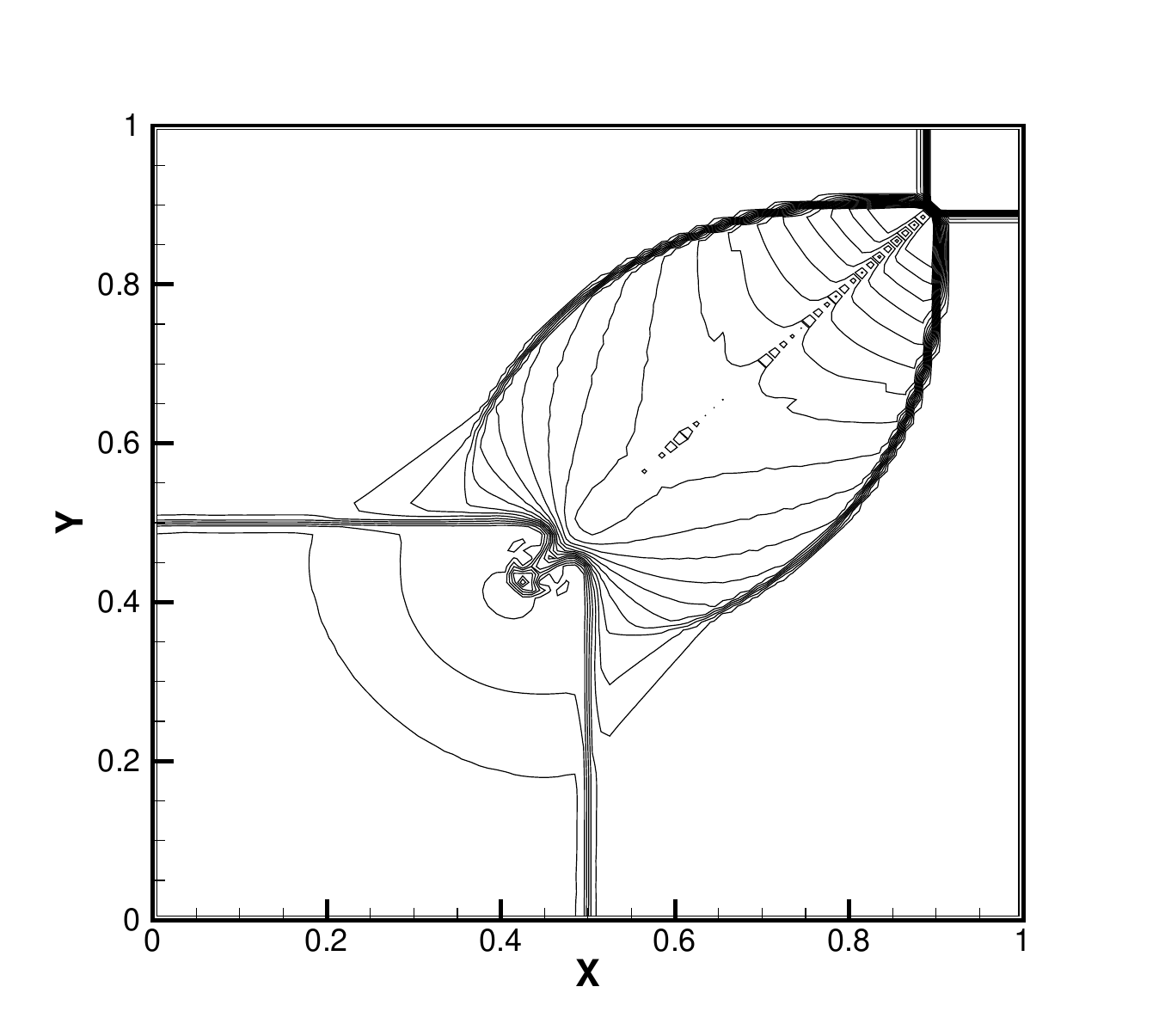}
    \includegraphics[width=0.48\textwidth]{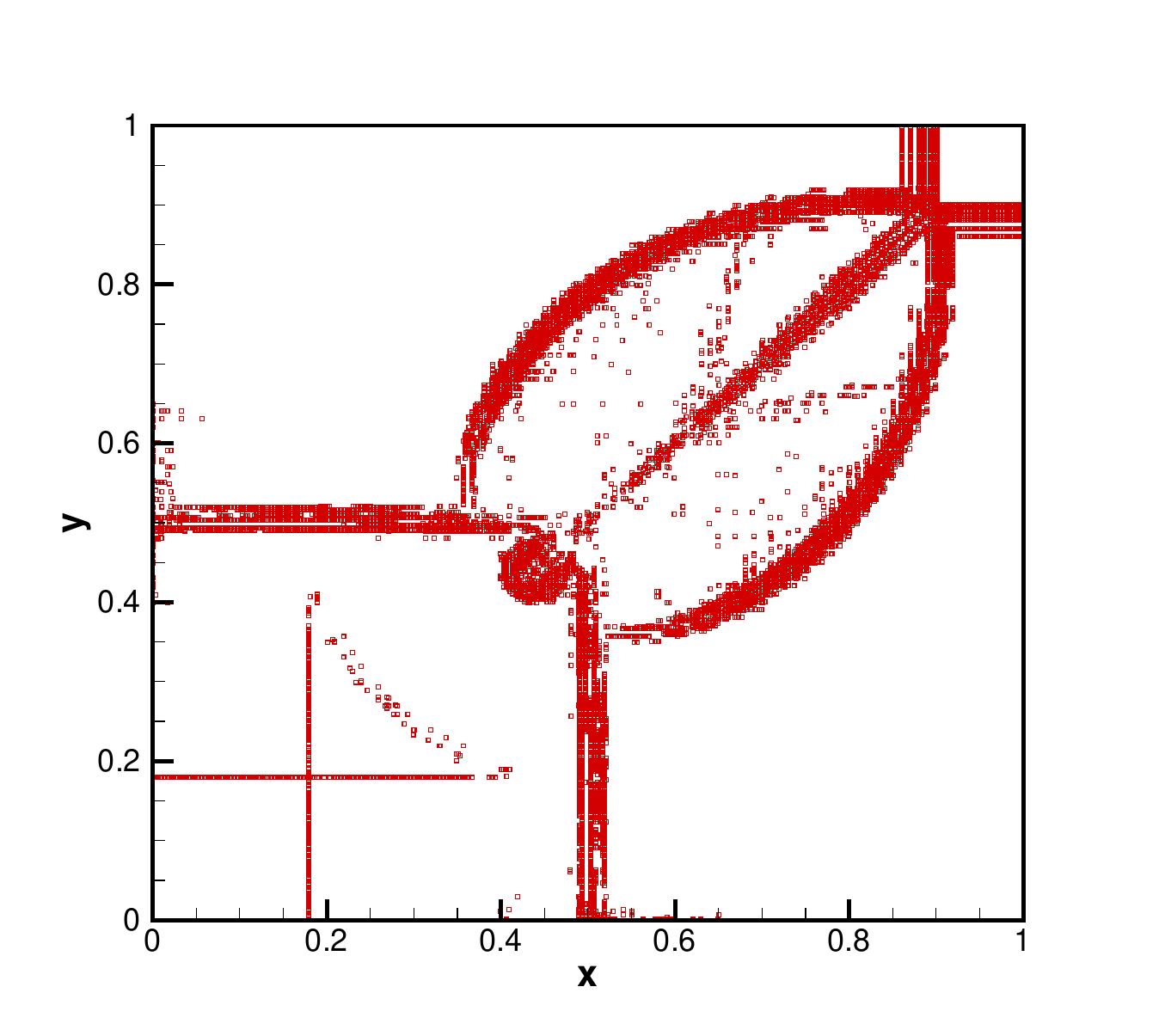}
     \caption{Example \ref{chap03:Riemann2d1}: 2D Riemannproblem \uppercase\expandafter{\romannumeral1}:
     30 equally density contour lines(left) from 0.6 to 1.6 and troubled cells(right) of \texttt{SV-cvMSWENO4} scheme when $M=0.01$(top),$M=100$(middle),$M=200$(bottom).}
    \label{Fig:riemann2d-1-k4}
 \end{figure}

\end{example}

\begin{example}[2D Riemannproblem \uppercase\expandafter{\romannumeral2}]\label{chap03:Riemann2d2}\rm
This Riemann problem\cite{lax1998} tests the resolution ability of numerical schemes for complex multi-
dimensional wave structures such as oblique shock waves, rarefaction wave, and contact
discontinuous coupling. The initial conditions are
 \begin{align*}
(\rho,U,V,P)(x,y,0)=\begin{cases}
 (1.0,         0.1,     -0.3,  1.0),  &   x>0.5,  y>0.5,\\
 (0.5197, -0.6259, -0.3,  0.4), &   x<0.5,  y>0.5,\\
 (0.8,    0.1, -0.3, 0.4),  &   x<0.5,  y<0.5,\\
 (0.5313,      0.1, 0.4276, 0.4), &   x>0.5,  y<0.5,  
\end{cases}\end{align*}
It contains a rightward fast rarefaction wave, a downward fast shock, and two strong shear slip
lines that form a vortex structure at the center.

Figure \ref{Fig:riemann2d-2-k3} and figure \ref{Fig:riemann2d-2-k4} respectively present the
density contour and troubled cells calculated using $100\times100$ uniform spectral volume grids at
$t=0.2$ by the \texttt{SV-cvMSWENO3} scheme and \texttt{SV-cvMSWENO4} scheme.
The contour figure takes $30$ uniform contour lines from 0.55 to 1.06. The figures show the
results of $M=0.01$(top),$M=100$(middle), and $M=200$(bottom) respectively.
The results show that when M=0.01, there are a lot of troubled cells. As M increases, the number
of detected troubled cells decreases gradually, and the resolution increases gradually in both {\tt SV-
cvMSWENO3} scheme and \texttt{SV-cvMSWENO4} scheme.
With the same M value, the \texttt{SV-cvMSWENO4} scheme has better resolution than the {\tt SV-
cvMSWENO3} scheme.

\begin{figure}[htbp]
    \centering
    \includegraphics[width=0.48\textwidth]{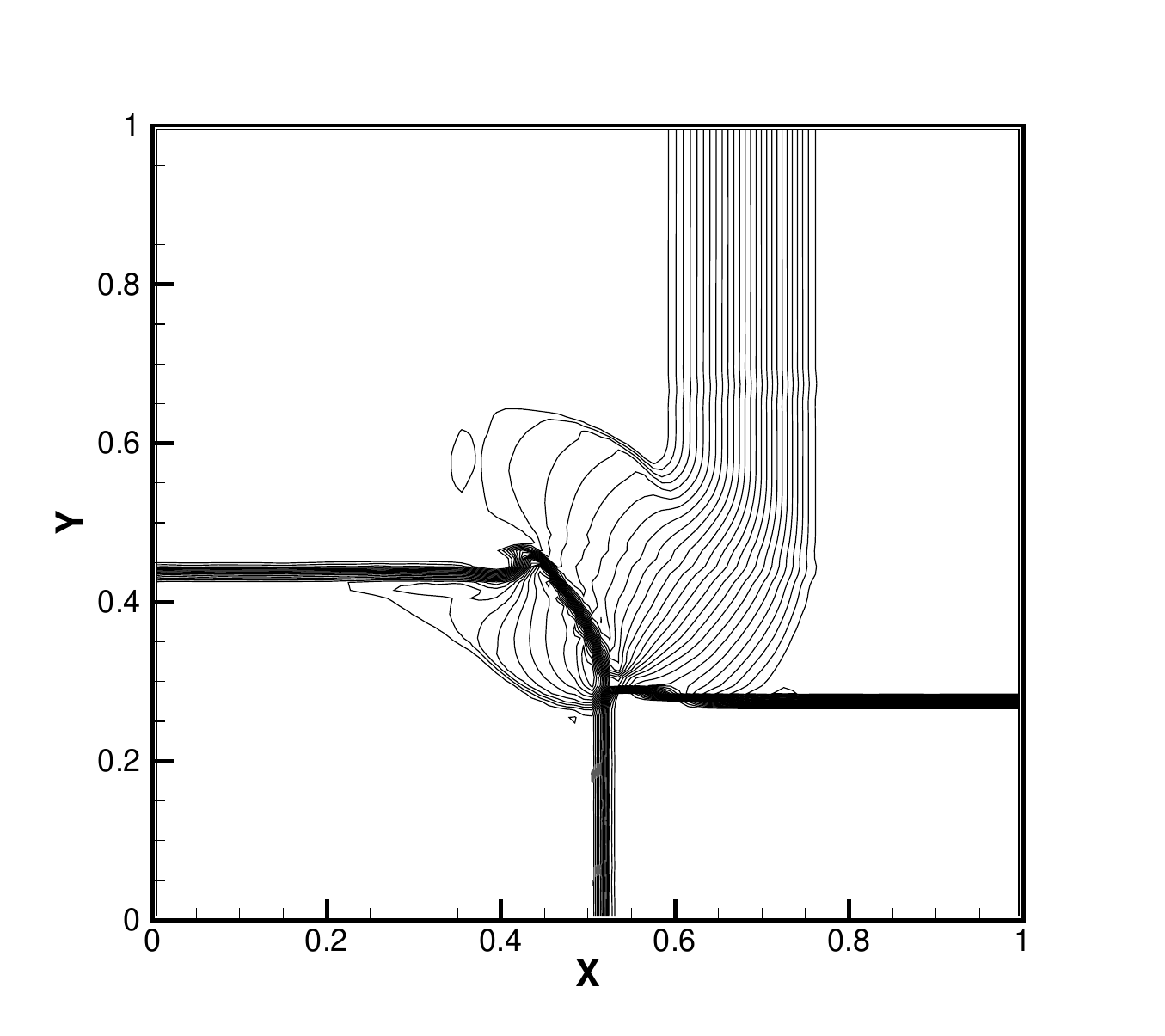}
    \includegraphics[width=0.48\textwidth]{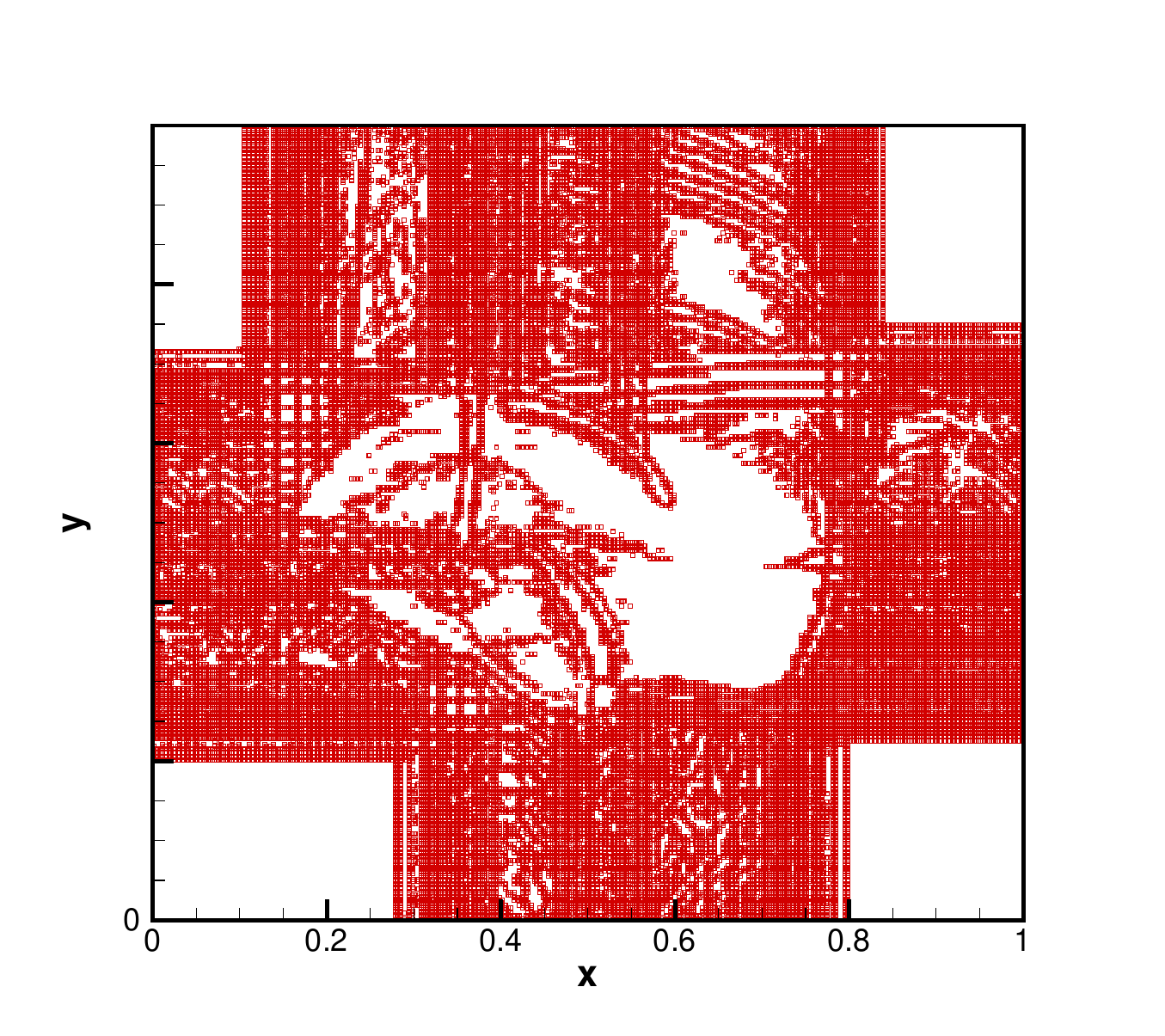}
    \includegraphics[width=0.48\textwidth]{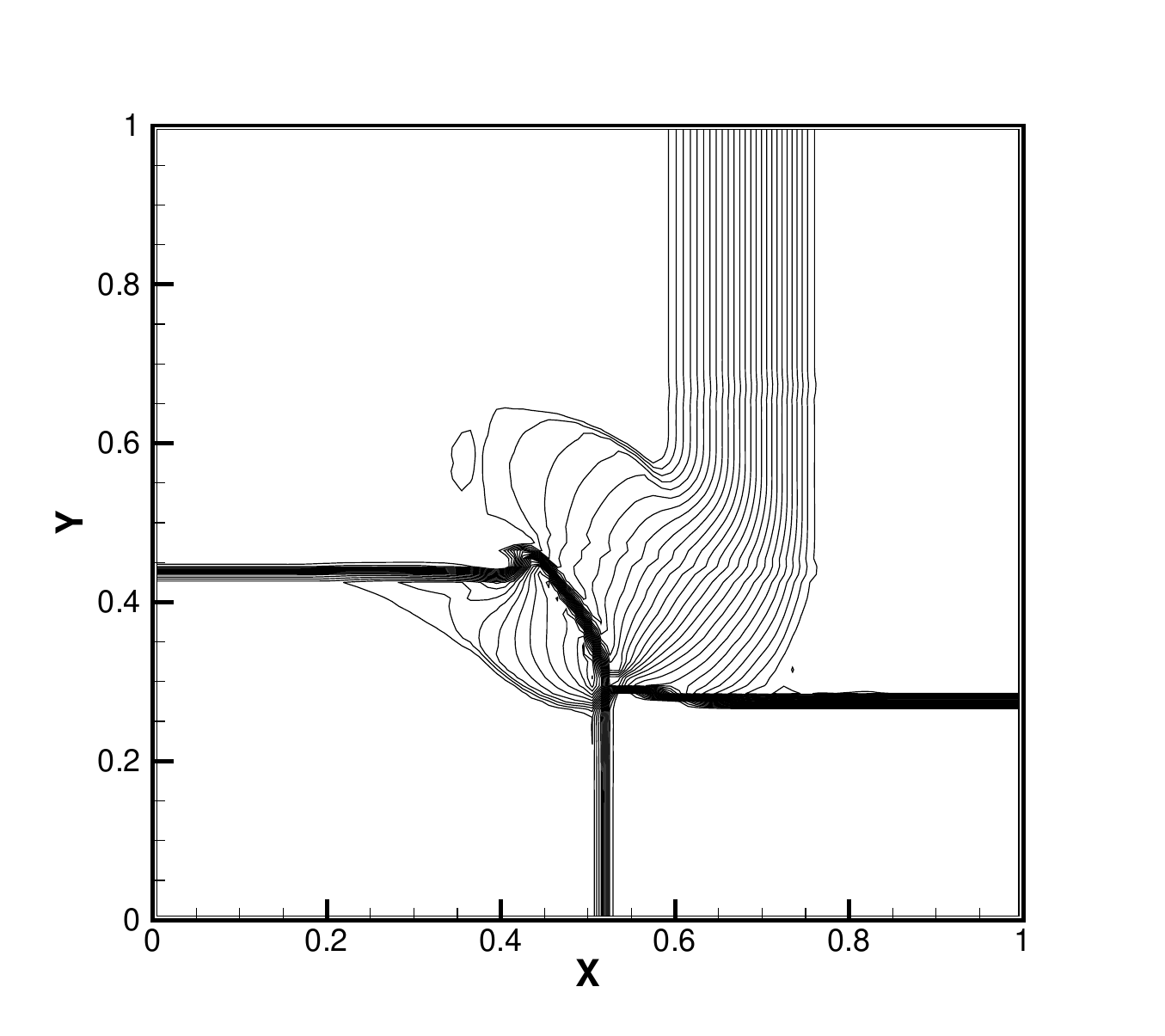}
    \includegraphics[width=0.48\textwidth]{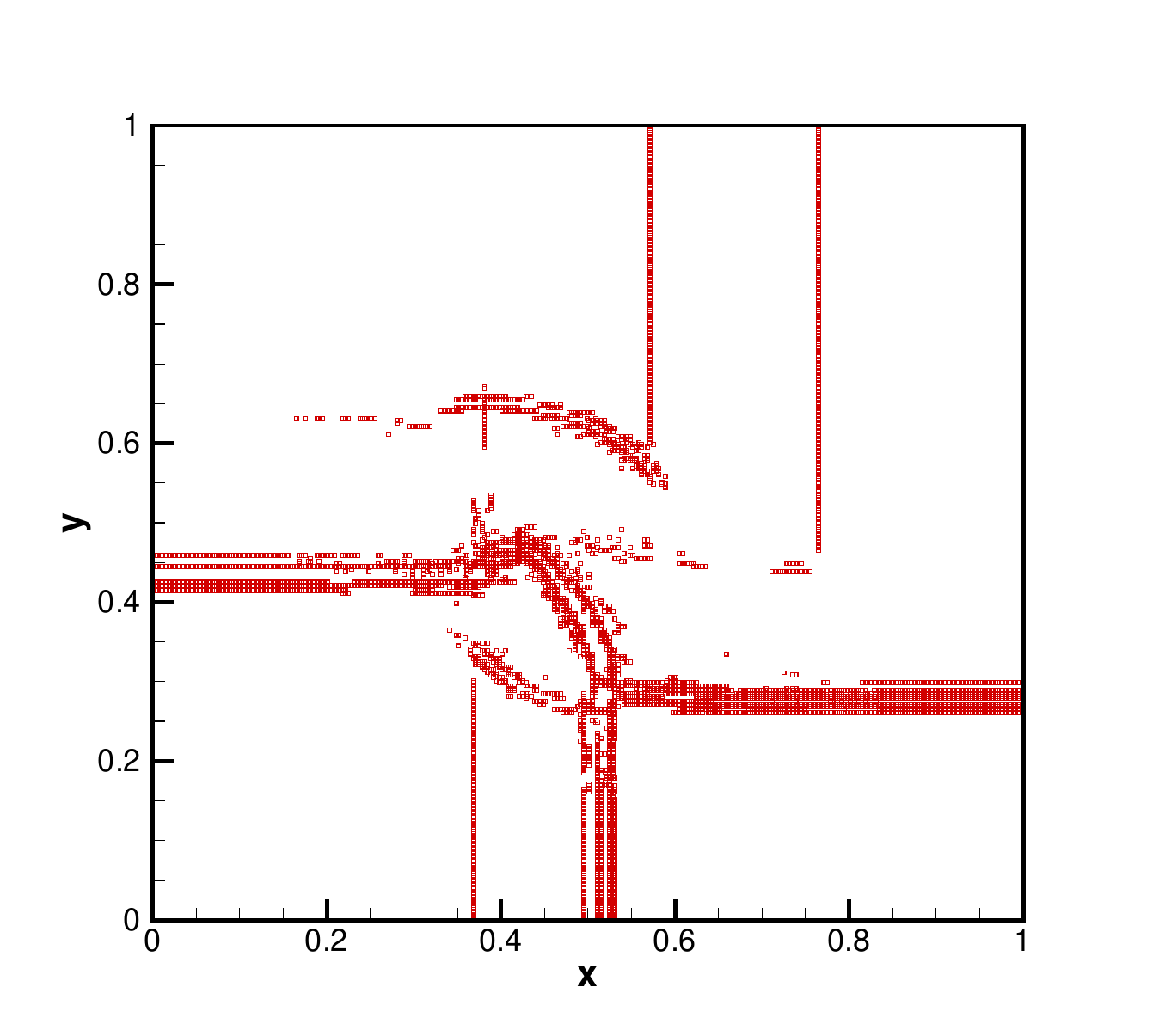}
    \includegraphics[width=0.48\textwidth]{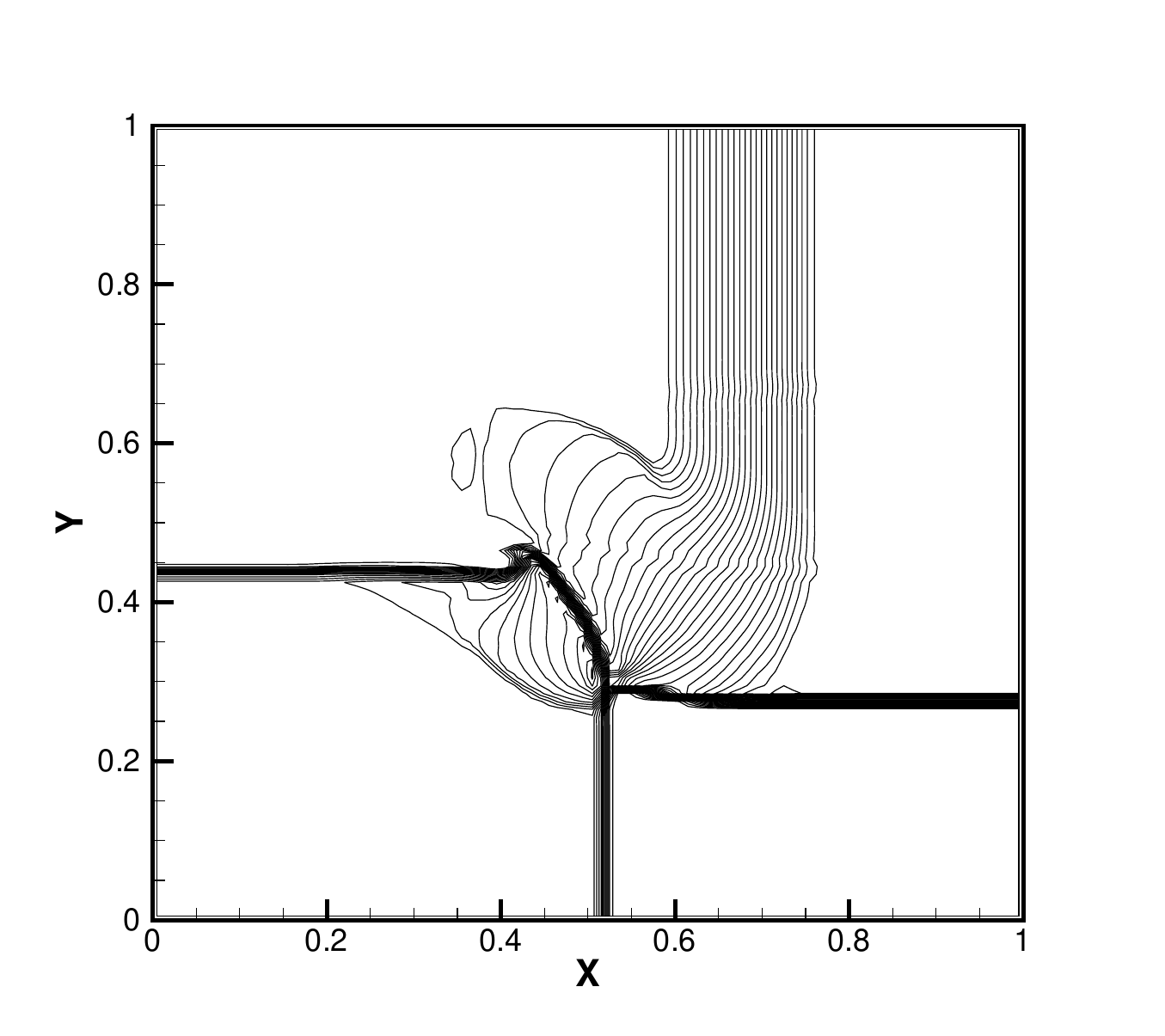}
    \includegraphics[width=0.48\textwidth]{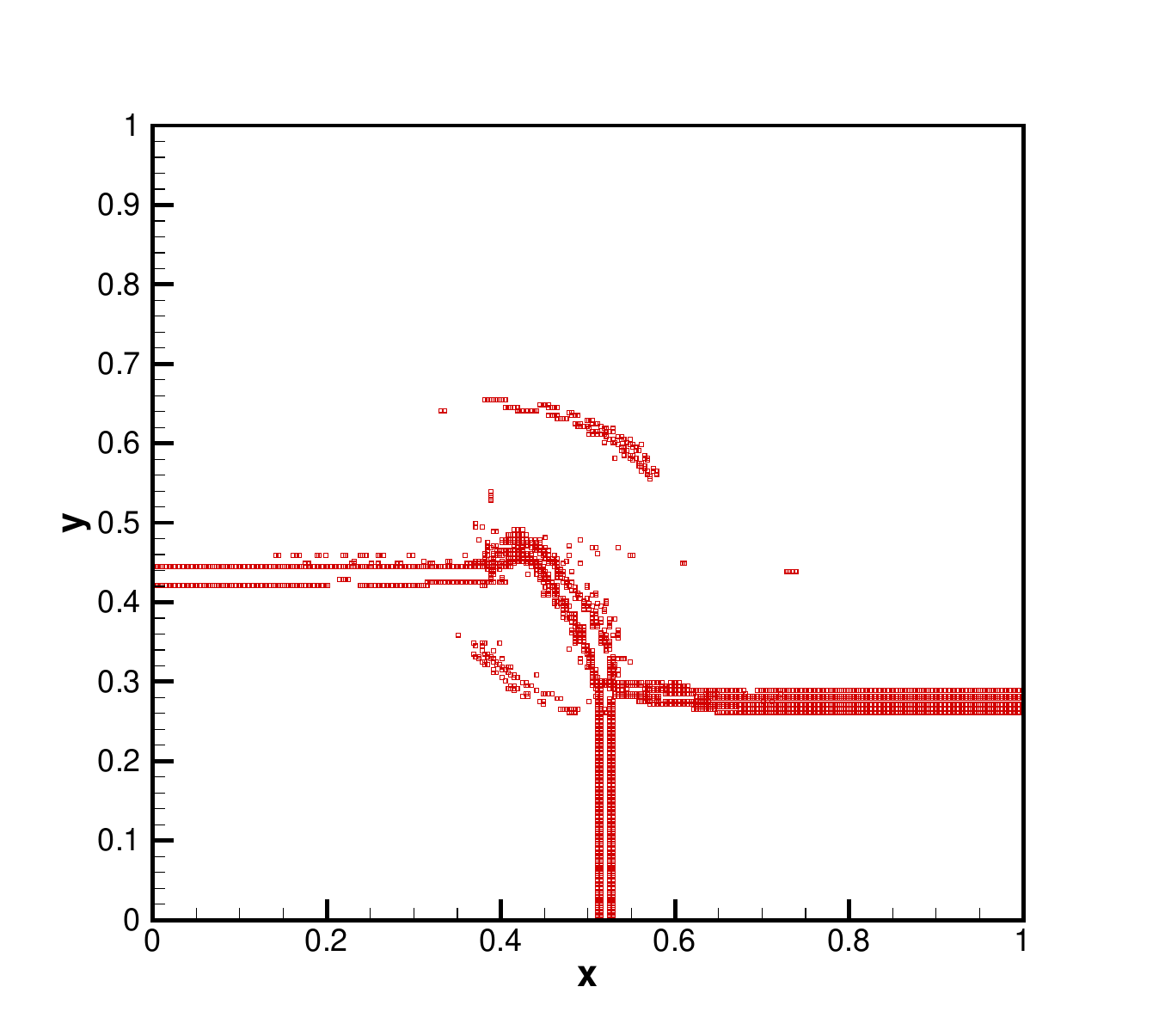}
     \caption{Example \ref{chap03:Riemann2d2}: 2D Riemannproblem \uppercase\expandafter{\romannumeral2}:
     30 equally density contour lines(left) from 0.55 to 1.06 and troubled cells(right) of \texttt{SV-cvMSWENO3} scheme when $M=0.01$(top),$M=100$(middle),$M=200$(bottom).}
    \label{Fig:riemann2d-2-k3}
 \end{figure}

 \begin{figure}[htbp]
    \centering
    \includegraphics[width=0.48\textwidth]{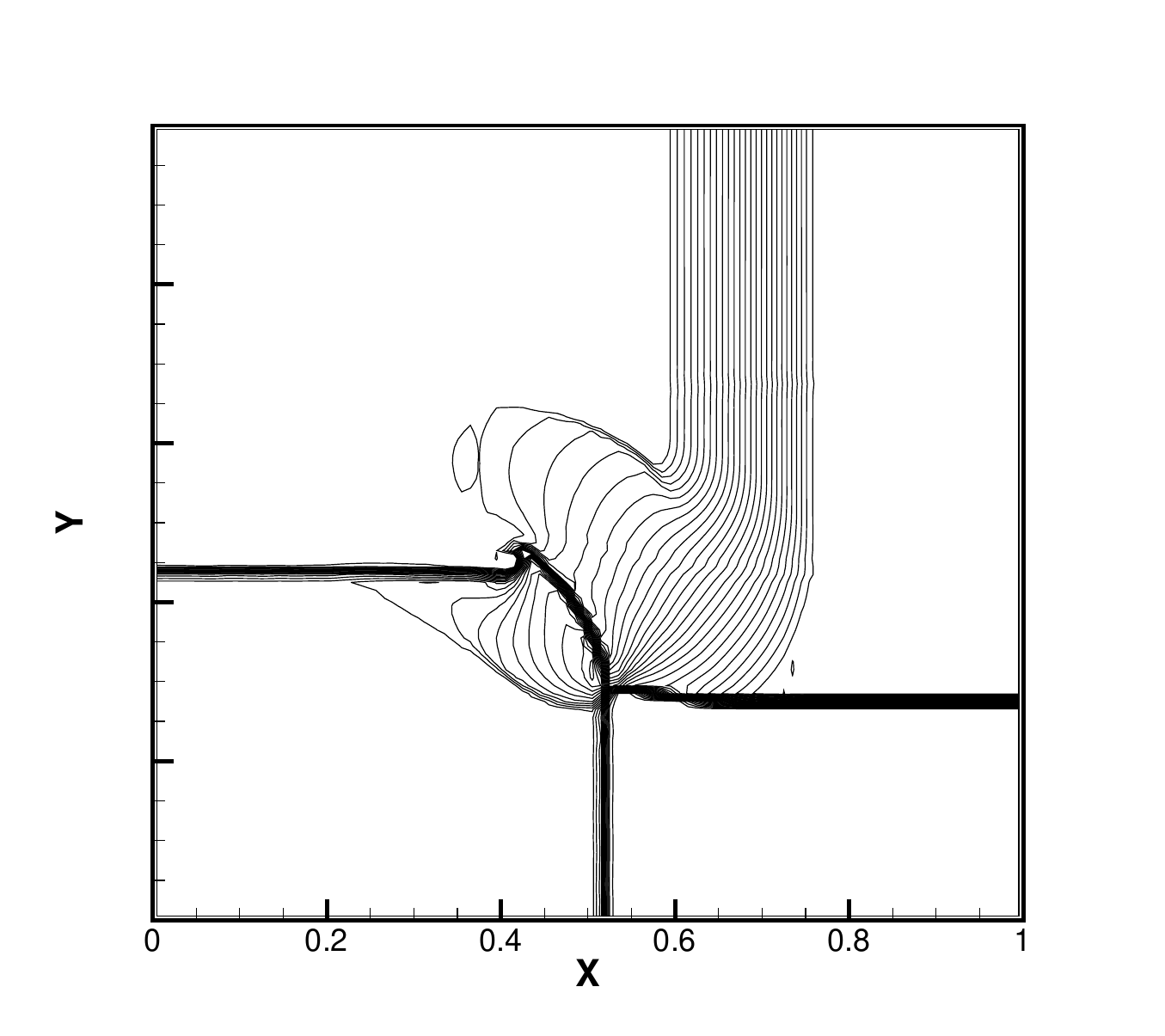}
    \includegraphics[width=0.48\textwidth]{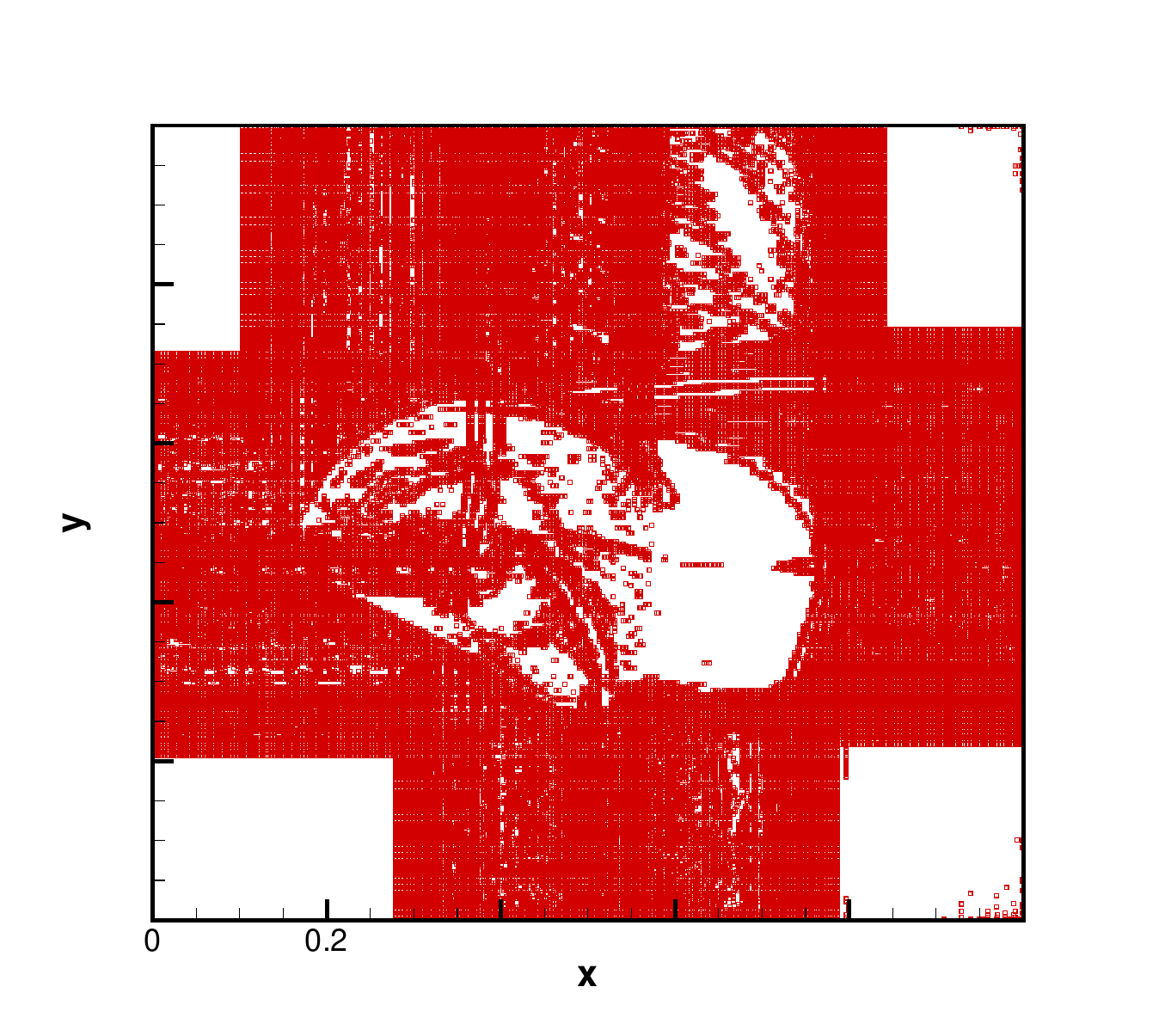}
    \includegraphics[width=0.48\textwidth]{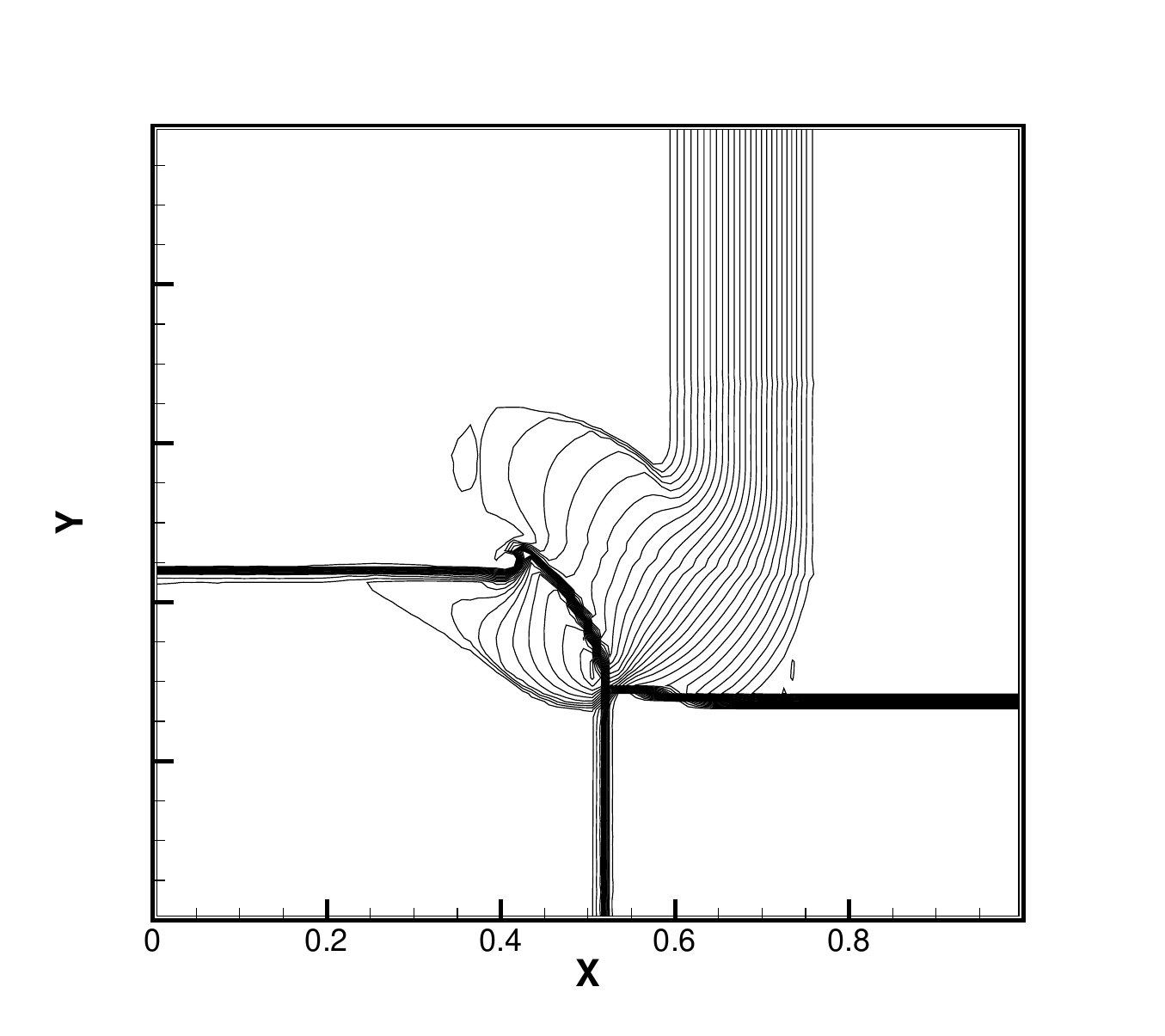}
    \includegraphics[width=0.48\textwidth]{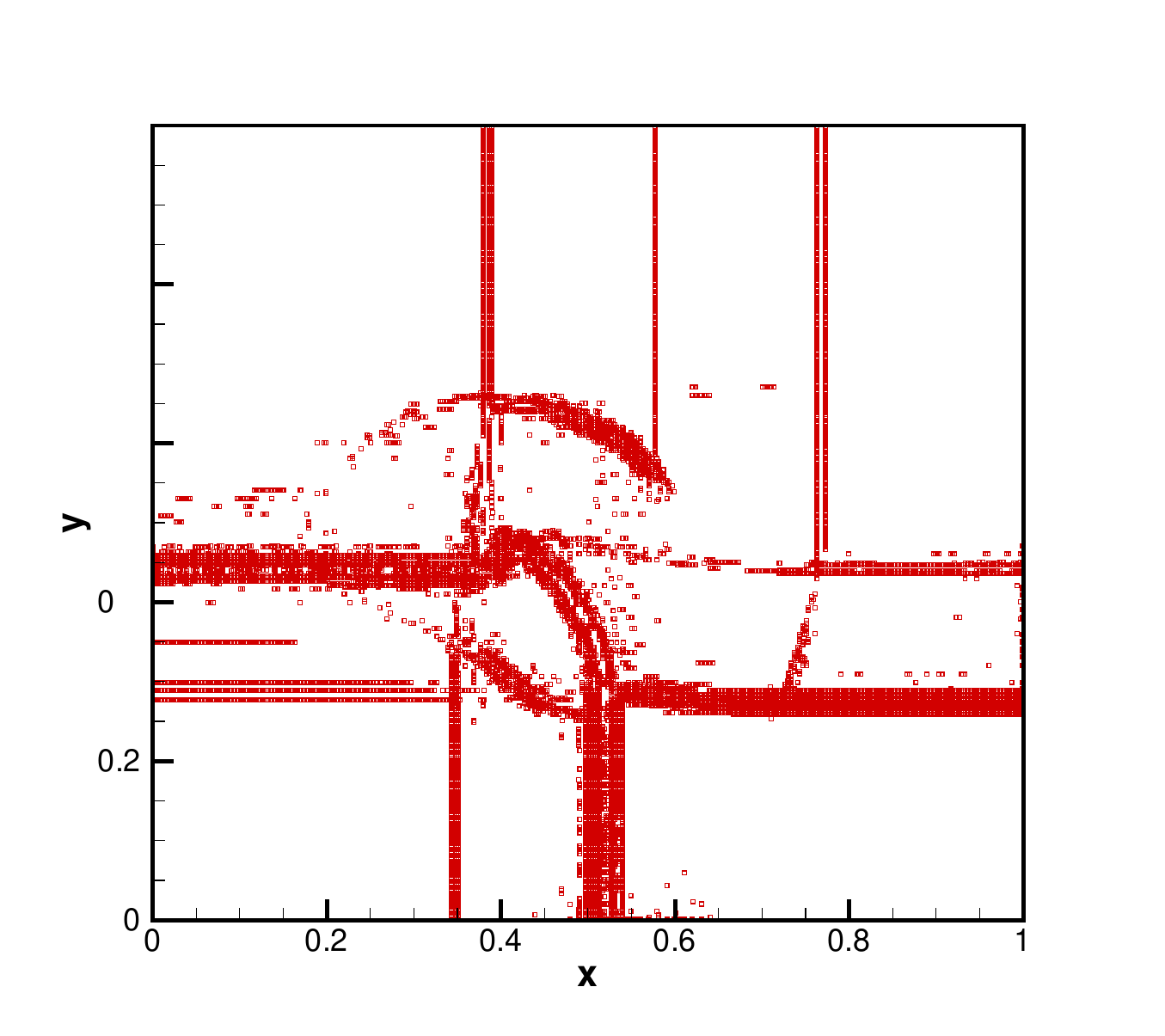}
    \includegraphics[width=0.48\textwidth]{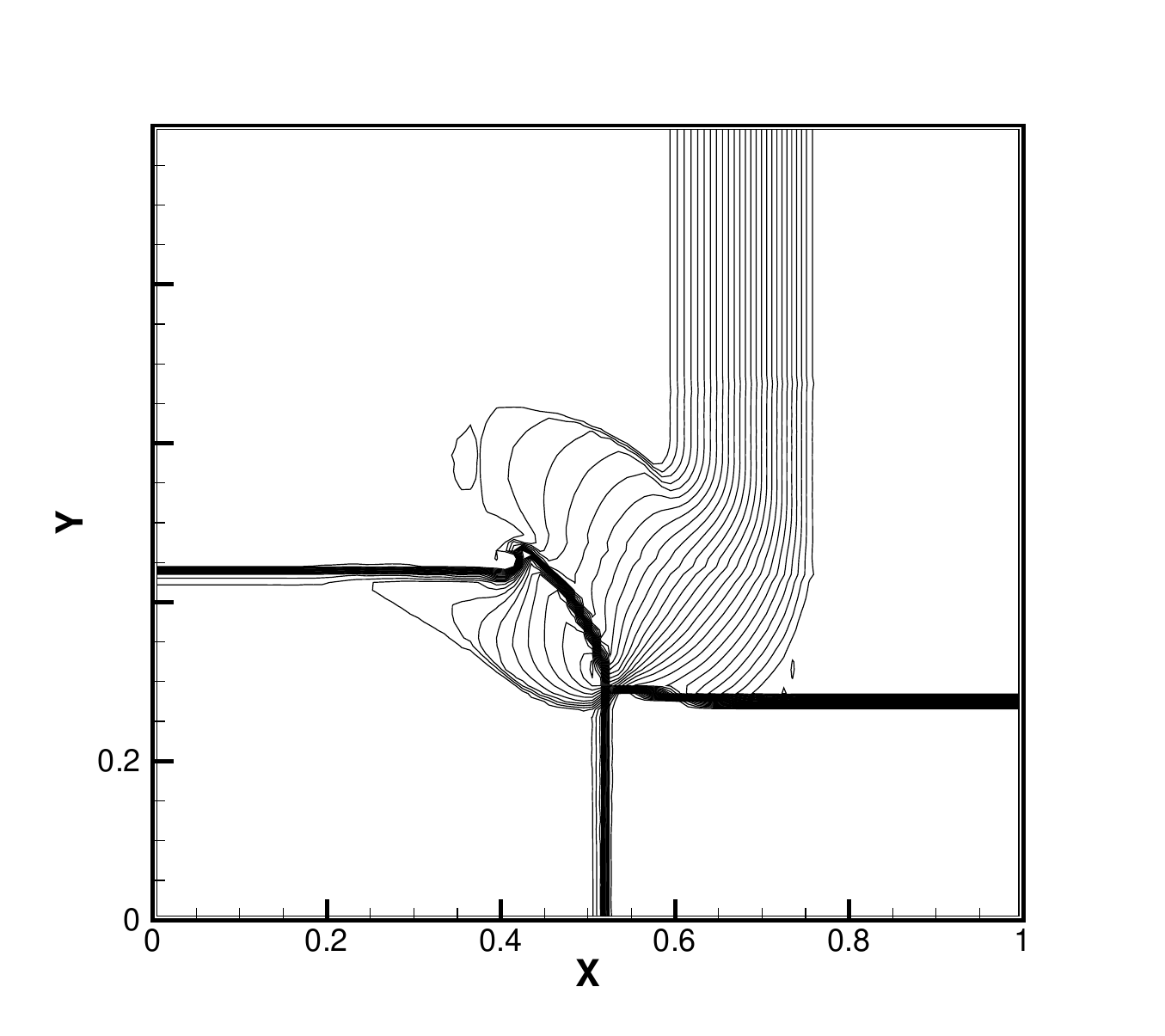}
    \includegraphics[width=0.48\textwidth]{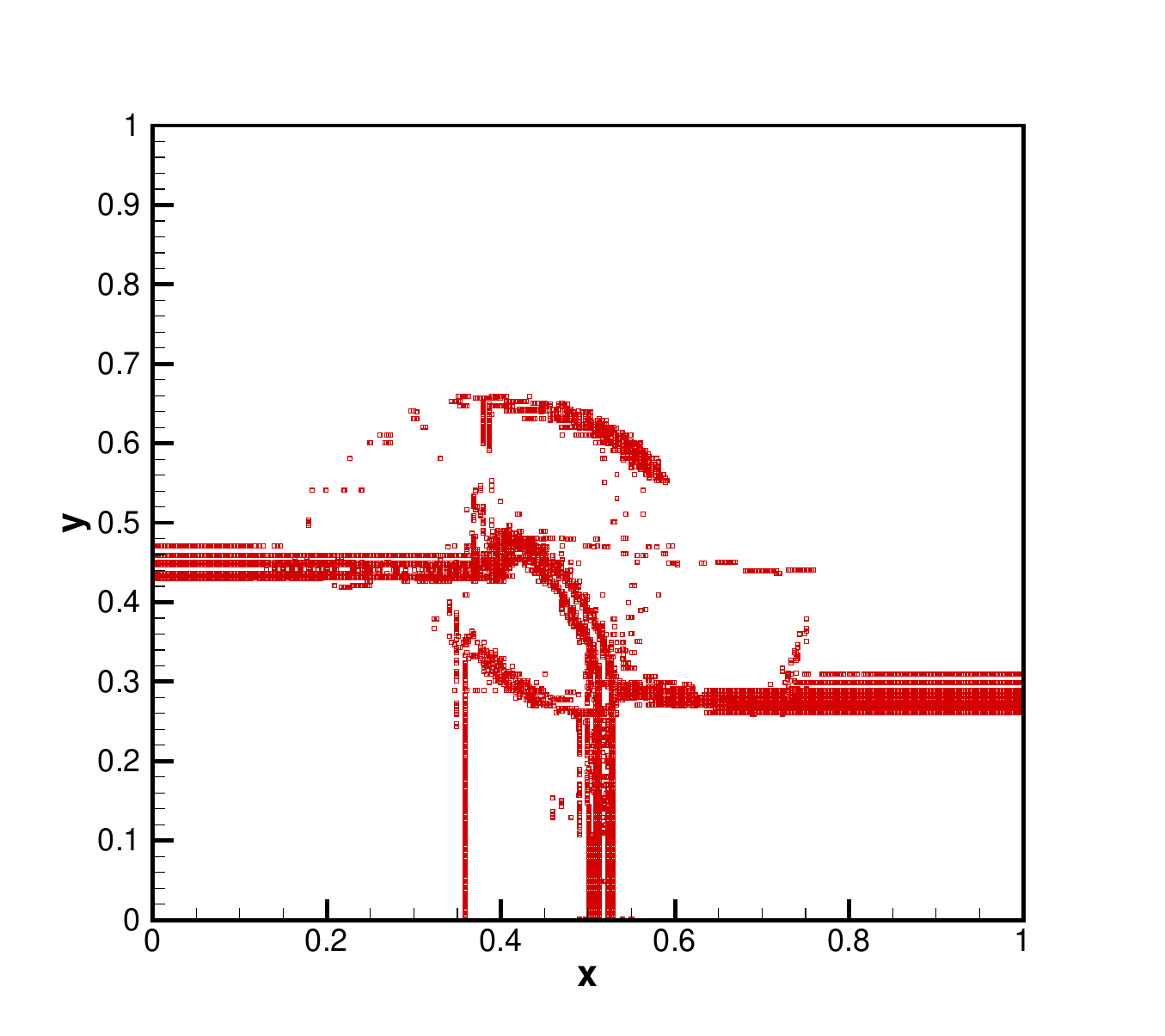}
     \caption{Example \ref{chap03:Riemann2d2}: 2D Riemannproblem \uppercase\expandafter{\romannumeral2}:
     30 equally density contour lines(left) from 0.55 to 1.06 and troubled cells(right) of \texttt{SV-cvMSWENO4} scheme when $M=0.01$(top),$M=100$(middle),$M=200$(bottom).}
    \label{Fig:riemann2d-2-k4}
 \end{figure}
\end{example}

\begin{example}[Double mach reflection problem]\label{chap03:doublemach}\rm
Initially, a right-lateral shock wave with a mach number of 10 is present at the point $(x,y)=(1/6, 0)$, forming at an angle of $60^\circ$ with the X-axis. Within the domain $\Omega=[0, 4]\times [0, 1]$, the initial conditions are defined as follows:
$$
\vec u(x,y,0)=\begin{cases}
\vec u_L, & y> h(x,0),\\
\vec u_R, &\ y<h(x,0),\end{cases}
$$
The $\vec u_L$ and $\vec u_R$ shock waves around the state are
$$
\vec u_L
=(8,57.1597,-33.0012,563.544)^T, \quad
\vec u_R
=(1.4,0,0,2.5)^T.
$$
The position of the shock wave at time t is given by $y=h(x,t)=\sqrt{3}(x-\frac{1}{6})-20t$. For the subsequent boundary conditions, when $x>1/6$, the reflection  boundary condition applies, whereas for $x<1/6$, the region is set to the state after the wave $\vec u_L$. The upper and left boundaries are defined as follows:
\begin{align*}
\vec u(0,y,t)=&\vec u_L,\quad \vec u(4,y,t)=\vec u_R, \quad 0\leq y\leq 1;\\
\vec u(x,1,t)=&\begin{cases}
\vec u_L, & 0\leq x< x_s,\\ \vec u_R, & x_s < x\leq  4,
\end{cases}
\end{align*}
where $x_s$ is the solution to the equation $h(x,t)=1$.

In Figure \ref{Fig:doublemach}, density contour figures derived from $960 \times 240$ spectral volumes at $t=0.2$ and $M=0.01,100$ and 200 are presented using  the third-order scheme, respectively. 29 contour lines with densities ranging from 1.3 to 23 are uniformly distributed. Figure \ref{Fig:doublemach2} illustrates a local magnification of the density contour lines within the corresponding dual mach region, highlighting the high resolution dual mach structure captured by the scheme. Figure \ref{Fig:doublemach3} depicts the distribution of troubled cells for $M=0.01, 100$ and 200. As $M$ increased from 0.01 to 200, there was a slight improvement in resolution, accompanied by a reduction in the number of troubled cells.

 \begin{figure}[htbp]
    \centering
    \includegraphics[width=0.8\textwidth,height=0.3\textwidth]{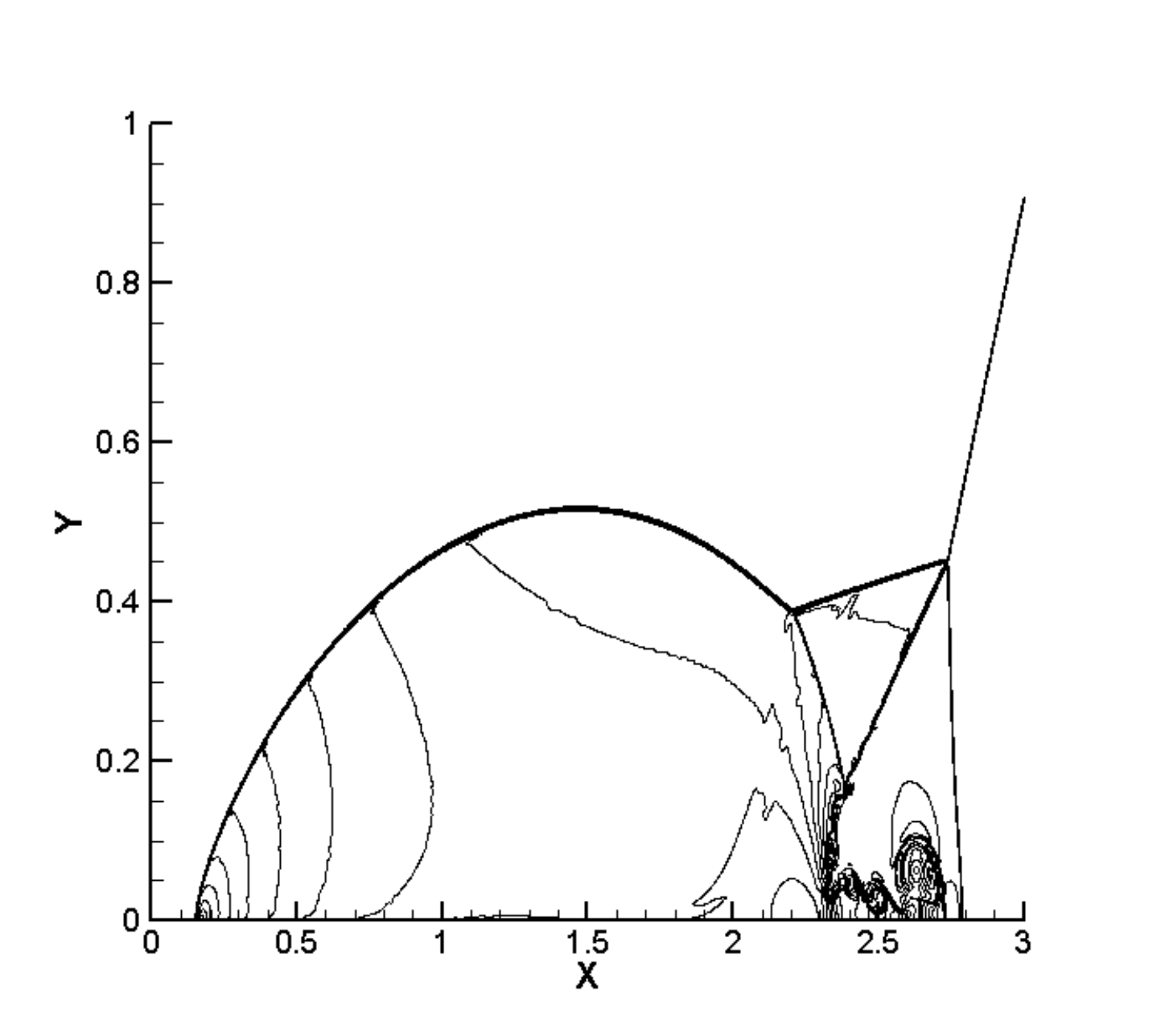}
    \includegraphics[width=0.8\textwidth,height=0.3\textwidth]{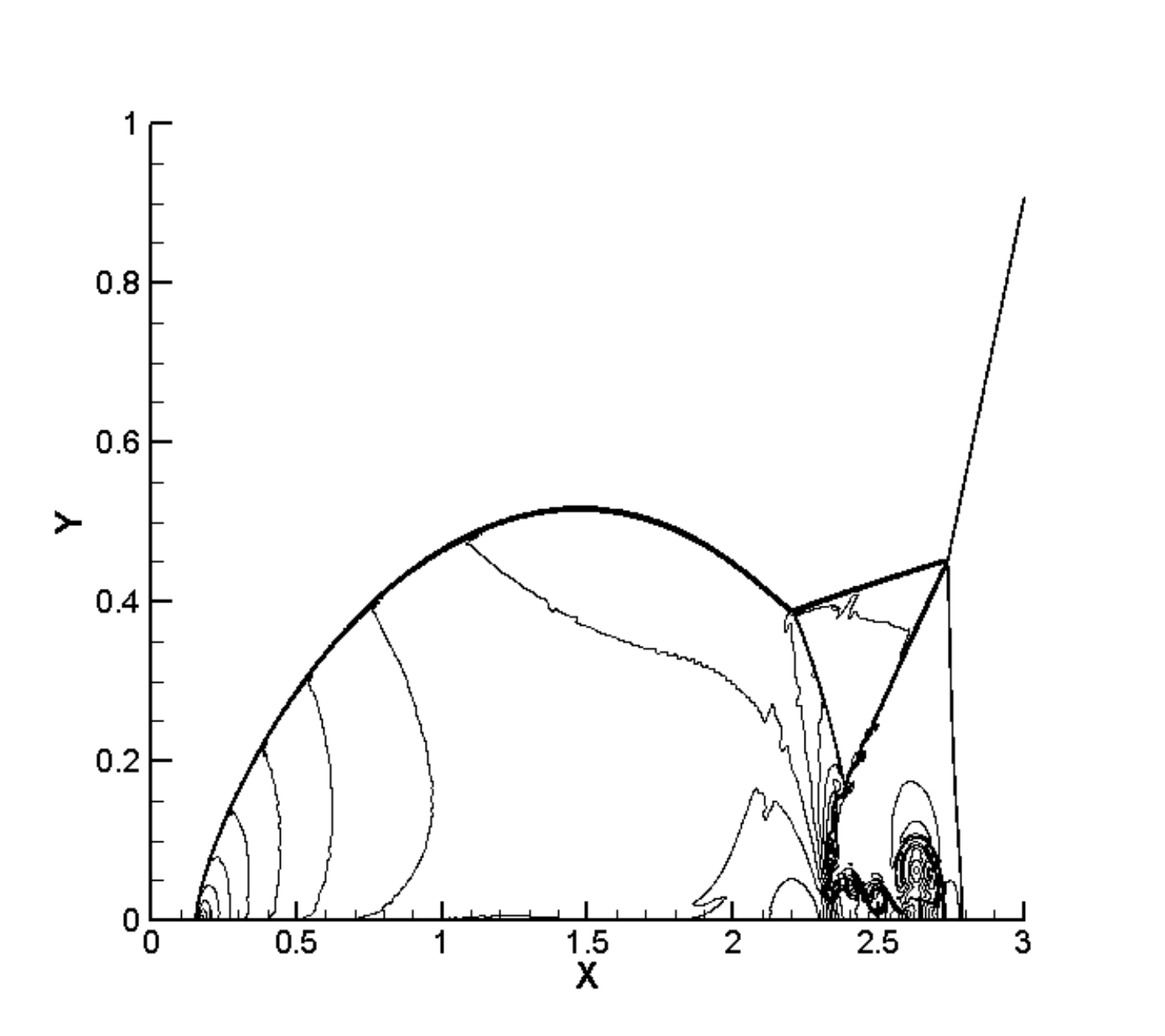}
    \includegraphics[width=0.8\textwidth,height=0.3\textwidth]{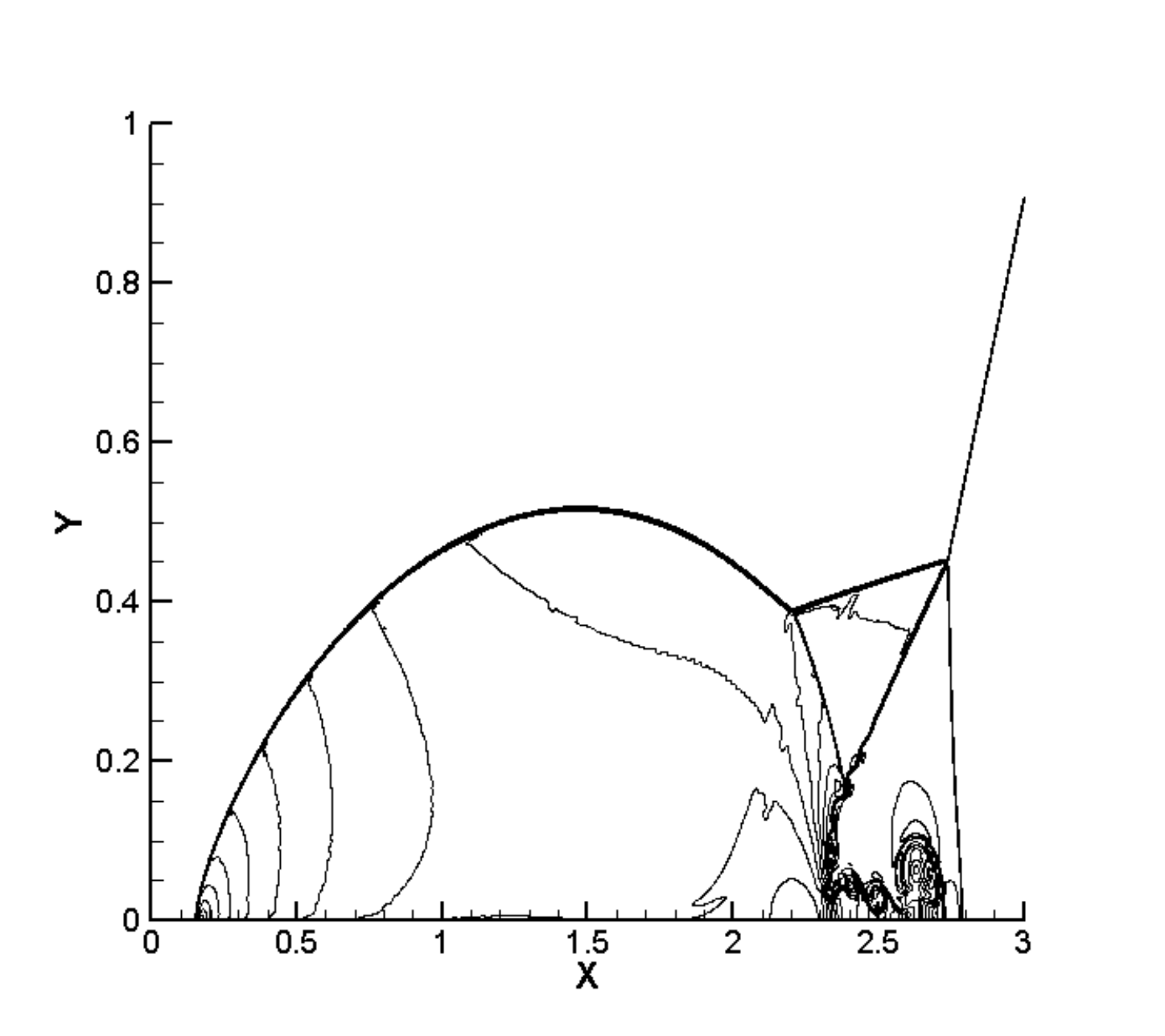}
     \caption{Example \ref{chap03:doublemach}: Double mach reflection problem: 29 equally density contour lines from 1.3 to 29 of \texttt{SV-cvMSWENO3} scheme when $M=0.01$(top),$M=100$(middle),$M=200$(bottom).}
    \label{Fig:doublemach}
 \end{figure}

  \begin{figure}[htbp]
    \centering
    \includegraphics[width=0.6\textwidth]{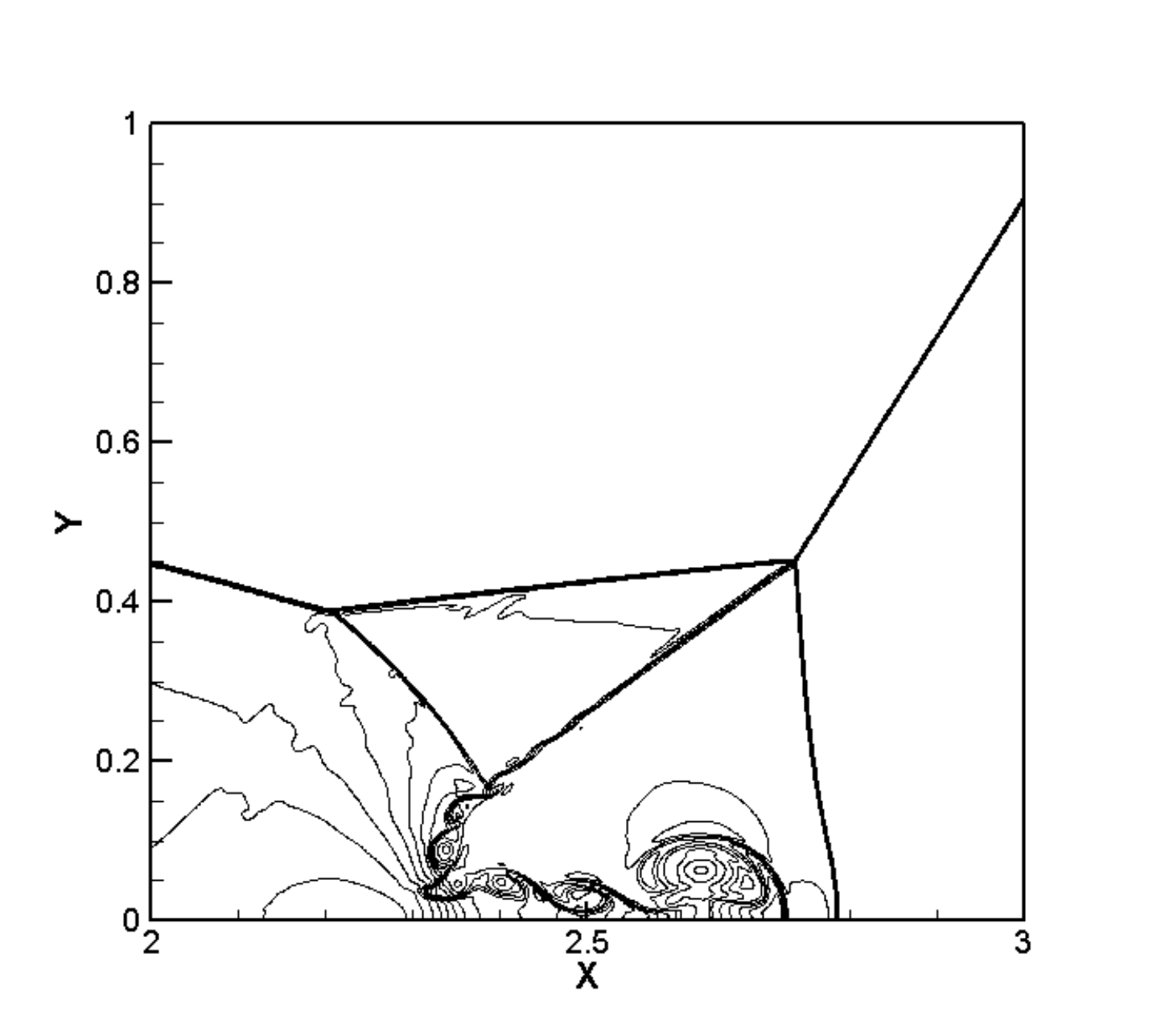}
    \includegraphics[width=0.6\textwidth]{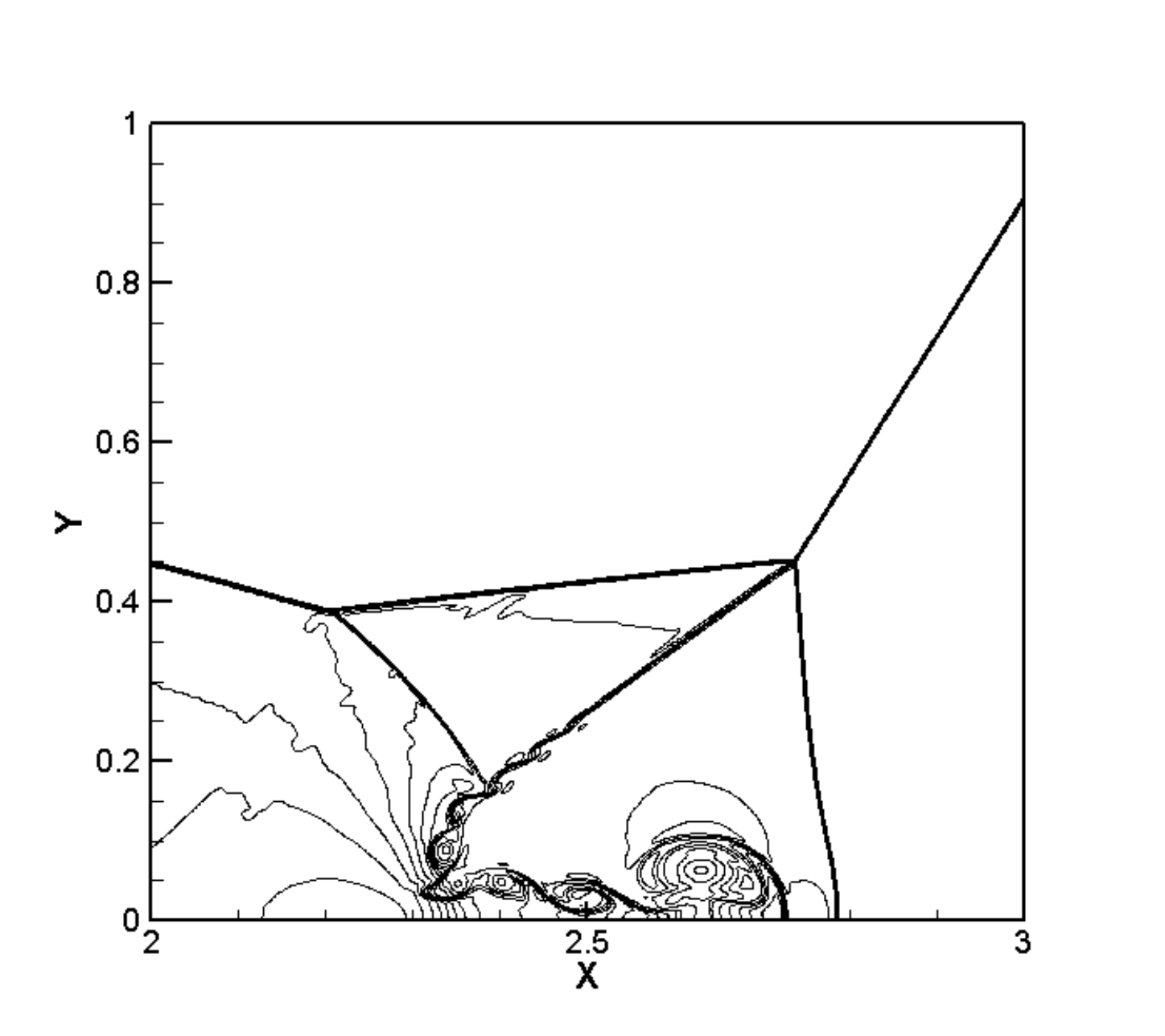}
    \includegraphics[width=0.6\textwidth]{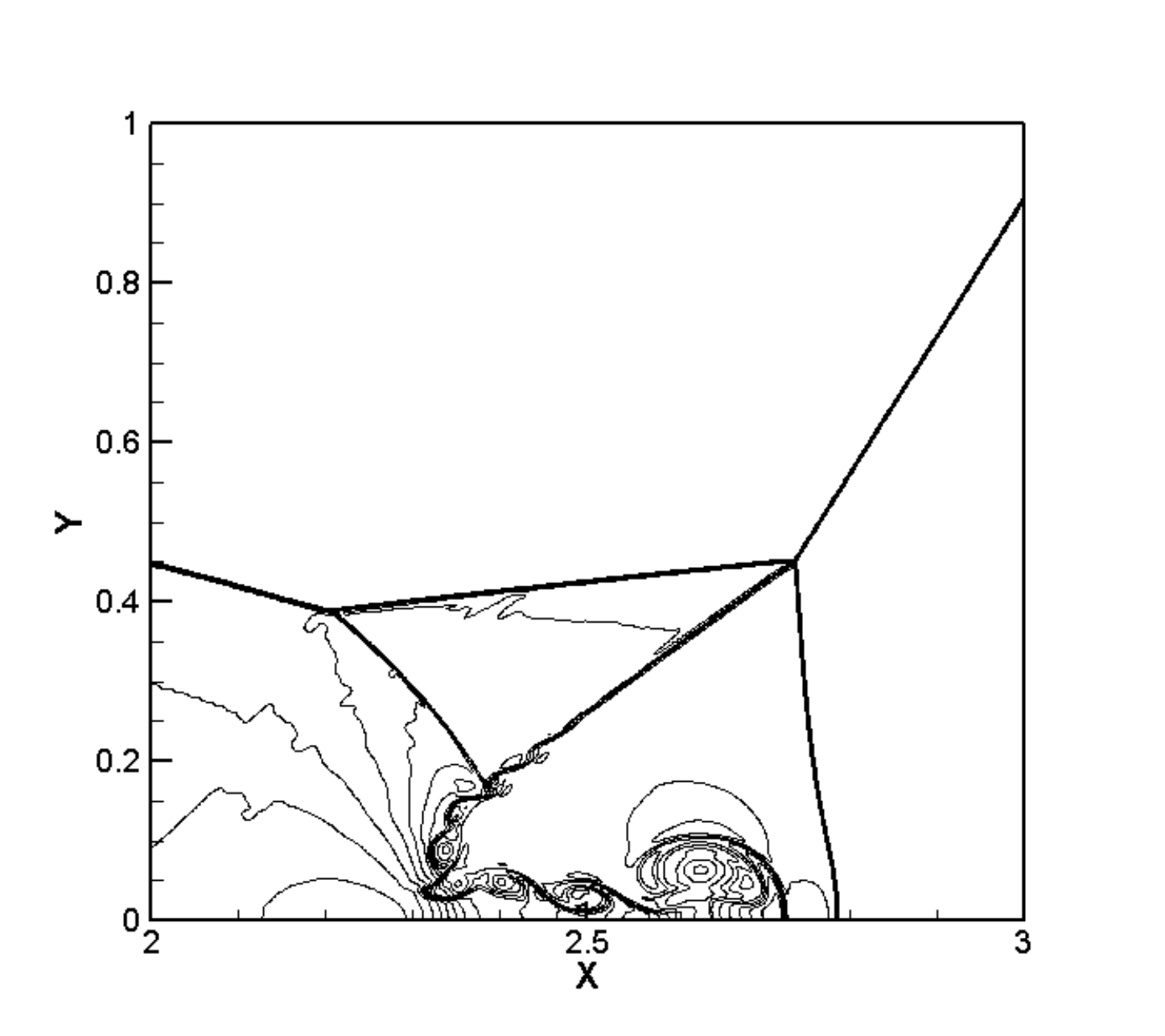}
     \caption{Example \ref{chap03:doublemach}: Double mach reflection problem: close-up view around the double Mach region of the density contour lines from 1.3 to 29 of \texttt{SV-cvMSWENO3} scheme when $M=0.01$(top),$M=100$(middle),$M=200$(bottom).}
    \label{Fig:doublemach2}
 \end{figure}

  \begin{figure}[htbp]
    \centering
    \includegraphics[width=0.8\textwidth,height=0.3\textwidth]{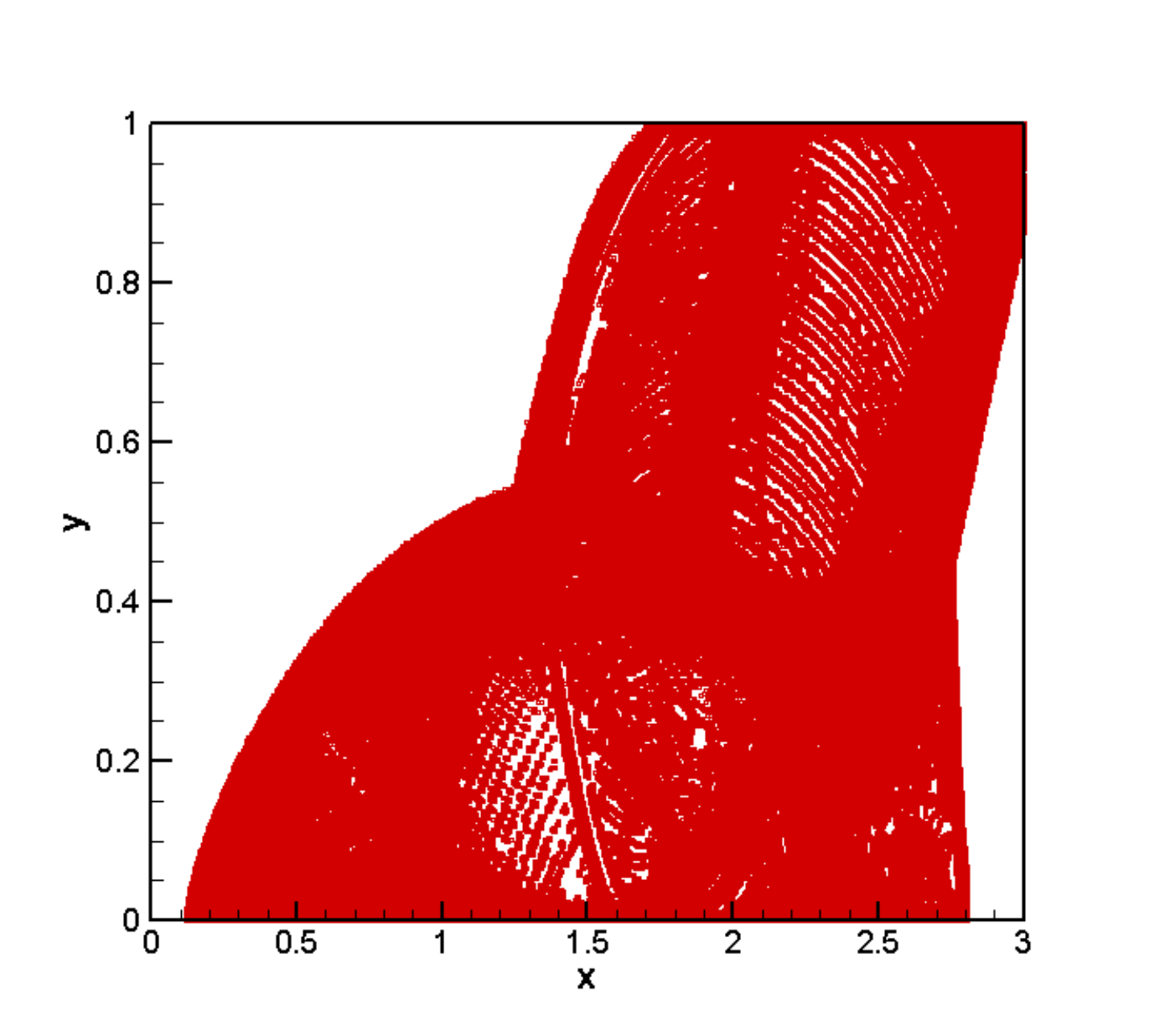}
    \includegraphics[width=0.8\textwidth,height=0.3\textwidth]{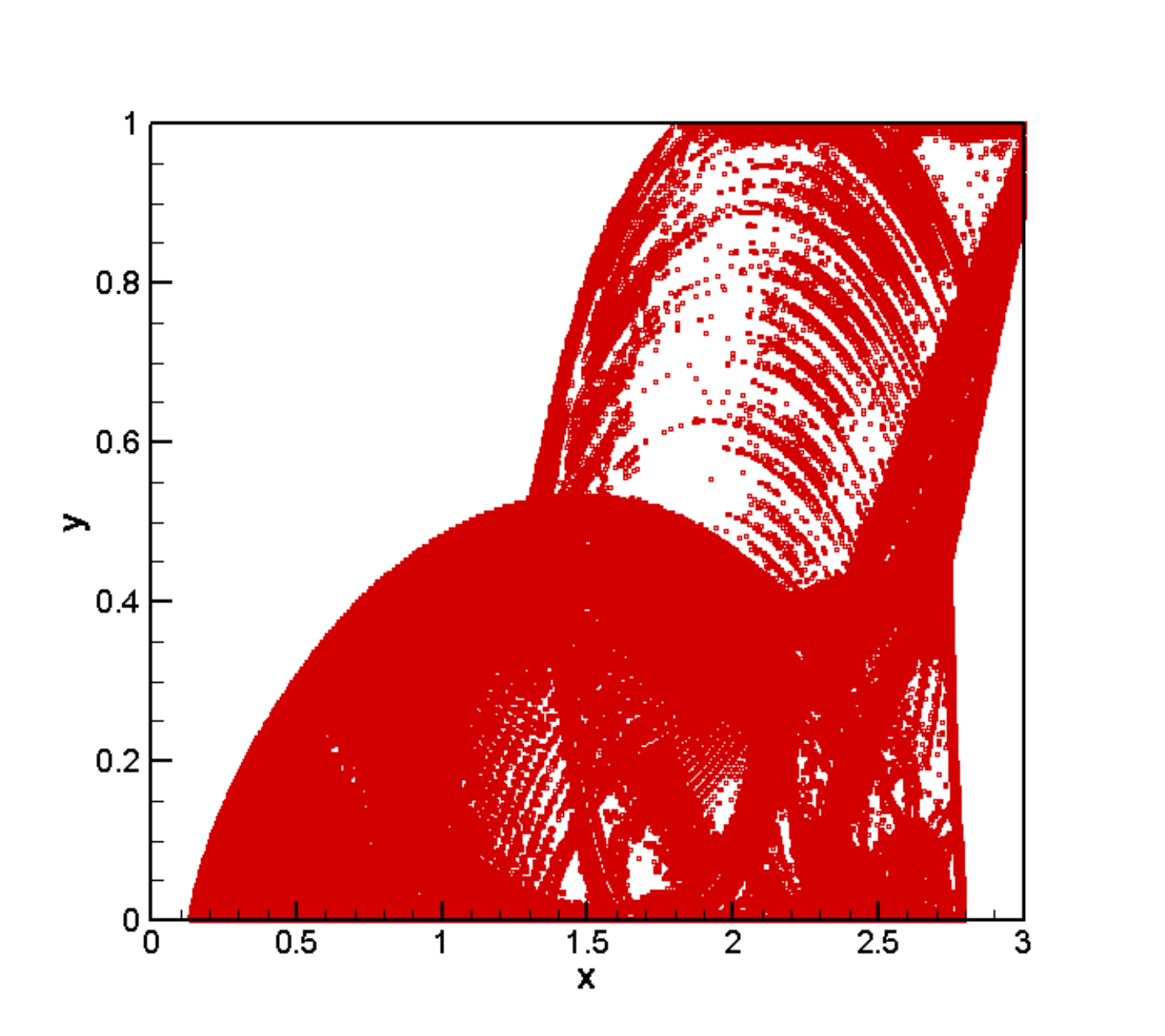}
    \includegraphics[width=0.8\textwidth,height=0.3\textwidth]{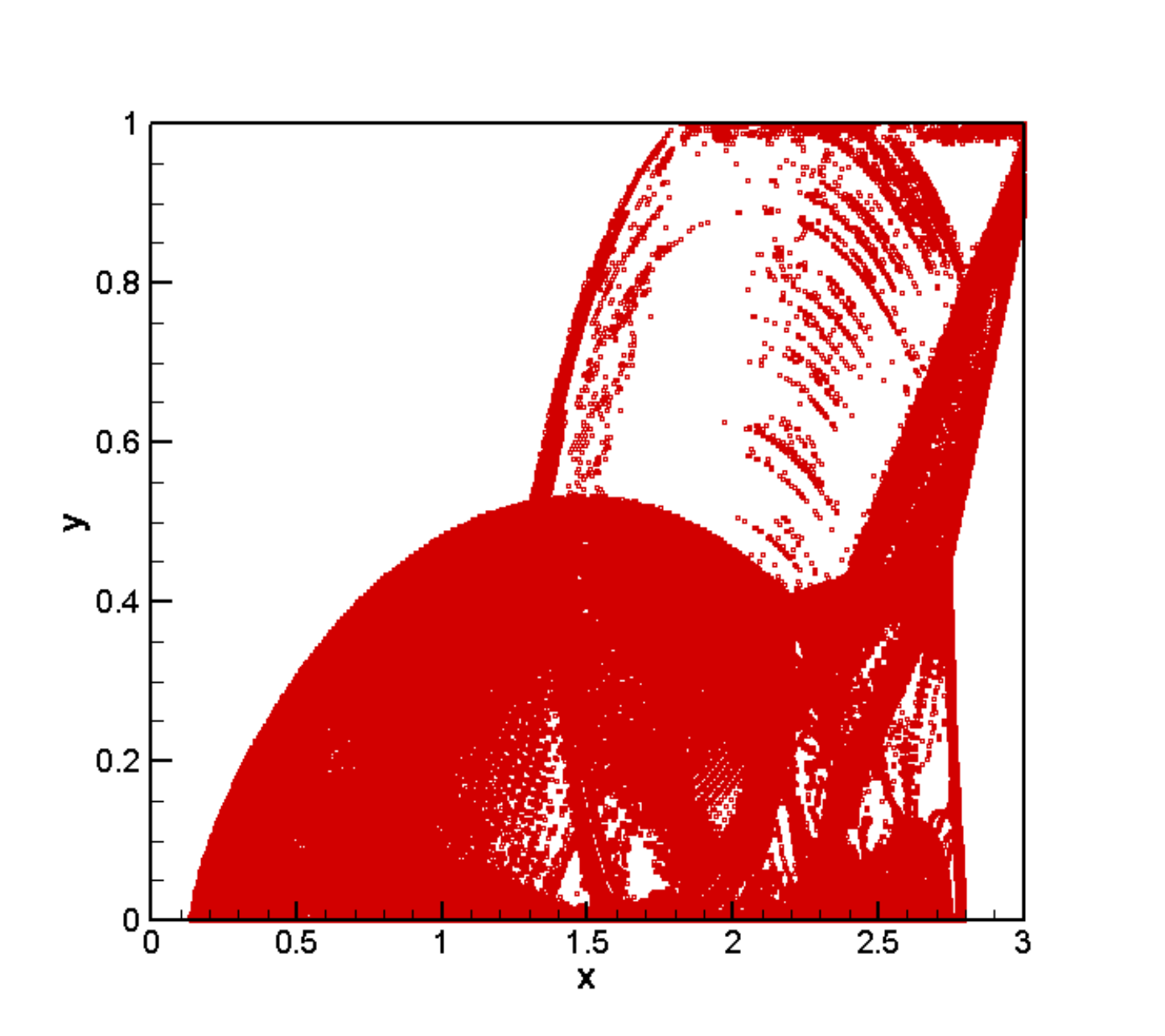}
     \caption{Example \ref{chap03:doublemach}: Double mach reflection problem: troubled cells of the \texttt{SV-cvMSWENO3} scheme when $M=0.01$(top),$M=100$(middle),$M=200$(bottom).}
    \label{Fig:doublemach3}
 \end{figure}

\end{example}

\section{Conclusion}
Based on a simple weighted non-oscillatory method, a high resolution limiter is developed for the spectral volume method to enhance the resolution of the control volume. This limiter preserves the resolution of the control volume while maintaining the compactness inherent in the spectral volume method. The proposed cvSWENO limiter is straightforward and easy to implement. Numerical results confirm the high resolution characteristics of the limiter. To our knowledge, this is the first instance of introducing a high resolution control volume limiter for the spectral volume method.

\section{Acknowledgements}
The work was partially supported by President Foundation of CAEP(Grant No.YZJJZQ2023024), Science Challenge Project
of China(Grant No.TZ2025007), and the National Natural Science Foundation of China (No.
11971070, 12471390).




\end{document}